\documentclass[10pt]{article} 
\def\date{20.5.2010}   
\usepackage{hyperref}
\usepackage{amssymb}
\usepackage{amsmath}
\input{liemacs10.sty} 

\newcommand{\sg}{{\mathfrak{sg}}}
\newcommand{\sk}{{\mathfrak{sk}}}

\newcommand{\Osc}{\mathop{\rm Osc{}}\nolimits}
 
\renewcommand{\L}{\mathop{\bf L{}}\nolimits}

\renewcommand{\mlabel}{\label}

\begin{document} 



\title{Semibounded representations of hermitian Lie groups} 

\author{Karl-Hermann Neeb\begin{footnote}{
Department  Mathematik, FAU Erlangen-N\"urnberg, Bismarckstrasse 1 1/2, 
91054-Erlangen, Germany; neeb@mi.uni-erlangen.de}
\end{footnote}
\begin{footnote}{Supported by DFG-grant NE 413/7-1, Schwerpunktprogramm 
``Darstellungstheorie''.} 
\end{footnote}}

\maketitle

\begin{abstract}  A unitary representation 
of a, possibly infinite dimensional, Lie group 
$G$ is called semibounded if the corresponding 
operators $i\dd\pi(x)$ from the derived representation  
are uniformly bounded from above on some non-empty open subset 
of the Lie algebra $\g$ of $G$. A hermitian Lie group is a central 
extension of the identity component of the 
automorphism group of a hermitian Hilbert symmetric space. 
In the present paper we classify the irreducible semibounded 
unitary representations of hermitian Lie groups corresponding 
to infinite dimensional irreducible symmetric spaces. 
These groups come in three essentially different types: those corresponding 
to negatively curved spaces (the symmetric Hilbert domains), 
the unitary groups acting on the duals of Hilbert domains, such as 
the restricted Gra\ss{}mannian, and the motion groups 
of flat spaces. \\
{\em Keywords:} infinite dimensional Lie group, unitary representation, 
semibounded representation, hermitian symmetric space, symmetric 
Hilbert domain. \\
{\em MSC2000:} 22E65, 22E45.  
\end{abstract} 

\section*{Introduction} 

This paper is part of a project concerned with a 
systematic approach to unitary representations 
of Banach--Lie groups in terms of conditions 
on spectra in the derived representation. For the derived 
representation to carry significant information, we have to 
impose a suitable smoothness condition. A unitary 
representation $\pi \: G \to \U(\cH)$ 
is said to be {\it smooth} if the subspace 
$\cH^\infty \subeq \cH$ of smooth vectors is dense. 
This is automatic for continuous 
representations of finite dimensional groups, but not 
for Banach--Lie groups (\cite{Ne10a}). For any smooth 
unitary representation, the {\it derived representation} 
\[ \dd\pi \: \g = \L(G)\to \End(\cH^\infty), \quad 
\dd\pi(x)v := \derat0 \pi(\exp tx)v\]  
carries significant information in the sense that the closure of the 
operator $\dd\pi(x)$ coincides with the infinitesimal generator of the 
unitary one-parameter group $\pi(\exp tx)$. We call $(\pi, \cH)$ 
{\it semibounded} if the function 
\[ s_\pi \: \g \to \R \cup \{ \infty\}, \quad 
s_\pi(x) 
:= \sup\big(\Spec(i\dd\pi(x))\big) \]
is bounded on the neighborhood of some point in $\g$.  
Then the set $W_\pi$ of all such points 
is an open $\Ad(G)$-invariant convex cone in the Lie algebra $\g$. 
We call $\pi$ {\it bounded} if $s_\pi$ is bounded on some $0$-neighborhood, 
i.e., $W_\pi = \g$. 
All finite dimensional continuous unitary representations are bounded 
and most of the unitary representations appearing in physics are semibounded 
(cf.\ \cite{Ca83}, \cite{Mi87, Mi89}, \cite{PS86}, 
\cite{SeG81}, \cite{CR87}, \cite{Se58}, \cite{Se78}, 
\cite{Bak07}).  

For finite dimensional Lie groups, the irreducible semibounded 
representations are precisely the unitary highest weight representations 
and one has unique direct integral decompositions \cite[X.3/4, XI.6]{Ne00}. 
For many other classes of groups such as 
the Virasoro group and affine Kac--Moody groups 
(double extensions of loop groups with compact target groups),  
the irreducible highest weight representations 
are semibounded, but to prove the converse 
is more difficult and requires a thorough understanding 
of invariant cones in the corresponding Lie algebras as well 
as of convexity properties of coadjoint orbits (\cite{Ne10c}). 

Finite dimensional groups only have faithful bounded representations 
if their Lie algebras are compact, which is equivalent to the 
existence of an $\Ad$-invariant norm on the Lie algebra. For 
infinite dimensional groups, the picture is much more colorful. 
There are many interesting bounded representations, in particular 
of unitary groups of $C^*$-algebras (cf.\ \cite{BN11}) and a central 
result of Pickrell (\cite{Pi88}), combinded with classification results 
of Kirillov, Olshanski and I.~Segal  
(\cite{Ol78}, \cite{Ki73}, \cite{Se57}), 
implies that all separable unitary representations 
of the unitary group $\U(\cH)$ of an infinite dimensional separable 
Hilbert space can can be classified 
in the same way as for their finite dimensional analogs by 
Schur--Weyl theory. In particular,  the irreducible ones are bounded 
(cf.\ Theorem~\ref{thm:8.1b}). 

To address classification problems one needs refined 
analytic tools based on recent results asserting that the space
$\cH^\infty$ of smooth vectors is a Fr\'echet space on which 
$G$ acts smoothly (\cite{Ne10a}). 
In \cite{Ne10d} we use these facts to develop some spectral 
theoretic tools concerning the space of smooth vectors. 
This does not only lead to a 
complete description of semibounded representations 
of various interesting classes of groups such as hermitian Lie groups 
which are  dealt with in the present paper. 
They also apply naturally to unitary representations 
of Lie supergroups generated by their odd part 
(\cite{NSa10}). 

Our goal is a classification of the irreducible semibounded 
representations and the development of 
tools to obtain direct integral decompositions 
of semibounded representations. The first part of this 
goal is achieved in the present paper for the class of hermitian Lie groups. 
These are triples $(G,\theta, d)$, where 
$G$ is a connected Lie group, $\theta$ an involutive 
automorphism of $G$ with the corresponding eigenspace 
decomposition $\g = \fk \oplus \fp$, 
$d \in \z(\fk)$ (the center of $\fk$) 
an element for which $\ad d\res_\fp$ is a complex 
structure, and $\fp$ carries an $e^{\ad \fk}$-invariant 
Hilbert space structure. 
We then write $K := (G^\theta)_0$ for the identity component of 
the group of $\theta$-fixed points in $G$ and observe that our 
assumptions imply that $G/K$ is a hermitian Hilbert 
symmetric space. 
The simply connected symmetric spaces arising from this construction 
have been classified by W.~Kaup by observing that 
$\fp$ carries a natural structure of a $JH^*$-triple, 
and these objects permit a powerful structure theory which 
leads to a complete classification in terms of orthogonal 
decompositions and simple objects (cf.\ \cite{Ka81,Ka83}). 

Typical examples of hermitian symmetric spaces
are symmetric Hilbert domains (the negatively curved case), 
their duals, such as the restricted Gra\ss{}mannian of a polarized 
Hilbert space $\cH = \cH_+ \oplus \cH_-$ (the positively curved case) 
(cf.\ \cite{PS86}), 
and all Hilbert spaces $\cH$ as quotients 
$G = (\cH \rtimes_\alpha K)/K$, where $\alpha$ is a norm-continuous 
unitary representation of $K$ on $\cH$ (the flat case). 
Our concept of a hermitian Lie group contains almost 
no restrictions on the group $K$, but the structure of $K$ is crucial 
for the classification of semibounded representations of~$G$. 
To make this more specific, we call $(G,\theta,d)$ {\it irreducible} 
if the unitary $K$-representation on $\fp$ is irreducible. 
In this case either $[\fp,[\fp,\fp]] = \{0\}$ 
(the flat case) or $\fp$ is a simple $JH^*$-triple. 
Then we say that $\g$ is {\it full} if 
$\ad \fk = \aut(\fp)$ is the full Lie algebra 
of  the automorphism group $\Aut(\fp)$ of the $JH^*$-triple~$\fp$. 
One of our key observations is that 
the quotient Lie algebra $\fk/\z(\fk)$ contains no non-trivial 
open convex invariant cones.

The structure of this paper is as follows.
Section~\ref{sec:herm} introduces the concept of a 
hermitian Banach--Lie group and in Section~\ref{sec:2} 
we explain their connection with $JH^*$-triples and recall 
Kaup's classification of infinite dimensional simple $JH^*$-triples. 
As we shall see in the process, hermitian Lie groups $G$ have 
natural central extensions $\hat G$ and these central extensions 
often enjoy a substantially richer supply of semibounded unitary 
representations than the original group~$G$. This phenomenon 
is also well-known for the group of diffeomorphisms of 
the circle (cf.\ \cite{Ne10c}) and loop groups (\cite{PS86}, 
\cite{Ne01b}). This motivates our detailed discussion 
of central extensions in Section~\ref{sec:centext}. 
Any semibounded representation $(\pi, \cH)$ 
defines the open convex invariant cone $W_\pi \subeq \g$. 
Therefore the understanding of semibounded representations 
requires some information on the geometry of open invariant 
cones in Lie algebras. In our context we mainly need the information 
that for certain Lie algebras $\fk$, all open invariant cones 
in $\fk/\z(\fk)$ are trivial. Typical examples 
with this property 
are the Lie algebras of the unitary groups of real, complex or 
quaternionic Hilbert spaces (Section~\ref{sec:4}). 
Sections~\ref{sec:5}-\ref{sec:8} are devoted to the 
classification of the irreducible semibounded representations 
of hermitian Lie groups. 
The main steps in this classification are the following results 
on semibounded representations $(\pi, \cH)$ of $G$: 
\begin{description}
\item[\rm(1)] If $\fk/\z(\fk)$ contains no open invariant cones, 
then we derive from the results on invariant cones developed in 
Section~\ref{sec:4} that $\pi\res_{Z(K)_0}$ is also semibounded. 
\item[\rm(2)] If $(G,\theta, d)$ is irreducible, 
then $d \in W_\pi \cup - W_\pi$. If $d \in W_\pi$, then we call 
$(\pi, \cH)$ a {\it positive energy representation}. In this 
case the maximal spectral value of the essentially selfadjoint operator 
$i\dd\pi(d)$ is an eigenvalue and the $K$-representation 
$(\rho,V)$ 
on the corresponding eigenspace is bounded and 
irreducible. Using the holomorphic induction techniques developed 
in \cite{Ne10d} for Banach--Lie groups, it follows that 
$(\pi, \cH)$ is uniquely determined by $(\rho,V)$ (Section~\ref{sec:5}). 
We call a bounded representation $(\rho,V)$ of $K$ 
{\it (holomorphically) inducible} 
if it corresponds as above to a unitary representation $(\pi, \cH)$ 
of~$G$. 
\item[\rm(3)] If $(G,\theta,d)$ is full or $G/K$ is flat 
(with $\fk/\fz(\fk)$ not containing open invariant cones), 
we derive an explicit 
characterization of the inducible bounded irreducible $K$-representations. 
We thus obtain a classification of all irreducible semibounded 
representations of $G$ in terms of the corresponding 
$K$-representations (Sections~\ref{sec:6}-\ref{sec:8}). 
For irreducible symmetric Hilbert domains, 
this explicit characterization is based on  the classification of 
unitarizable highest weight modules of 
locally finite hermitian Lie algebras obtained in \cite{NO98} 
(Section~\ref{sec:7}). 
For the $c$-dual spaces we obtain the surprisingly simple result 
that inducibility can be characterized by 
an easily verifiable positivity condition 
(Section~\ref{sec:8}).  
\item[\rm(4)] A central point in our characterization 
is that the group $K$ has a normal closed subgroup 
$K_\infty$ (for $K = \U(\cH)$ this is the group $\U_\infty(\cH)$) 
with the property that each bounded irreducible 
representation $(\rho,V)$ is a tensor product 
of two bounded irreducible 
representations $(\rho_0, V_0)$ and $(\rho_1, V_1)$, 
where $\rho_0\res_{K_\infty}$ is irreducible and 
$K_\infty \subeq \ker \rho_1$. We show that 
$(\rho,V)$ is inducible if and only if $(\rho_0, V_0)$ has this property.
Here the main point is that, even though 
no classification of the representations $(\rho_1, V_1)$ is known, 
the representations $(\rho_0, V_0)$ 
can be parameterized easily by ``highest weights'' 
(cf.\ Definition~\ref{def:tensrep}).
If $V$ is separable, then the representation $\rho_1$ is trivial, 
which can be derived from the fact that all continuous 
separable unitary representations of the Banach--Lie group 
$\U(\cH)/\U_\infty(\cH)$ are trivial if $\cH$ is separable 
(cf.\ \cite{Pi88} and also Theorem~\ref{thm:8.1b} below). 
\end{description}

We collect various auxiliary results in appendices. 
In Appendix~\ref{app:a} we discuss operator-valued positive 
definite functions on Lie groups. The main result 
is Theorem~\ref{thm:extension} asserting that analytic local 
positive definite functions extend to global ones. 
This generalizes the corresponding result for the scalar case 
in \cite{Ne10b}. Its applications to holomorphically induced 
representations are developed in Appendix~\ref{app:c}, 
where we show that holomorphic inducibility of $(\rho,V)$ 
can be characterized 
in terms of positive definiteness of a $B(V)$-valued function  
on some identity neighborhood of $G$. This is a key tool in 
Sections~\ref{sec:7} and \ref{sec:8}. 

For the convenience of the reader we provide in Appendix~\ref{app:d} 
a description of various classical groups 
of operators on Hilbert spaces over $\K \in \{\R, \C, \H\}$. 
Appendix~\ref{app:e} provides a complete discussion 
of the bounded unitary representations of 
the unitary group $\U_p(\cH)$ of 
an infinite dimensional real, complex or quaternionic Hilbert space 
for $1 < p \leq \infty$. 
The irreducible representations are 
parametrization in terms of highest weights. 
For $\K = \C$ this was done in \cite{Ne98}, 
and for $\K = \R$ and $\H$, these results are new but quite direct 
consequences of the complex case. We show in particular that 
all bounded representatations of these groups are direct sums 
of irreducible ones. If $\cH$ is separable, this is true 
for any continuous unitary representation of $\U_\infty(\cH)$ 
(\cite{Ki73}), but for $p < \infty$, the topological groups 
$\U_p(\cH)$ are not of type I (cf.\ \cite{Bo80}, \cite{SV75}). 
In Appendix~\ref{app:g} we finally recall the special 
features of separable representations of unitary groups. 
We conclude this paper with a discussion of some open problems 
and some comments on variations of the concept of a hermitian 
Lie group. 

The classification results of the present paper also 
contribute to the Olshanski--Pickrell program of classifying the unitary 
representations of automorphism groups of Hilbert symmetric spaces 
(cf.\ \cite{Pi87, Pi88, Pi90, Pi91}, \cite{Ol78, Ol84, Ol88, Ol89, Ol90}). 
For symmetric spaces $M = G/K$ of finite rank, Olshanski has shown in 
\cite{Ol78, Ol84} that the so-called {\it admissible} 
unitary representations of $G$, i.e., representations whose restriction to 
$K$ is tame,  lead to so-called {\it holomorphic} representations 
of the automorphism group $G^\sharp$ of a hermitian symmetric space 
$M^\sharp = G^\sharp/K^\sharp$ 
containing $M$ as a totally real submanifold. 
The irreducible holomorphic representations of the automorphism group of 
$M^\sharp$ turn out to be highest weight representations, 
which correspond to the representations showing up 
for type I$_{\rm fin}$ in our context. 
For spaces $M$ of infinite rank, the classification of the admissible 
representations is far less complete, although Olshanski 
formulates in \cite{Ol90} quite precise conjectures. 
These conjectures suggest that one should also try to understand the 
semibounded representations of mapping groups of the form 
$C^\infty(\bS^1, G)$, where $G$ is a hermitian Lie group.


{\bf Acknowledgement:} We thank Daniel Belti\c t\u a, Hasan G\"undogan, 
St\'ephane Merigon, Stefan Wagner and Christoph Zellner 
for various comments on this paper and for several discussions 
on its subject matter.

{\bf Notation and conventions:} 
For an open or closed convex cone $W$ in a real 
vector space $V$ we write 
$H(W) := \{ x \in V \: x + W = W \}$ for the {\it edge of $W$}. 

If $\g$ is a real Lie algebra and $\g_\C$ its complexification, we 
write $\oline z := x -iy$ for $z = x + iy$, $x,y \in \g$, 
and $z^* := - \oline z$. 

For a Lie group $G$ with Lie algebra $\g$ 
and a topological vector space $V$, we associate to each 
$x \in \g$ the left invariant 
differential operator on $C^\infty(G,V)$ defined by 
\[ (L_x f)(g) := \derat0 f(g\exp(tx)) \quad \mbox{ for } \quad 
x \in \g. \]
By complex linear extension, we define the operators 
\[L_{x+ iy} := L_x + i L_y \quad \mbox{ for } \quad z = x + iy \in \g_\C, 
x,y \in \g. \]
Accordingly, we define 
\[ (R_x f)(g) := \derat0 f(\exp(tx)g) \quad \mbox{ for } \quad 
x \in \g \] 
and $R_z$ for $z \in \g_\C$ by complex linear extension. 
For a homomorphism $\phi \: G \to H$ of Lie groups, 
we write $\L(\phi) \: \L(G) = \g \to \L(H) = \fh$ for the 
derived homomorphism of Lie algebras. 

For a set $J$, we write $|J|$ for its cardinality. 

\tableofcontents

\section{Hermitian Lie groups} \mlabel{sec:herm}

Let $G$ be a connected 
Banach--Lie group with Lie algebra $\g$, endowed with an 
involution~$\theta$. We write $K := (G^\theta)_0$ for the identity 
component of the subgroup of $\theta$-fixed points in 
$G$ and 
\[ \g = \fk \oplus \fp, \quad 
\fk := \ker(\L(\theta)-\1), \quad 
\fp := \ker(\L(\theta)+\1) \] 
for the corresponding eigenspace decomposition of the Lie algebra. 
Then $\fp \subeq \g$ is a closed $\Ad(K)$-invariant subspace 
and we write $\Ad_\fp \: K \to \GL(\fp)$ for the corresponding 
representation of $K$ on $\fp$. 

\begin{defn}  \mlabel{def:1.1}
(a) We call a triple $(\g,\theta,d)$ consisting of a Banach--Lie algebra 
$\g$, an involution $\theta$ of $\g$ and an element 
$d \in \z(\fk)$ a {\it hermitian} Lie algebra if the following 
conditions are satisfied: 
\begin{description}
\item[\rm(H1)] $\fp$ is a complex Hilbert space. 
\item[\rm(H2)] $\ad_\fp \: \fk \to \gl(\fp)$ is  
a representation of $\fp$ by bounded skew-hermitian operators. 
\item[\rm(H3)] $[d,x] = ix$ for every $x \in \fp$.   
\item[\rm(H4)] $\ker \ad_\fp \subeq \z(\fk)$. 
\end{description}
We say that the hermitian Lie algebra $(\g,\theta,d)$ is {\it irreducible} 
if, in addition, the representation of $\fk$ on $\fp$ is irreducible. 

(b) A triple $(G,\theta,d)$ of a connected Lie group $G$, an involutive 
automorphism $\theta$ of $G$ and an element $d \in \g$ 
is called a {\it hermitian} Lie group if 
the corresponding triple $(\g, \L(\theta), d)$ is hermitian. 
Since $K$ is connected, (H2) means that 
\begin{description}
\item[\rm(H2')] $\Ad_\fp$ is a unitary representation of $K$. 
\end{description}

Since $\fp$ is a closed complement of $\fk = \L(K)$ in $\g$, the quotient 
$M := G/K$ carries a natural Banach manifold structure for which 
the quotient map $q \: G \to M$ is a submersion. Moreover, 
the complexification 
$\g_\C$ decomposes into $3$-eigenspaces of $\ad d$: 
\begin{equation}
  \label{eq:tridec}
\g_\C = \fp^+ \oplus \fk_\C \oplus \fp^-, \quad 
\fp^\pm := \ker (\ad d \mp i \1),  
\end{equation}
and all these eigenspaces are $\Ad(K)$-invariant. 
Therefore $M = G/K$ carries a natural 
complex structure for which the tangent space $T_{\1K}(M)$ 
in the base point $\1 K$ can be identified with 
$\fp^+ \cong \g_\C/(\fk_\C + \fp^-)$, resp., with 
$\fp$, endowed with the complex structure defined by $\ad d$ 
(\cite[Thm.~15]{Bel05}). 
As the scalar product on 
$\fp$ is $\Ad(K)$-invariant, it defines a Riemannian structure on $M$. 
Endowed with these structures, $M$ becomes a 
{\it hermitian symmetric space} (cf.\ \cite{Ka83}). 
The fact that the triple bracket $[x,[y,z]]$ corresponds to the 
curvature now implies that, for $x \in \fp$, the operator  
$(\ad x)^2 \: \fp  \to \fp$ is hermitian. 
The symmetric space is said to be {\it flat} if 
$[\fp,[\fp,\fp]]=0$, {\it positively curved} if the operators 
$(\ad x)^2\res_\fp$, $x \in \fp$, are negative semidefinite, 
and {\it negatively curved} if these operators 
are positive semidefinite (cf.\ \cite{Ne02c}). 
\end{defn}

\begin{rem} \mlabel{rem:1.2} (a) 
If $G$ is simply connected, then $M = G/K$ is also simply 
connected because the connectedness of $K$ implies that the natural 
homomorphism $\pi_1(G) \to \pi_1(G/K)$ is surjective. 

(b) We neither require that $\fk = [\fp,\fp]$ nor that 
$[\fp,\fp]$ is dense in $\fk$. 

(c) If $(\g,\theta, d)$ is a hermitian Lie algebra, then 
$\theta = e^{\pi \ad d}$ follows immediately from the definitions. 
In particular, every connected Lie group $G$ with Lie algebra 
$\g$ carries the involution $\theta_G(g) := (\exp\pi d)g(\exp-\pi d)$ 
with $\L(\theta_G) = \theta$. 

(d) From $\fk = \ker(\ad d) \supeq \z(\g)$ (H3), we derive that 
$\z(\g) = (\ker\ad_\fp) \cap \z(\fk)$, so that (H4) 
is equivalent to $\ker \ad_\fp = \z(\g)$. If, in addition,
$(\g,\theta, d)$ is irreducible, then 
Schur's Lemma implies that 
$\ad_\fp(\z(\fk)) = \R i \1 = \R \ad_\fp d$, which leads to 
\begin{equation}
  \label{eq:centk}
\fz(\fk) = \R d \oplus \fz(\g).
\end{equation}
\end{rem}

\begin{exs} \mlabel{ex:1.3} (Flat case) 
(a) Suppose that $\cH$ is a complex Hilbert space, 
$K$ a connected Banach--Lie group and  
$\pi \: K \to \U(\cH)$ a norm continuous (hence smooth) 
representation of $K$ on $\cH$ with 
$\ker(\dd\pi) \subeq \z(\fk)$. Then 
$G := \cH \rtimes_\pi K$ is a symmetric Banach--Lie group with respect to 
the involution $\theta(v,k) := (-v,k)$. In this case 
$G^\theta = \{0\} \times K$, so that $G/K \cong \cH$. 
If there exists an element $d \in \z(\fk)$ with $\dd\pi(d) = i \1$, 
then $(G,\theta,d)$ is a hermitian Lie group. 
If the latter condition is not satisfied, we may replace 
$K$ by $\R \times K$ and extend the given representation 
by $\pi(t,k) := e^{it}\pi(k)$. 

(b) We can modify the construction under (a) as follows. 
We define the {\it Heisenberg group of $\cH$} by 
\[  \Heis(\cH) := \R \times \cH, \]
endowed with the multiplication 
$$ (t,v) (t',v') := (t + t' + \shalf\omega(v,v'), v + v'), \quad 
\omega(v,v') := 2\Im \la v,v'\ra. $$
Then $K$ acts smoothly by automorphisms on $\Heis(\cH)$ via 
$\alpha_k(t,v) = (t, \pi(k)v)$, and we obtain a 
Lie group $\hat G := \Heis(\cH) \rtimes_\alpha K$ with an involution 
defined by $\theta(t,v,k) := (t,-v,k)$ and 
$\hat K := \hat G^\theta = \R \times \{0\} \times K \cong \R \times K$. 

For $K = \T \1 \subeq \U(\cH)$, we  obtain in particular 
the {\it oscillator group} $\Osc(\cH) = \Heis(\cH) \rtimes \T$ 
of~$\cH$. 
\end{exs}

\begin{exs} \mlabel{ex:1.4} 
Let $\cH$ be a complex 
Hilbert space which is a direct sum $\cH = \cH_+ \oplus \cH_-$. 

(a) (Positively curved case) The {\it restricted unitary group} 
$\U_{\rm res}(\cH) = \U(\cH) \cap \GL_{\rm res}(\cH)$ 
from Appendix~\ref{app:d.3} is hermitian with 
$d := \shalf\diag(i\1,-i\1)$ and the involution 
\[ \theta\pmat{a & b \\ c &d}
:= \pmat{\1 & 0 \\ 0 &-\1}\pmat{a & b \\ c &d}\pmat{\1 & 0 \\ 0 &-\1}
= \pmat{a & -b \\ -c &d} \] 
whose group of fixed points is the connected group 
$\U_{\rm res}(\cH)^\theta 
\cong \U(\cH_+) \times \U({\cH}_{-}).$
We also see that 
\[ \fp = \Big \{ \pmat{0 & z \\ -z^* & 0} \: z \in B_2(\cH_-, \cH_+)\Big\} 
\cong B_2(\cH_-, \cH_+) \cong \cH_+ \hat\otimes \cH_-^* \] 
is a complex Hilbert space, and the adjoint representation 
$\Ad_\fp(k_+, k_-) z = k_+ z k_-^{-1}$ 
is irreducible because it is a tensor product of two irreducible 
representations. 

(b) (Negatively curved case) The restricted hermitian group 
$\U_{\rm res}(\cH_+, \cH_-)$ (Appendix~\ref{app:d.3}) 
is also hermitian 
with respect to the same involution and the same element $d \in \fk$ 
as in (a). We also note that 
$\U_{\rm res}(\cH_+, \cH_-)^\theta  
= \U_{\rm res}(\cH)^\theta 
\cong \U(\cH_+) \times \U(\cH_-).$
\end{exs}

\begin{defn} \mlabel{def:1.5} 
(a) A {\it symmetric Lie algebra} is a pair $(\g,\theta)$, where
$\g$ is a Lie algebra and $\theta \in \Aut(\g)$ 
(the group of topological automorphisms) 
with $\theta^2 = \id_\g$.
If $\g = \fk \oplus \fp$ is the decomposition into $\theta$-eigenspaces, 
then $\g^c := \fk + i \fp \subeq \g_\C$ with the involution 
$\theta^c(x + iy) := x - iy$, $x \in \fk$, $y \in \fp$, is called 
the {\it $c$-dual symmetric Lie algebra} $(\g^c,\theta^c)$ of 
$(\g,\theta)$. 

(b) If $(\g,\theta,d)$ is hermitian, then the $c$-dual 
Lie algebra $(\g^c,\theta^c,d)$ is also hermitian, 
called the {\it $c$-dual hermitian Lie algebra}. 

(c) For any hermitian Lie group $(G,\theta, d)$, the semidirect 
product $G^m := \fp \rtimes_{\Ad_\fp} K$ is called the 
associated {\it motion group}. It is hermitian with respect to the 
involution $\theta^m(x,k) := (-x,k)$ (Example~\ref{ex:1.3}). 
\end{defn}

\begin{rem} Note that $c$-duality changes the sign of the curvature.
In fact, multiplying $x \in \fp$  with $i$ leads to the 
operator $(\ad ix)^2 = - (\ad x)^2$ on $i\fp\cong \fp$. 
\end{rem}

\begin{ex} \mlabel{ex:1.7} The restricted unitary group 
$\U_{\rm res}(\cH_+\oplus \cH_-)$ and the restricted 
pseudo-unitary group $\U_{\rm res}(\cH_+, \cH_-)$ 
from Example~\ref{ex:1.4} are $c$-duals of each other: 
\[ \fu_{\rm res}(\cH_+\oplus \cH_-)^c = \fu_{\rm res}(\cH_+, \cH_-).\] 
\end{ex}

\section{$JH^*$-triples and their classification} \mlabel{sec:2}

A classification of the simply connected hermitian symmetric 
spaces has been obtained by W.~Kaup in \cite{Ka83} 
(cf.\ Theorem~\ref{thm:classif} below). 
Kaup calls a connected Riemannian Hilbert manifold 
$M$ a {\it symmetric hermitian manifold} (which is the same 
as a hermitian symmetric space) if it carries a 
complex structure and for each $m \in M$ there exists a 
biholomorphic involutive isometry $s_m$ having $m$ as an isolated 
fixed point. According to \cite{Ka77}, the category of simply connected 
symmetric hermitian manifolds with base point is equivalent to the 
category of $JH^*$-triples, 
and these are classified in \cite{Ka83}. We briefly recall the 
cornerstones of Kaup's classification. 

\begin{defn} (a) A complex Hilbert space $U$ endowed with a 
real trilinear map 
\[ \{ \cdot,\cdot,\cdot\} \: U^3 \to U, \quad 
(x,y,z) \mapsto (x \square y)z := \{x,y,z\} \]  
which is complex linear in $x$ and $z$ and antilinear in $y$ 
is called a  {\it $JH^*$-triple} if the following conditions are satisfied 
for $a,b,c,x,y,z \in U$ (\cite[Sect.~1]{Ka81}): 
\begin{description}
\item[\rm(JH1)] $\{x,y,z\} = \{z,y,x\}$. 
\item[\rm(JH2)] $\{x,y,\{a,b,c\}\} = \{\{x,y,a\}, b,c\}  
- \{ a, \{y,x,b\},z\} + \{ a,b,\{x,y,c\}\}$. 
\item[\rm(JH3)] The operators $x \square x$ are hermitian. 
\end{description}

We write $\Aut(U)$ for the {\it automorphism group} of the $JH^*$-triple 
$U$, i.e., for the group of all unitary operators $g$ on $U$ satisfying 
$g\{x,y,z\} = \{gx,gy,gz\}$ for $x,y,z \in U$. 
We define the {\it spectral norm} on $U$ by $\|z\|_\infty := \sqrt{\|z \square z\|}$. 

(b) For every $JH^*$-triple, we define the $c$-dual $U^c$ as the same 
complex Hilbert space $U$, endowed with the triple product 
\[ \{x,y,z\}^c := - \{x,y,z\}.\] 
This turns $U^c$ into a $JH^*$-triple with 
$x \square^c x = - x \square x$. 
\end{defn}

\begin{rem} \mlabel{rem:2.2} 
In \cite[Thm.~3.9]{Ka83}, Kaup shows that 
every $JH^*$-triple $U$ is the orthogonal direct sum 
\[ U = \z(U) \oplus \big(\hat\oplus_{j \in J} U_j\Big),\] 
where the $U_j$, $j \in J$, are the simple triple ideals. 
A triple ideal is a closed subspace invariant under 
$U \square U$, and the simple triple ideals are the minimal 
non-zero triple ideals on which the triple product is non-zero.
Further, 
\[ \z(U) := \{ z \in U \: (\forall x,y \in U)\, \{x,y,z\} = 0\},\] 
the {\it center of $U$}, is flat, 
i.e., its triple product vanishes. 
Since each automorphism of $U$ 
permutes the simple ideals $(U_j)_{j \in J}$ 
and preserves the center $\z(U)$, it is easy to see that 
the Lie algebra of the 
Banach--Lie group $\Aut(U)$ preserves each simple 
ideal $U_j$, and therefore the identity component 
$\Aut(U)_0$ also preserves each simple ideal 
$U_j$.\begin{footnote}{For $X \in \aut(U)$ each operator $e^{tX}$, 
$t \in \R$,  
permutes the spaces $U_j$. Let $P_j \: U \to U_j$ denote the orthogonal 
projection onto $U_j$. Then, for each $j$, the operator 
$P_j e^{tX} P_j^*$ on $U_j$ is non-zero if $t$ is sufficiently small 
because it converges uniformly to $\1$. For any such $t$ we then have
$e^{tX}U_j \subeq U_j$, and this implies $e^X U_j = U_j$. 
}\end{footnote}
\end{rem}

\begin{rem} \mlabel{rem:2.3} 
(a) (From hermitian Lie groups to $JH^*$-triples; \cite{Ka77}) 
We claim that, for each hermitian Lie algebra $(\g,\theta,d)$, 
$\fp$ carries a natural structure of a $JH^*$-triple. 
First we use the embedding 
$\fp \into \g_\C, x \mapsto x - i Ix$, where $Ix = [d,x],$ 
to obtain an isomorphism of the complex Hilbert space $\fp$ 
with the complex abelian subalgebra $\fp^+ \subeq \g_\C$ from 
\eqref{eq:tridec}. 
Therefore it suffices to exhibit a natural $JH^*$-structure on $\fp^+$. 

For $x,y,z \in \fp^+,$ we put 
\[ \{x,y,z\} := [[x,\oline y],z] \in \fp^+. \] 
This map is complex linear in the first and third and antilinear 
in the second argument.
Using the grading 
$\g_\C = \fp^- \oplus \fk_\C \oplus \fp^+$, the two conditions 
(JH1/2) are easily verified.\begin{footnote}
  {In \cite{Ka83}, Kaup uses the formula 
$\{x,y,z\} := \shalf[[x,\oline y],z]$.}
\end{footnote} We also note that 
\[ x \square x = \ad([x,\oline x]),\] 
so that (H2) and $[x,\oline x] \in i\fk$ imply that 
$x \square x$ is hermitian. Therefore $(\fp^+, \{\cdot,\cdot,\cdot\})$ 
is a $JH^*$-triple. 

(b) (From $JH^*$-triples to hermitian Lie groups) 
If, conversely, $U$ is a $JH^*$-triple and 
$\oline U$ denotes $U$, endowed with the opposite complex structure, 
then 
the Kantor--Koecher--Tits 
construction leads to the complex Banach--Lie algebra 
\[ \g_\C(U) := U \oplus \aut(U)_\C \oplus \oline U,\] 
endowed with the bracket 
\[ [(x,A,y), (x',A',y')] = (Ax'-A'x, x \square y' - x' \square y 
+ [A,A'], -A^*y'+(A')^* y).\]
Since 
$(x \square x)^* = x \square x$, polarization leads to 
$(x \square y)^* = y \square x$ for $x,y \in U$. 
This implies that 
$\sigma(x,A,y) := (y,-A^*,x)$ defines an antilinear involution 
on $\g_\C(U)$, which leads to the real form 
\[ \g := \g(U) := \g_\C(U)^\sigma 
= \{ (x,A,x) \in \g_\C(U) \: x \in U, A^* = -A\}.\] 
This real form has an involution $\theta$, defined by 
$\theta(x,A,x) := (-x,A,-x)$ with 
$\g^\theta \cong \aut(U)$ and $\g^{-\theta} \cong U$. With 
the element $d \in \fz(\aut(U))$ given by $dx = ix$, we 
now obtain a hermitian Lie algebra $(\g,\theta,d)$ with 
\[ \fp^+ = U \times \{(0,0)\}, \quad \fk_\C = \aut(U)_\C 
\quad \mbox{ and }  \quad \fp^- = \{(0,0)\} \times U.\] 
Note that, for $x,y,z \in U \cong \fp^+$, we have 
 $[[x,\oline y],z] = (x\square y)z = \{x,y,z\}$.
 \begin{footnote}
   {For $A \in \aut(U)$ we have 
$A\{x,y,z\} = \{Ax,y,z\} + \{x,Ay,z\} + \{x,y,Az\},$ 
so that complex linear extension leads for $A \in \aut(U)_\C$ to 
$A\{x,y,z\} = \{Ax,y,z\} - \{x,A^*y,z\} + \{x,y,Az\}.$
}
 \end{footnote}
\end{rem}

\begin{defn}
We call the irreducible hermitian Lie algebra $(\g,\theta, d)$ 
{\it full} if 
\begin{description}
\item[\rm(F1)]  $[\fp,[\fp,\fp]]\not=0$, i.e., 
$\fp$ is a simple $JH^*$-triple (Remark~\ref{rem:2.2}), and 
\item[\rm(F2)] $\ad_\fp(\fk) = \aut(\fp)$ 
is the Lie algebra of the automorphism 
group $\Aut(\fp)$ of the $JH^*$-triple $\fp$. 
\end{description}
\end{defn}

\begin{rem} \mlabel{rem:zent} (a) Suppose that $(\g,\theta, d)$ is full. 
In view of Schur's Lemma, the 
irreducibility of the $\fk$-module $\fp$ implies that 
$\z(\aut(\fp)) = \R i \1 = \R \ad_\fp d \subeq \ad_\fp(\z(\fk))$. 
Condition (F2) also implies the converse inclusion, so that 
\[ \z(\aut(\fp)) = \ad_\fp(\z(\fk)) = \R i \1.\] 
In particular, we obtain 
\begin{equation}
  \label{eq:transfer}
\aut(\fp)/\z(\aut(\fp)) \cong \fk/\z(\fk).
\end{equation}

(b) If $(G,\theta, d)$ is a hermitian 
Lie group satisfying (F2), then 
$\z(\g) = \ker \ad_\fp$ (Remark~\ref{rem:1.2}(d)) 
leads to 
$\ad \g \cong \aut(\fp) \oplus \fp = \g(\fp)$, endowed with the Lie bracket 
\[ [(A,x), (A',x')] := ([A,A'] + \ad_\fp([x,x']), Ax' - A'x)\]  
(cf.\ Remark~\ref{rem:2.3}(b)). 
\end{rem}

For an irreducible hermitian Lie group $(G,\theta,d)$,  
the connectedness of 
$K$ and the irreducibility of the $K$-representation on $\fp$ implies  
that either $\fp$ is a simple $JH^*$-triple or flat. 
Kaup's classification also implies that a simple $JH^*$-triple 
either is {\it positive}, i.e., all operators $x \square x$ are positive,  
or {\it negative}, i.e., all these operators are negative 
(\cite[p.~69]{Ka83}). Since we shall need it later on, we recall 
the infinite dimensional part of Kaup's classification. 

\begin{thm} \mlabel{thm:classif} 
Each simple infinite dimensional positive $JH^*$-triple $U$ is isomorphic 
to one of the following, where 
$\cH_\pm$ and $\cH$ are complex Hilbert spaces and 
$\sigma$ is a conjugation on $\cH$, i.e., an antilinear isometric involution. 
\begin{description}
\item[\rm(I)] The space $B_2(\cH_-, \cH_+)$ of Hilbert--Schmidt 
operators $\cH_- \to \cH_+$ with 
\begin{equation}
  \label{eq:struc}
\la A,B \ra = \tr(AB^*) \quad \mbox{ and } \quad 
\{A,B,C\} = \shalf(AB^*C + CB^*A).
\end{equation}
In this case $\|A\|_\infty = \|A\|$ is the operator norm. 
For $\cH_+ \not\cong \cH_-$, 
the automorphism group consists of the maps of the form 
$\alpha(A) = g_+ A g_-$ for $g_\pm \in \U(\cH_\pm)$, and 
for $\cH_+ = \cH_- = \cH$, we obtain an additional automorphism 
by $A \mapsto A^\top := \sigma A^*\sigma$. 
\item[\rm(II)] The space $\Skew_2(\cH) := \{ A \in B_2(\cH) 
\: A^\top= - A\}$ of skew-symmetric Hilbert--Schmidt operators, 
which is a $JH^*$-subtriple of $B_2(\cH)$. The automorphisms are the maps  
$\alpha(A) = g A g^{-1}$ for $g \in \U(\cH)$. 
\item[\rm(III)] The space $\Sym_2(\cH) := \{ A \in B_2(\cH) 
\: A^\top=  A\}$ of symmetric Hilbert--Schmidt operators, 
which is another $JH^*$-subtriple of $B_2(\cH)$. 
The automorphisms are the maps  
$\alpha(A) = g A g^{-1}$ for $g \in \U(\cH)$. 
\item[\rm(IV)] A complex Hilbert space $U$ with a conjugation $\sigma$ 
and 
\[ \{x,y,z\} := \la x,y\ra z + \la z, y\ra x - \la x,\sigma(z)\ra 
\sigma(y).\] 
Then $\|z\|_\infty^2 = \la z,z \ra + \sqrt{\la z,z\ra^2 - |\la 
z, \sigma(z)\ra|^2}$ 
and the automorphisms are of the form 
$\alpha(v) = z g(v)$ for $z \in \T$ and 
$g \in \OO(\cH^\sigma)$. 
\end{description}
The negative simple $JH^*$-triples are the $c$-duals of the 
positive ones and have the same automorphism group. 
\end{thm}

\begin{prf} The classification can be found in \cite[Thm.~3.9]{Ka83}; 
see also \cite[p.~474]{Ka81} for the type IV case. 

For a description of the automorphism groups of 
simple $JH^*$-triples, we refer to \cite[pp.~199/200]{Ka97} 
and \cite[p.~91]{Ka75}. 
To see that Kaup's list also describes the automorphism 
groups of the associated $JH^*$-triples, we 
note that any automorphism of $U$ preserves 
the spectral norm on $U$, hence extends to the bidual 
$U_\infty''$ of the completion $U_\infty$ of $U$ 
with respect to the spectral norm 
(see also \cite[Thm.~V.11]{Ne01a}). Conversely, 
Kaup's list shows that all automorphisms of 
$U_\infty''$ preserve the subspace $U$ on which they 
induce $JH^*$-triple automorphisms. It is clear that $\Aut(U^c) = \Aut(U)$. 
\end{prf}

\subsection*{Automorphisms of symmetric Hilbert domains}

\begin{rem}\mlabel{rem:2.7} For every positive $JH^*$-triple $U$, the set 
\[ \cD := \{ x \in U \: \|x\|_\infty < 1 \} \] 
is a {\it symmetric Hilbert domain}, i.e., for each point 
$m \in \cD$, there exists an involutive biholomorphic map 
$s_m$ having $m$ as an isolated fixed point (\cite[p.~71]{Ka83}). 
Since this is clear 
for $m = 0$ with $s_0(x) = -x$, the main point is that 
the group $\Aut(\cD)$ of biholomorphic maps acts transitively 
on $\cD$. This group is a Banach--Lie group and 
$\theta(g) := s_0 g s_0$ defines an involution on 
$\Aut(\cD)$ with $\Aut(\cD)^\theta = \Aut(U)$ and 
the corresponding decomposition 
\[ \aut(\cD) := \L(\Aut(\cD)) \cong \aut(U) \oplus U, \quad \aut(\cD)^\theta
 = \aut(U), 
\quad \fp =  U. \] 
It is easy to see that 
$(\Aut(\cD)_0, \theta, i\1)$ is a hermitian Lie group, 
and if $U$ is simple, then it is full. We also derive 
that $\aut(\cD) \cong \g(U)$ (Remark~\ref{rem:2.3}(b)). 
\end{rem}

We now describe the automorphism groups of infinite dimensional 
irreducible symmetric Hilbert domains which correspond to the 
four families of infinite dimensional 
positive simple $JH^*$-triples. 

\nin {\bf Type I}: The group $G := \U_{\rm res}(\cH_+, \cH_-)$ 
(cf.\ Example~\ref{ex:1.4}(b)) 
acts transitively on 
\[ \cD := \{ Z \in B_2(\cH_-,\cH_+) \: \|Z\| < 1\}\]  
by fractional linear transformations 
\[ gZ = (aZ + b)(cZ + d)^{-1} \quad \mbox{ for } 
\quad g =\pmat{a & b \\ c & d}\]   
(cf.\ \cite{NO98}). 
The stabilizer of $0$ in $G$ is $K:= \U(\cH_+) \times \U(\cH_-)$, 
so that $\cD \cong G/K$. 
With Theorem~\ref{thm:classif} we see that 
$\Aut(\cD)_0 \cong \U_{\rm res}(\cH_+, \cH_-)/\T \1$, so that  
$\U_{\rm res}(\cH_+,\cH_-)$ is a central $\T$-extension of 
$\Aut(\cD)_0$ which is easily seen to be a hermitian Lie group 
with respect to $d := \frac{i}{2}\diag(\1,-\1)$. 
From Theorem~\ref{thm:classif} it also follows that 
$\U_{\rm res}(\cH_+,\cH_-)$ is full. 

\nin {\bf Type II}: For a conjugation $\sigma$ on the complex Hilbert space 
$\cH$, the subgroup $G := \OO^*_{\rm res}(\cH) \subeq \U_{\rm res}(\cH,\cH)$ 
(cf.\ Appendix~\ref{app:d.3}) acts transitively on 
\[ {\cal D} = \{ Z \in B_2({\cal H}) \: Z^\top = -Z, \|Z\| < 1 \}. \] 
The stabilizer of $0$ is 
$ K := \left\{ \pmat{ g & 0 \cr 0 & g^{-\top} \cr} \:
g \in \U({\cal H}) \right\} \cong \U(\cH),$ 
so that $G/K \cong \cD$. In this case 
$\Aut(\cD)_0 \cong \OO^*_{\rm res}(\cH)/\{\pm \1\}$ 
and $\OO^*_{\rm res}(\cH)$ is a hermitian Lie group 
with respect to $d := \frac{i}{2}\diag(\1,-\1)$. 

\nin {\bf Type III}: Likewise, the subgroup 
$\Sp_{\rm res}(\cH)\subeq \U_{\rm res}(\cH,\cH)$ 
(cf.\ Appendix~\ref{app:d.3}) acts transitively on 
\[ {\cal D} = \{ Z \in B_2({\cal H}) \: Z^\top = Z, \|Z\| < 1 \}. \] 
The stabilizer $K$ of $0$ is isomorphic to $\U(\cH)$, 
$\Aut(\cD)_0 \cong \Sp_{\rm res}(\cH)/\{\pm \1\}$, 
and $\Sp_{\rm res}(\cH)$ is a hermitian Lie group 
with respect to $d := \frac{i}{2}\diag(\1,-\1)$. 

\nin {\bf Type IV}: We consider the 
orthogonal group 
$G := \OO(\R^2,\cH_\R)_0$ of the indefinite quadratic form 
\[ q(x,v) := \|x\|^2 - \|v\|^2 \] 
on $\R^2 \oplus \cH_\R$, where $\cH_\R$ is an infinite dimensional 
real Hilbert space. Then 
\[ d := \Big(\pmat{0 & -1 \\ 1&0},0\Big) \in \g := \fo(\R^2 \oplus \cH_\R)\] 
satisfies $(\ad d)^3 + \ad d=0$, so that 
$\theta := e^{\pi \ad d}$ defines an involution on $G$ 
with 
\[ K := G^\theta  = \SO_2(\R) \times \OO(\cH_\R). \]  
Elements of the Lie algebra $\g = \fo(\R^2, \fo(\cH_\R))$ 
have a block structure 
\[ x = \pmat{a & b^\top \\ b &c},
\quad \mbox{ with } \quad 
b = (b_1, b_2) \in \cH_\R^2 \cong (\cH_\R)_\C.\] 
From 
\[ \Big[ \pmat{a & 0 \\ 0 & c}, \pmat{0 & b^\top \\ b & 0}\Big] 
= \pmat{0 & (cb-ba)^\top \\ cb-ba & 0} \] 
we derive that $a\in \so_2(\R)$ acts on 
$\fp \cong \cH_\R^2$ by $b \mapsto -ba$. 
Therefore $\ad d$ defines a complex structure on 
$\fp \cong \cH_\R^2$, acting by 
\[ I(v,w) = -(v,w)\pmat{0 & -1 \\ 1 & 0} 
= (-w,v).\] 
This leads to a natural identification of 
$\fp$ with $(\cH_\R)_\C$, where 
$\fo(\cH_\R) \subeq \fk$ acts in the natural way by complex 
linear operators. Since the representation of 
$\fo(\cH_\R)$ on $\cH_\R$ is absolutely irreducible, 
i.e., its commutant is $\R \1$, the representation of 
$\fk$ on the complex Hilbert space $\fp$ is irreducible. 
We therefore obtain the hermitian Lie group $(G,\theta,d)$. 

Putting $\cH := (\cH_\R)_\C \cong \fp$, we then obtain 
$\fp_\C \cong \cH^2$ with 
\[ \fp^\pm \cong \{ b \mp i I b \: b \in \cH \} 
= \{ (b,\mp i b) \: b \in \cH\}.\] 
In this realization, an easy matrix calculation leads 
with the isomorphism $\cH \to \fp^+, x \mapsto (x,-ix)$ 
for $x,y,z \in \fp^+ \cong \cH$ to 
\[ [[x,\oline y],z] = 2 (\la x,y \ra z + \la z,y\ra x - \la x, \oline z\ra \oline y).\] 
Therefore the associated $JH^*$-triple is of type IV 
(Theorem~\ref{thm:classif}).

\section{Central extensions of 
hermitian Lie groups} \mlabel{sec:centext} 

In Section~\ref{sec:6}-\ref{sec:8} below 
we show that for certain hermitian Lie groups $G$ 
there are natural central extensions $\hat G$ with a substantially 
richer semibounded unitary representation theory. 
This motivates our detailed discussion of central 
extensions in the present section. Its goal is to 
show that for every simple $JH^*$-triple 
$U$, the Lie algebra $\g(U)$ has a universal central 
extension $\hat\g(U)$ whose center is at most $2$-dimensional 
and which is {\it integrable} in the sense that it is the Lie algebra 
of a simply connected Banach--Lie group $\hat G(U)$. 
For these central Lie group extensions we give rather 
direct constructions based on \cite{Ne02}.

\begin{defn} \mlabel{def:h2} 
\rm (a) Let $\z$ be a topological vector space and 
$\g$ be 
a topological Lie algebra. A {\it continuous $\z$-valued $2$-cocycle} is a
continuous skew-symmetric function $\omega \: \g \times \g \to \z$ with 
$$ \omega([x,y], z) + \omega([y,z], x) + \omega([z,x], y) = 0 
\quad  \mbox{ for } \quad x,y,z \in \g. $$
It is called a {\it coboundary} if there exists  a continuous linear map 
$\alpha \: \g \to \z$ 
with $\omega(x,y) = \alpha([x,y])$ for all $x, y \in
\g$. We write $Z^2_c(\g,\z)$ for the space of continuous $\z$-valued 
$2$-cocycles, $B^2_c(\g,\z)$ for the subspace of coboundaries defined by
continuous linear maps, and define the {\it second continuous Lie
algebra cohomology space} 
$$ H^2_c(\g,\z) := Z^2_c(\g,\z)/B^2_c(\g,\z). $$

(b) If $\omega$ is a continuous $\z$-valued cocycle on $\g$, then 
we write $\fz \oplus_\omega \g$ for the topological Lie algebra
whose underlying topological vector space is the product space 
$\fz \times \g$, and whose bracket is defined by 
$$ [(z,x),(z',x')] = \big(\omega(x,x'),[x,x']\big). $$
Then $q \: \fz \oplus_\omega \g \to \g, (z,x) \mapsto x$ is a central
extension and $\sigma \:\g\to{\fz \oplus_\omega\g}, x \mapsto (0,x)$
is a continuous linear section of $q$. 

(c) A {\it morphism of central extensions} 
\[ \fz_1 \ssmapright{i_1} \hat\g_1 \ssmapright{q_1} \g, \qquad 
\fz_2 \ssmapright{i_2}  \hat\g_2 \ssmapright{q_2} \g \] 
is a continuous homomorphism $\phi \: \hat\g_1 \to \hat\g_2$ with 
$q_2 \circ \phi = q_1$ and $\phi \circ i_1 = i_2 \circ \psi$ 
for a continuous linear map $\psi \: \fz_1 \to \fz_2$. 
For $\fz_1 = \fz_2$ and $\psi = \id_{\fz_1}$, the morphism 
$\phi$ is called an {\it equivalence}. 
The space $H^2_c(\g,\fz)$ classifies equivalence classes of topologically 
split central extensions of $\g$ by $\fz$ 
(cf.\ \cite[Rem.~1.2]{Ne02b}). 

(d) A central extension $\fz \into \hat \g \to \g$ is called 
{\it universal} if, for every other central extension 
$\fz_1 \into \hat\g_1 \to \g$, there exists a unique 
morphism $\phi \: \hat\g \to \hat\g_1$ of central extensions. 
This condition readily implies that universal central extensions 
are unique up to isomorphism. 
According to \cite[Prop.~I.13]{Ne02}, a perfect 
Banach--Lie algebra
\begin{footnote}{We call a Lie algebra $\g$ {\it perfect} if 
$\g = [\g,\g]$, i.e., if it is spanned by all brackets. 
A topological Lie algebra $\g$ is called {\it topologically perfect} 
if the subspace $[\g,\g]$ is dense in $\g$.}
\end{footnote}
for which $H^2_c(\g,\R)$ is finite dimensional has a universal 
central extension $\fz \into \hat\g \onto \g$ satisfying 
$\fz \cong H^2_c(\g,\R)^*$. 

(e) We say that $\g$ is {\it centrally closed} if $H^2_c(\g,\R) = \{0\}$. 
From \cite[Lemma~1.11, Cor.~1.14]{Ne02b} it follows that a 
perfect central extension 
$q \:  \hat\g \to \g$ with finite dimensional 
kernel which is centrally closed is universal. 
\end{defn}

\begin{defn} \mlabel{def:3.1} (a) Let $(\g,\theta,d)$ be a hermitian 
Lie algebra, write $\la \cdot,\cdot\ra$ for  the scalar product 
on $\fp$ and let $\omega_\fp(x,y) := 2\Im \la x,y\ra$ be the corresponding 
symplectic form on $\fp$. We extend $\omega_\fp$ to a skew-symmetric 
form $\omega_\fp$ on $\g$ satisfying $i_x \omega_\fp =0$ for every $x \in \fk$. 
We claim that $\omega_\fp$ is a $2$-cocycle, i.e., an element of 
$Z^2(\g,\R)$. As $\omega_\fp$ is $\fk$-invariant, 
$i_x \dd \omega_\fp = -\dd(i_x \omega_\fp) = 0$ for $x \in \fk$. 
Hence it remains to see that 
\[ (\dd\omega_\fp)(x,y,z) = -\sum_\cyc \omega_\fp([x,y],z)\] 
vanishes for $x,y,z \in \fp$, but this follows from 
$[\fp,\fp] \subeq \fk$. 

We call $\omega_\fp$ the {\it canonical cocycle} of the  
hermitian Lie algebra $(\g,\L(\theta), d)$. 

(b) Let $(\g,\theta, d)$ be a hermitian Lie algebra and 
$\omega_\fp \in Z^2(\g,\R)$ be its canonical  $2$-cocycle 
(Definition~\ref{def:3.1}(a)). We write 
$\hat\g = \R \oplus_{\omega_\fp} \g$ for the corresponding central 
extension with the bracket 
\[ [(z,x), (z',x')] := (\omega_\fp(x,x'), [x,x']).\] 
Since $\omega_\fp$ is $\theta$-invariant, 
$\hat\theta(z,x) := (z,\L(\theta)x)$ 
defines an involution on the Banach--Lie algebra $\hat\g$ with 
eigenspace decomposition 
\[ \hat\fk = \R \oplus \fk 
\quad \mbox{ and } \quad \hat\fp = \fp.\] 
We thus obtain a hermitian Lie algebra $(\hat\g, \hat\theta, d)$. 
In Theorem~\ref{thm:exist} below we shall see how to obtain 
a corresponding Lie group. 
\end{defn}

\begin{lem} \mlabel{lem:cocyc} Let $(\g,\theta,d)$ be irreducible 
hermitian and 
$\omega \in Z^2_c(\g,\R)$ be a continuous $2$-cocycle. 
Then the following assertions hold: 
\begin{description}
\item[\rm(i)] The class $[\omega] \in H^2_c(\g,\R)$ can be represented by a 
cocycle vanishing on $\fk \times \fp$. 
\item[\rm(ii)] If $\omega$ vanishes on $\fk \times \fk$, 
then  $[\omega] \in \R [\omega_\fp]$, where $\omega_\fp$ is the canonical 
cocycle. 
\end{description}
\end{lem}

\begin{prf} (i) In view of \cite[Thm.~9.1]{Ne10e}, 
the group $e^{\R \ad d} = e^{[0,2\pi]\ad d}$ acts trivially on 
$H^2_c(\g,\R)$, so that 
\[ \tilde\omega(x,y) 
:= \frac{1}{2\pi} \int_0^{2\pi} \omega(e^{t \ad d}x, e^{t \ad d}y)\, dt \] 
is a $2$-cocycle with the same cohomology class as $\omega$ 
and which coincides with $\omega$ on $\fk \times \fk$. 
Since $\tilde\omega$ is $\ad d$-invariant, we obtain for $x \in \fk$ and 
$y \in \fp$: 
\[0 = \tilde\omega([d,x],y)+ \tilde\omega(x,[d,y])
= \tilde\omega(x,[d,y]),\] 
so that $[d,\fp] = \fp$ leads to $\tilde\omega(\fk,\fp) = \{0\}$. 

Suppose that $\omega$, and hence also $\tilde\omega$, vanishes  
on $\fk \times \fk$. From the Cartan formula we further obtain for 
$x \in \fk$ the relation 
$\cL_x \tilde\omega = \dd(i_x\tilde\omega) = 0,$
so that $\tilde\omega\res_{\fp \times \fp}$ is a $\fk$-invariant 
alternating form, hence extends to a $\fk$-invariant 
skew-hermitian form 
$$\tilde\omega_\C \: \fp \times \fp \to \C, \quad 
\tilde\omega_\C(x, y) := \tilde\omega(x,y)+ i \tilde\omega(-Ix,y), \quad 
 Ix = [d,x]. $$
This form can be represented by a skew-hermitian 
operator $A = - A^* \in B(\fp)$ via 
$$ \tilde\omega_\C(x,y) = \la Ax,y \ra \quad \mbox{ for } 
\quad x, y \in \fp. $$
Now the $\fk$-invariance of $\tilde\omega_\C$ implies that 
$A$ commutes with $\ad \fk$, so that the irreducibility 
of the $\fk$-module $\fp$ leads to 
$A \in i \R \1$. We conclude that there exists a 
$\lambda \in \R$ with 
$$ \tilde\omega(x,y) = \Re \la \lambda i x, y \ra 
= -\lambda \Im \la x, y \ra = - \frac{\lambda}{2} \omega_\fp(x,y), $$
so that $\tilde\omega = - \frac{\lambda}{2} \omega_\fp$. 
\end{prf}

\begin{rem} \mlabel{rem:3.5} 
(a) If $\g = \g_1 \oplus \g_2$ is a direct sum of topological 
Lie algebras, 
where $H^2_c(\g_j,\R)= \{0\}$ for $j =1,2$ and $\g_1$ is topologically 
perfect, then $H^2_c(\g,\R) = \{0\}$. 

In fact, if $\omega \in Z^2_c(\g,\R)$ is a cocycle, then the triviality 
of the cohomology of $\g_1$ and $\g_2$ implies that we may subtract 
a coboundary, so that we can 
assume that $\omega$ vanishes on $\g_j \times \g_j$ for $j =1,2$. 
The cocycle property further implies that 
\[  \omega([x_1, x_2], y) = 0 \quad \mbox{ for } 
\quad x_1, x_2 \in \g_1, y \in \g_2. \]
Since $\g_1$ is topologically perfect, we obtain $\omega = 0$. 

(b) \cite[Prop.~I.5]{Ne02}: If $\cH$ is a complex Hilbert space 
of any dimension, 
then the second continuous cohomology $H^2_c(\gl(\cH),\C)$ of the 
complex Banach--Lie algebra $\gl(\cH)$ vanishes.
This implies in particular that it also vanishes for all its real forms: 
\[ H^2_c(\fu(\cH),\R) = \{0\} \quad \mbox{ and } \quad 
H^2_c(\fu(\cH_+, \cH_-),\R) = \{0\}.\] 
\end{rem} 

We now apply all this to the Lie algebras $\fk$ corresponding to the 
hermitian operator groups from Section~\ref{sec:2}. 

\begin{ex} \mlabel{ex:3.5} \nin{\bf Type I:} 
(a) For $\g = \fu_{\rm res}(\cH_+, \cH_-)$,  
we have $\fk = \fu(\cH_+) \oplus \fu(\cH_-)$. 
If at least one space $\cH_\pm$ is infinite dimensional, 
then Remark~\ref{rem:3.5}(b) 
implies that $H^2_c(\fu(\cH_\pm),\R) =\{0\}$, 
and since $\fu(\cH_+)$ or $\fu(\cH_-)$ 
is perfect (\cite[Lemma~I.3]{Ne02}), (a) leads to 
$H^2_c(\fk,\R) = \{0\}$. 

Therefore every class $[\omega] \in H^2_c(\fu_{\rm res}(\cH_+, \cH_-),\R)$ 
can be represented by a cocycle vanishing on $\fk \times \fk$, so 
that Lemma~\ref{lem:cocyc} shows that 
$H^2_c(\fu_{\rm res}(\cH_+, \cH_-),\R) = \R [\omega_\fp]$. 
If not both $\cH_\pm$ are infinite dimensional, then 
$\fu_{\rm res}(\cH_+, \cH_-) = \fu(\cH_+, \cH_-)$, 
so that Remark~\ref{rem:3.5}(b) implies that its second 
cohomology is trivial. In this case we therefore have 
$H^2_c(\fu_{\rm res}(\cH_+, \cH_-),\R) = 0.$

If both spaces $\cH_\pm$ are infinite dimensional, then the 
cocycle $\omega_\fp$ leads to a central extension 
$\hat\g = \hat\fu_{\rm res}(\cH_+, \cH_-)$ 
of $\g = \fu_{\rm res}(\cH_+, \cH_-)$ with 
$\hat\fk \cong \R \oplus \fu(\cH_+) \oplus \fu(\cH_-)$, 
so that $\dim \z(\hat\fk) = 3$ and $\dim \z(\hat\g) = 2$. 

We claim that the Lie algebra 
$\hat\g$ is centrally closed. 
To this end, we first recall  the central $\R$-extension 
of $\g$ described in \cite{Ne02} whose cocycle is defined 
in terms of the $(2 \times 2)$-block 
matrix structure of the elements 
$x = \pmat{x_{11} & x_{12} \\ x_{21} & x_{22}}$ 
of $\g$ by  
\[ \omega(x,y) := \tr([x,y]_{11}) = \tr(x_{12} y_{21} - y_{12} x_{21}). \] 
(cf.\ \cite[Rem.~IV.5]{Ne02}).  
Identifying a Hilbert--Schmidt operator 
$z \in B_2(\cH_-, \cH_+)$ with the corresponding matrix 
\[ \tilde z := \pmat{0 & z \\ z^* & 0} \in \fp,\]  
this leads to 
\[ \omega(\tilde z, \tilde w) = \tr(zw^*-wz^*) 
= 2i \Im \tr(zw^*).\] 
This means that $\omega = i\omega_\fp$ and thus 
\cite[Prop.~IV.8]{Ne02}
implies that $\hat\g$ is perfect and centrally closed, 
hence a universal central extension 
(Definition~\ref{def:h2}(e)). 

To obtain a (simply connected) Lie group $\hat G = 
\hat \U_{\rm res}(\cH_+, \cH_-)$ with Lie algebra 
$\hat\g$, we first consider the Banach--Lie algebra 
\[ \g_1 := 
\Big\{ \pmat{a & b \\ c & d} \in \g \: 
a \in B_1(\cH_+), d \in B_1(\cH_-) \Big\}\]  
(cf.\ Appendix~\ref{app:d.2}). Then 
\[ \sigma \: \g_1 \to \hat\g, \quad 
\sigma(x) := (-i\tr x_{11}, x) \] 
is a homomorphism of Lie algebras 
because the linear  functional 
$\alpha(x) := -i\tr(x_{11})$ on $\g_1$ 
satisfies 
$\alpha([x,y]) = -i\omega(x,y) = \omega_\fp(x,y)$ for $x,y \in \g_1$. 
As $\tr \: B_1(\cH_+) \to \C$ vanishes on  
$[B(\cH_+), B_1(\cH_+)]$, the homomorphism 
$\sigma$ can be combined with the inclusion 
$\fk \into \hat\g, x \mapsto (0,x)$ to a surjective homomorphism, 
\[ \g_1 \rtimes \fk \to \hat\g, \quad 
(x,y) \mapsto (-i\tr(x_{11}), x + y).\] 
Passing from $\g_1$ to 
\[ \sg_1 := \ker(\tr\res_{\g_1}) 
= \Big\{ \pmat{a & b \\ c & d} \in \g_1 \:  \tr(a) + \tr(d) = 0 \Big\},\] 
we also obtain a surjective homomorphism 
\[ \gamma_\g \: \sg_1 \rtimes \fk \to \hat\g, \quad 
(x,y) \mapsto (-i\tr(x_{11}), x + y).\] 

(b) For the positive $JH^*$-triple structure on $\fp = 
B_2(\cH_-, \cH_+)$, the situation 
on the group level looks as follows. 
The subgroup $K = \U(\cH_+) \times \U(\cH_-) 
\subeq \U_{\rm res}(\cH_+, \cH_-)$ is simply connected 
by Kuiper's Theorem, and the group 
\begin{equation}
  \label{eq:g1,2}
G_1 := \{ g \in G \: \|\1 - g_{11}\|_1 < \infty, \|\1 - g_{22}\|_1 < \infty\} 
\end{equation}
is a Banach--Lie group with Lie algebra 
$\g_1$ and a polar decomposition 
$G_1 = K_1 \exp(\fp)$, where 
\[ K_1 := K \cap G_1 \cong \U_1(\cH_+) \times \U_1(\cH_-)\]  
(cf.\ Remark~\ref{rem:polar}). 
For the subgroup $SG_1 := \ker(\det \: G_1 \to \T)$ we have
\[ SK_1 := K \cap SG_1 \cong 
\{ (k_1, k_2) \in  \U_1(\cH_+) \times \U_1(\cH_-) 
\: \det(k_1) \det(k_2) = 1 \}.\] 
Since $\tilde\U_1(\cH_\pm) \cong \R \ltimes \SU(\cH_\pm)$ 
holds for the universal covering group of $\U_1(\cH_\pm)$, 
we have the universal covering group 
\[ \tilde{SK}_1 \cong \R \ltimes (\SU(\cH_+) \times \SU(\cH_-)),\] 
which, according to the polar decomposition 
(cf.\ Remark~\ref{rem:polar}), 
is contained in the universal covering group 
$\tilde{SG}_1$ of $SG_1$. 
If $\hat G$ is a connected Lie group with Lie algebra $\hat\g$, 
the homomorphism $\gamma_\g$ integrates to a surjective 
homomorphism of Lie groups 
\[ \gamma_G \: \tilde{SG}_1 \rtimes K \to \hat G \] 
which restricts to 
\[ \gamma_G \: \tilde{SK}_1 \rtimes K \to \hat K \cong 
\R \times K, \quad ((t,k_1,k_2), (u_1, u_2)) \mapsto 
(t, k_1 u_1, k_2 u_2)\] 
(cf.\ \cite[Sect.~IV]{Ne02}, where this is used to construct 
$\hat G$ as a quotient group). 
The kernel of $\gamma_G$ is isomorphic to 
$\SU(\cH_+) \times \SU(\cH_-)$. 

(c) For the negative $JH^*$-triple structure on $\fp = B_2(\cH_-, \cH_+)$ 
and $G = \U_{\rm res}(\cH)$, 
we define $G_1$ and $SG_1$ as in \eqref{eq:g1,2}. 
Then the inclusion $\SU(\cH) \to SG_1$ is a weak 
homotopy equivalence, 
and in particular $SG_1$ is simply connected 
(use \cite[Prop.~III.2(c)]{Ne02} and the polar decomposition). 
We thus obtain a surjective homomorphism of Lie groups 
$\gamma_G \: SG_1 \rtimes K \to \hat G$  
which restricts to 
\[ \gamma_G \: SK_1 \rtimes K \to \hat K \cong 
\T \times K, \quad ((t,k_1,k_2), (u_1, u_2)) \mapsto 
(t, k_1 u_1, k_2 u_2).\] 

(d) We also note that the positive $JH^*$-triple structure 
$\{x,y,z\} = \shalf(xy^*z + zy^*x)$ on 
$U_+ := B_2(\cH_-, \cH_+)$ leads to the Lie algebra 
\[ \g(U_+) 
\cong \fu_{\rm res}(\cH_+, \cH_-)/i\R \1 
\cong \ad(\fu_{\rm res}(\cH_+, \cH_-)), \] 
and if both $\cH_\pm$ are infinite dimensional, then we obtain the 
non-trivial central extension 
\[ \hat\g(U_+) 
\cong \hat\fu_{\rm res}(\cH_+, \cH_-).\] 
Accordingly, the negative $JH^*$-triple structure 
$\{x,y,z\} = -\shalf(xy^*z + zy^*x)$ on 
$U_- := B_2(\cH_-,, \cH_+)$ leads to 
\[ \g(U_-) \cong \fu_{\rm res}(\cH)/i\R \1 
\cong \ad(\fu_{\rm res}(\cH)), \] 
and if both $\cH_\pm$ are infinite dimensional, we have the 
non-trivial central extension 
\[ \hat\g(U_-) \cong \hat\fu_{\rm res}(\cH).\] 
\end{ex}

\begin{ex}
\mlabel{ex:3.6}  {\bf Type II/III:} (a) For $\fk = \fu(\cH)$ and $\dim \cH = \infty$,  
Remark~\ref{rem:3.5}(b) 
implies $H^2_c(\fk,\R) = \{0\}$, and we obtain as above with 
Lemma~\ref{lem:cocyc} that 
$H^2_c(\sp_{\rm res}(\cH),\R) = \R [\omega_\fp]$ and likewise that 
$H^2_c(\fo^*_{\rm res}(\cH),\R) = \R [\omega_\fp]$. 
According to \cite[Prop.~I.11]{Ne02}, $[\omega_\fp] \not=0$ 
in both cases, and \cite[Prop.~IV.8]{Ne02} implies that 
$\hat\g$ is centrally closed. Note that 
$\dim \z(\hat\fk) = 2$ and $\dim \z(\hat\g) = 1$. 

(b) On the group level, we have for the positive $JH^*$-triple structure 
in a similar fashion as for type I a Lie algebra 
$\g_1 =  \fk_1 \oplus \fp$ with 
$\fk_1 \cong \fu_1(\cH)$ and a surjective homomorphism 
\[ \gamma \: \g_1 \rtimes \fk \to \hat\g, \quad  
(x,y) \mapsto (-i\tr(x_{11}), x + y). \] 
Again, 
$K = \U(\cH)$ is simply connected and the group 
\begin{equation}
  \label{eq:g1,2b}
G_1 := 
\{ g \in G \: \|\1 - g_{11}\|_1 < \infty, \|\1 - g_{22}\|_1 < \infty \} 
\end{equation}
is a Banach--Lie group with Lie algebra 
$\g_1$ and a polar decomposition 
$G_1 = K_1 \exp(\fp)$, where 
$K_1 := K \cap G_1 \cong \U_1(\cH).$ 
From the polar decomposition we obtain in particular 
$\pi_1(G_1) \cong \pi_1(K_1) \cong \Z$. Hence 
the universal covering group $\tilde G_1$ contains 
$\tilde K_1 \cong \tilde\U_1(\cH) \cong \R \ltimes \SU(\cH)$. 
Now $\gamma$ integrates to a surjective 
homomorphism of Lie groups \break 
$\gamma_G \: \tilde{G}_1 \rtimes K \to \hat G$ 
which restricts to 
\[ \gamma_G \: \tilde{K}_1 \rtimes K \to \hat K \cong 
\R \times K, \quad ((t,k), u) \mapsto 
(t, k u)\] 
(cf.\ \cite[Sect.~IV]{Ne02}). 
The kernel of $\gamma_G$ is isomorphic to $\SU(\cH)$. 

(c) For the negative $JH^*$-triples of type II/III we 
write 
\[ \beta_\pm((x,y), (x',y')) := \la x, \sigma(y')\ra 
\mp \la x', \sigma(y) \ra \] 
and obtain 
\[ G = \OO(\cH^2,\beta_+) \cap \U_{\rm res}(\cH^2), \quad  \mbox{ resp.,  } 
\quad 
\Sp(\cH^2,\beta_-) \cap \U_{\rm res}(\cH^2). \]  
In both cases 
$K = (G^\theta)_0 
= \{\diag(g,g^{-\top}) \: g \in \U(\cH) \} \cong \U(\cH).$ 
We define $G_1$ as in \eqref{eq:g1,2} and 
obtain $K_1 := G_1 \cap G \cong \U_1(\cH)$. 
From the polar decomposition and  
\cite[Prop.~III.15]{Ne02}, we derive that 
the inclusion $G_1 \to \OO_2(\cH^2,\beta_+)$, resp., 
$G_1 \to \Sp_2(\cH^2,\beta_-)$ are weak homotopy equivalences. 
This shows in particular that $\pi_1(G_1) \cong \Z/2\Z$ for type II 
and that $G_1$ is simply connected for type III (\cite[Thm.~II.14]{Ne02}). 
In both cases we have a surjective homomorphism of Lie groups 
$\gamma_G \: \tilde G_1 \rtimes K \to \hat G$, but 
the two types lead to different pictures for the subgroup 
$\hat K_1 = \la \exp_{\tilde G_1} \fk_1 \ra$. 
For type III we have $\tilde G_1 = G_1$ and $\hat K_1  = K_1$, 
but for type II, the group $\hat K_1$ is the unique $2$-fold 
covering of $K_1 \cong \T \ltimes \SU(\cH)$, i.e., 
$\hat K_1 \cong \T \ltimes \SU(\cH)$ with the universal covering 
map $q_{G_1} \: \tilde G_1 \to G_1$ satisfying 
$q_{G_1}(z,k) = (z^2, k) \in \T \ltimes \SU(\cH) \cong K_1$. 
Put differently, the projection 
$\sqrt{\det} \: \hat K_1 \to \T, (z,k) \mapsto z$ 
is a square root of the pullback of the determinant to $\hat K_1$ which 
satisfies $\L(\sqrt{\det}) = \shalf \tr$. 
\end{ex} 

\begin{ex}
  \mlabel{ex:3.7} {\bf Type IV:} 
Since $\g_\C = \fo(\R^2, \cH_\R)_\C 
\cong \fo(\C^2 \oplus \cH)$, with  $\cH = (\cH_\R)_\C$, is the full 
orthogonal Lie algebra of the complex Hilbert space $\C^2 \oplus \cH$, 
all its central extensions are trivial (\cite[Prop.~I.7]{Ne02}).
\end{ex}

\begin{rem} \mlabel{rem:3.5b} 
(a) From the preceding discussion, we 
conclude that if $\g$ is one of the hermitian Lie algebras 
\[  \fu_{\rm res}(\cH_+,\cH_-), \quad\sp_{\rm res}(\cH), \quad 
\fo^*_{\rm res}(\cH) \quad \mbox{ or } \quad \fo(\R^2, \cH_\R),\]  
the canonical central extension $\hat\g$ is either trivial 
(type I with $\cH_+$ or $\cH_-$ finite dimensional, and type IV), or 
$\hat\g$ is the universal central extension of $\g$ 
(cf.\ Definition~\ref{def:h2}(d)). 
The same holds for the $c$-dual 
Lie algebras 
\[  \fu_{\rm res}(\cH_+\oplus\cH_-),\quad \sp_{\rm res}(\cH^2_\H), \quad 
\fo_{\rm res}(\cH^\R) \quad \mbox{ and } \quad \fo(\R^2\oplus \cH_\R)\]  
(cf.\ Appendix~\ref{app:d.1b}).

(b) If $(\g,\theta, d)$ is a full hermitian Lie algebra, then 
$\ad \g \cong \g(\fp)$ holds for the corresponding $JH^*$-triple 
$\fp$. 
If $\hat\g(\fp)$ is the universal central extension of 
$\g(\fp)$ from (a), then the universal 
property (cf.\ Definition~\ref{def:h2}) 
implies the existence of a morphism 
$\alpha \: \hat\g(\fp) \to \g$ of central extensions of $\g(\fp)$. 
Then $\alpha$ maps the finite dimensional center 
$\z(\hat\g(\fp))$ into $\z(\g)$, and with any decomposition 
$\z(\g) = \z \oplus \alpha(\z(\hat\g(\fp)))$, 
we obtain a direct decomposition 
$\g \cong \z \oplus \alpha(\hat\g(\fp))$, where
$\alpha(\hat\g(\fp)) = [\g,\g]$ is the commutator algebra. 
\end{rem}

For the classification of irreducible semibounded representations 
it therefore suffices to consider the simply connected Lie group 
$\hat G(\fp)$ with Lie algebra $\hat\g(\fp)$. This is the 
approach we shall follow in Sections~\ref{sec:7} and \ref{sec:8}. 
The existence of $\hat G(\fp)$ is a consequence of the following 
theorem.

\begin{thm} \mlabel{thm:exist} 
If $(G,\theta, d)$ is full, then the 
Lie algebra $\hat\g = \R \oplus_{\omega_\fp}\g$ 
is integrable, i.e., there exists a Lie group $\hat G$ with 
Lie algebra $\hat\g$. 

Moreover, $\hat\g(\fp)_\C$ is integrable for every simple 
$JH^*$-triple $\fp$. 
\end{thm}

\begin{prf} The Lie algebra $\hat\g$ is clearly integrable if 
$[\omega_\fp]=0$, because this implies that 
$\hat\g \cong \R \oplus \g$. 

Since the adjoint representation yields a morphism 
$\Ad\: G \to G(\fp)$ of the simply connected Lie groups with 
the Lie algebras $\g$, resp., $\g(\fp)$, it suffices to show that 
$\g^\sharp := \R \oplus_{\omega_\fp} \g(\fp)$ is integrable. 
Then the pullback of the corresponding central 
Lie group extension $G^\sharp$ of $G(\fp)$ by $\Ad$ is a central extension 
of $G$ with Lie algebra~$\hat\g$. 

We recall that the cocycle $\omega_\fp$ is non-trivial for $\g(\fp)$ 
only if the $JH^*$-triple $\fp$ is of type I with $\dim \cH_\pm= \infty$ 
or of type II or III (Remark~\ref{rem:3.5b}). 
We also know from \cite[Def.~IV.4, Prop.~IV.8]{Ne02} that the 
universal central extension $\hat\g(\fp)$ of 
$\g(\fp)$ and its complexification are integrable. 
For Type II/III, the assertion follows from 
$\hat\g(\fp) \cong \g^\sharp$. 

For type I with either $\cH_+$ or $\cH_-$ finite dimensional,  we 
have for $\fp$ positive $\hat G(\fp) 
\cong \tilde\U_{\rm res}(\cH_+, \cH_-)
= \tilde\U(\cH_+, \cH_-)$ with 
$\hat G(\fp)_\C \cong \tilde\GL(\cH_+\oplus \cH_-)$, 
and for type IV we have for $\fp$ positive 
$\hat G(\fp)\cong \tilde\OO(\R^2,\cH_\R)$ with 
$\hat G(\fp)_\C\cong \tilde\OO(\C^2 \oplus \cH)$ 
(Example~\ref{ex:3.7}). 

For Type I with $\cH_\pm$ both infinite dimensional and 
$\fp$ positive, the identity component of the 
center of the simply connected Lie group 
$\hat G(\fp)$ with Lie algebra $\hat\g(\fp)$ is isomorphic to 
\[ \R \times \T \cong 
\R \times \{ (t\1,t^{-1}\1) \: t \in \T\}
\subeq \R \times \U(\cH_+) \times \U(\cH_-), \] 
where factorization of the circle subgroup yields a 
group with Lie algebra $\R \oplus_{\omega_\fp} \g(\fp)$. 
If $\fp$ is negative, then 
\[ Z(\hat G(\fp))_0 \cong \T^2
\cong \T \times \{ (t\1,t^{-1}\1) \: t \in \T\}
\subeq \T \times \U(\cH_+) \times \U(\cH_-), \] 
where factorization of the second  $\T$-factor yields a 
group with Lie algebra $\R \oplus_{\omega_\fp} \g(\fp)$ 
(cf.\ the discussion of the corresponding groups in Remark~\ref{rem:2.7}). 
The construction of the complex Lie group 
$\hat\GL_{\rm res}(\cH_+,\cH_-) \cong \hat\U_{\rm res}(\cH_+,\cH_-)_\C$ 
can be found in \cite[Def.~IV.4]{Ne02}. 
\end{prf}

\begin{rem} (a) If $[\fp,\fp] = \{0\}$, 
then the assertion of the preceding theorem is true for trivial 
reasons. In this case $\g \cong \fp \rtimes \fk$ leads to 
$\hat\g \cong \hat\fp \rtimes \fk$, where 
$\hat\fp = \R \oplus_{\omega_\fp} \fp$ is a Heisenberg algebra 
(cf.\ Example~\ref{ex:1.3}(b)). Therefore 
$\hat G := \Heis(\fp) \rtimes K$ is a corresponding Lie group. 

(b) The Lie algebra cocycle 
$\omega_\fp$ defines a $G$-invariant 
non-degenerate closed $2$-form $\Omega$ on $M$, for which 
$\Omega_{\1K}$ corresponds to $\omega_\fp$ under the natural identification 
$T_{\1K} \cong \fp$. Here the closedness of this $2$-form 
follows from the closedness of its pullback $q^*\Omega$ to $G$ under the 
quotient map $G \to G/K = M$ which is a 
left invariant $2$-form on $G$ with $(q^*\Omega)_\1 = \omega_\fp$ 
(cf.~\cite[Thm.~10.1]{CE48}). 

On the central extension $\hat\g = \R \oplus_{\omega_\fp} \g$,  
the coadjoint action of $\hat G$ factors through an action 
of $G$, and the linear functional $\lambda(t,x) := t$ satisfies 
$G_\lambda = K$, so that its orbit 
$\cO_\lambda = \Ad^*(G)\lambda$ 
is isomorphic to $G/K$. Therefore 
the canonical action of $\hat G$ on the symplectic 
K\"ahler manifold $(M,\Omega)$ is hamiltonian. 
\end{rem}

\begin{rem} \mlabel{rem:3.9} 
Suppose that $\fp$ is a simple infinite dimensional 
$JH^*$-triple and let $\g = \fk \oplus \fp = \hat\g(\fp)$ be 
the universal central extension of $\g(\fp)$ 
(Remark~\ref{rem:3.5b}).  
We now exhibit an inclusion of a Banach-Lie algebra 
$\fk_1 \into \fk$ such that 
$\g_1 := \fk_1 \oplus \fp$ carries a natural Banach--Lie algebra 
structure with a continuous inclusion 
$\g_1 \into \g$ such that the summation map 
$\g_1 \rtimes \fk \to \g, (x,y) \mapsto x + y$ 
is a quotient homomorphism. We have the following 
cases: 

\begin{description}
\item[\rm(I$_\infty$)] $\fp = B_2(\cH_-,\cH_+)$ with 
both $\cH_\pm$ infinite dimensional, 
\[ \fk = \R \oplus \fu(\cH_+) \oplus \fu(\cH_-),  \qquad 
\fk_1 := \fu_1(\cH_+) \oplus \fu_1(\cH_-), \] and 
$\g_1 := \fk_1 \oplus \fp \subeq \gl(\cH)$ with 
the inclusion $\g_1 \into \g, x \mapsto (-i\tr x_{11}, x)$. 
\item[\rm(I$_{\rm fin}$)] $\fp = B_2(\cH_-,\cH_+)$ with 
$\cH_+$ or $\cH_-$ finite dimensional. Then 
$\fk = \fu(\cH_+) \oplus \fu(\cH_-)$, 
$\fk_1 := \fu_1(\cH_+) \oplus \fu_1(\cH_-)$ and 
$\g_1 := \fk_1 \oplus \fp \subeq \gl(\cH)$ 
with the canonical inclusion $\g_1 \into \g \subeq \g(\cH)$.
\item[\rm(II)] $\fp = \Skew_2(\cH) = \{ A \in B_2(\cH) 
\: A^\top= - A\}$. Then 
$\fk \cong \R \oplus \fu(\cH)$, 
$\fk_1 := \fu_1(\cH)$ and 
$\g_1 := \fk_1 \oplus \fp \subeq \gl(\cH^2)$ 
with the inclusion $\g_1 \into \g$ given by 
$x \mapsto (-i\tr x_{11}, x)$. 
\item[\rm(III)] $\fp = \Sym_2(\cH) = \{ A \in B_2(\cH) 
\: A^\top=  A\}$. Then 
$\fk \cong \R \oplus \fu(\cH)$, 
$\fk_1 := \fu_1(\cH)$ and 
$\g_1 := \fk_1 \oplus \fp \subeq \gl(\cH^2)$ 
with the inclusion $\g_1 \into \g$ given by 
$x \mapsto (-i\tr x_{11}, x)$. 
\item[\rm(IV)] $\fp = \cH_\R^2\cong \R^2 \otimes \cH_\R$ 
with a real Hilbert space 
$\cH_\R$. Then 
$\fk = \so_2(\R) \oplus \fo(\cH_\R)$, 
$\fk_1 := \so_2(\R) \oplus \fo_1(\cH_\R)$ and 
$\g_1 := \fk_1 \oplus \fp \subeq \gl(\R^2 \oplus \cH_\R)$ 
with the canonical inclusion $\g_1 \into \g$. 
\end{description}
\end{rem}

\section{Open invariant cones in hermitian Lie algebras} 
\mlabel{sec:4}

If each open convex invariant cone $W$ in a Banach--Lie algebra 
$\g$ is trivial, then this holds in particular for the cone 
$W_\pi$ of any semibounded unitary representation. As a consequence, 
every semibounded representation is bounded. Therefore 
we study in this section criteria for the triviality of open 
invariant cones in Banach--Lie algebras.

\subsection{Lie algebras with no open invariant cones} 

\begin{defn} We say that the Banach--Lie algebra $\fk$ has 
{\it no open invariant cones} if each non-empty open invariant 
convex cone $W \subeq\fk$ coincides with $\fk$. 
\end{defn}

\begin{lem} \mlabel{lem:conext} 
{\rm(a)} If $\fn \trile \fk$ is a closed ideal, then 
$\fk/\fn$ has no open invariant cones if and only if 
all open invariant cones in $\fk$ intersect~$\fn$.  

{\rm(b)} If $\fn \trile \fk$ is a closed ideal 
such that neither $\fn$ nor $\fk/\fn$ have open invariant cones, 
then $\fk$ has no open invariant cones. 

{\rm(c)} If $\fk = \bigoplus_{i = 1}^n \fk_i$ is a direct sum of Lie algebras 
with no open invariant cones, then $\fk$ has no open invariant cones. 
\end{lem}

\begin{prf} (a) Let $q \: \fk \to \fk/\fn$ denote the quotient map. 
If $W \subeq \fk$ is an open invariant cone intersecting 
$\fn$ trivially, then $W + \fn = q^{-1}(q(W))$ 
is a proper open invariant cone in 
$\fk$, hence $q(W)$ is a proper open invariant cone in $\fk/\fn$. 

Now suppose that $\fk/\fn$ has no open invariant cones and 
let $W \subeq \fk$ be an open invariant cone. Since 
$q(W)$ is an open invariant cone in $\fk/\fn$, it contains $0$, 
which means that $W \cap \fn \not=\eset$. 

(b) Let $\eset\not=W \subeq \fk$ be an open invariant cone. 
From (a) we derive that $W \cap \fn$ is a non-empty open 
invariant cone in $\fn$. As $\fn$ contains no open invariant cones, 
we obtain $0 \in W \cap \fn \subeq W$, and thus $W = \fk$. 

(c) follows from (b) by induction. 
\end{prf}

\begin{ex} \mlabel{ex:cpt} 
If $K$ is a connected Lie group with compact Lie algebra, 
then $\Ad(K)$ is compact, so that averaging with respect to Haar measure 
implies that every open invariant cone $\eset\not=W \subeq \fk$ intersects 
$\z(\fk)$. In particular, $\fk/\fz(\fk)$ has no open invariant cones because 
it is semisimple. 
\end{ex}

Below we shall see that there is an interesting family of infinite dimensional 
Banach--Lie algebras behaving very much like compact ones. 
An important point of the triviality of all open invariant cones is 
that this property permits us to draw conclusions 
concerning boundedness of unitary representations, such as the following. 
For the definition of boundedness and semiboundedness of a unitary 
representation, we refer to the introduction. 

\begin{prop} \mlabel{prop:4.5} 
Let $K$ be a Lie group for which $\fk/\z(\fk)$ has no 
open invariant cone. If $(\rho, V)$ is a semibounded unitary representation 
of $K$ for which $\rho\res_{Z(K)_0}$ is bounded, then $(\rho,V)$ is bounded. 
In particular, every irreducible semibounded representation of $K$ is bounded. 
\end{prop}

\begin{prf} Since $\rho$ is semibounded, the open cone $W_\rho \subeq \fk$ 
is non-trivial, and Lemma~\ref{lem:conext}(a) implies that 
$W_\rho \cap \z(\fk)\not=\eset$. The assumption that $\rho\res_{Z(K)_0}$ 
is bounded further implies that $\z(\fk) + W_\rho = W_\rho$, which 
leads to $0 \in W_\rho$, hence to $W_\rho = \fk$, i.e., $\rho$ is bounded. 

If, in addition, $\rho$ is irreducible, then 
$\rho(Z(K)) \subeq \T \1$ follows from Schur's Lemma, and the preceding 
argument implies that $\rho$ is bounded. 
\end{prf}

\subsection{Open invariant cones in unitary Lie algebras} 

In this section we show that, for any real, complex or quaternionic 
Hilbert space $\cH$, the Lie algebra $\fu(\cH)/\fz(\fu(\cH))$ has 
no open invariant cones. Note that the center $\z(\fu(\cH))$ 
is only non-trivial for complex Hilbert spaces, where it is one-dimensional, 
and for $2$-dimensional real Hilbert spaces.  

\begin{lem} \mlabel{lem:5.1a} 
Let $J$ be a set and $S_{J}$ be the group of all permutations 
of $J$ acting canonically on $\ell^\infty(J,\R)$. 
Then every $S_{J}$-invariant non-empty open convex cone 
in $\ell^\infty(J,\R)$ contains a constant function. 
\end{lem}

\begin{prf} {\bf Step 1:} If $J$ is finite, then $S_{J}$ is finite, so that the 
assertion follows by averaging. We may therefore assume that 
$J$ is infinite. Let $\eset\not= W \subeq \ell^\infty(J,\R)$ be 
an open $S_{J}$-invariant convex cone. Then $W$ contains a constant 
function if and only if $0 \in W + \R \1$, where 
$\1$ stands for the constant function $1$. We may therefore 
assume that $W = W + \R \1$ and show that $0 \in W$. 

 {\bf Step 2:} Since the functions 
$x$ with only finitely many values are dense in $\ell^\infty(J,\R)$, 
the cone $W$ contains such an element. 
Let $J = J_0 \cup \cdots \cup J_N$ be the corresponding decomposition 
of $J$ for which $x\res_{J_k}$ is constant $x_k$ for 
$k =1,\ldots, N$. We write this as 
\[ x = \sum_{k = 0}^N x_k \chi_{J_k},\] 
where $\chi_{J_k}$ is the characteristic function of $J_k$. 
Since  $J$ is infinite, there exists a $k$ with 
$|J_k| = |J|$. We may w.l.o.g.\ assume that 
$k = 0$ and observe that $y := x - x_0 \1 \in W$.

Fix $k > 0$ and put $J_k' := J \setminus (J_0 \cup J_k)$. 
As $|J_0| = |J_0 \times \N|$, 
there exist pairwise disjoint subsets 
$J_k^n\subeq J_0$, $n \in \N$, and involutions 
$\sigma_n \in S_J$ with 
\[\sigma_n(J_k') = J_k^n, \quad \sigma_n(J_k^n) = J_k' \quad \mbox{ and } 
\quad \sigma_n\res_{J_k \cup J_0\setminus J_k^n} = \id.\] 
Then $y^M := \frac{1}{M} \sum_{n = 1}^M y \circ \sigma_n$ 
coincides with $y$ on $J_k$, and 
on $J_0 \cup J_k'$ it is bounded by $\frac{1}{M}\|y\|_\infty$. 
Therefore 
$y_k \chi_{J_k} = \lim_{M \to \infty} y^M \in \oline W$ 
holds for $k =1,\ldots, N$. 

 {\bf Step 3:} Let $M \subeq J$ be any subset with 
$|J\setminus M| = |J|$ and 
observe that there exists a subset $M' \subeq J_0$ with the same 
cardinality for which $J_0 \setminus M'$ still has the same cardinality 
as $J$. Since $W$ is open, there exists an $\eps > 0$ with 
$z := y \pm \eps \chi_{M'} \in  W$  and $\eps < |y_j|$ whenever $y_j\not=0$. 
Then the argument from Step $2$ applies to $z$ instead of $y$, 
which leads to 
$\pm \eps \chi_{M'} \in \oline W$. We conclude that 
$\chi_{M'} \in H(\oline W) = H(W)$. 
As 
\[ (|M|, |J\setminus M|)  
= (|M'|, |J\setminus M'|),\] 
there exists a permutation $\sigma \in S_J$ with 
$\sigma(M) = M'$, so that 
we also obtain $\chi_M \in H(W)$. 

{\bf Step 4:} From $|J \setminus J_k| \geq |J_0| = |J|$, 
we further derive that 
$\chi_{J_k} \in H(W)$ for $k =1,\ldots, N$, and this eventually leads to 
$0 = y - \sum_{k=1}^N y_k \chi_{J_k} \in W$. 
\end{prf}

\begin{thm} \mlabel{thm:coneinuh} If $\cH$ is a complex Hilbert space, then 
the Banach--Lie algebra $\pu(\cH) = \fu(\cH)/\R i \1$ 
has no non-trivial open invariant cones, i.e., 
each non-empty open invariant cone in $\fu(\cH)$ intersects $\R i \1$. 
\end{thm}

\begin{prf} If $\dim \cH < \infty$, then the assertion follows from the 
fact that $\pu(\cH)$ is a compact Lie algebra with trivial center. 
We may therefore assume that $\dim \cH = \infty$ and that 
$W \subeq \fu(\cH)$ is an open invariant convex cone. 

Since the elements with finite spectrum are dense in $\fu(\cH)$ 
by the Spectral Theorem, $W$ contains such an element $X$. 
Let $(e_j)_{j \in J}$ be an orthonormal basis of $\cH$ 
consisting of eigenvectors for $X$. 
Then the space $\ft$ of all diagonal operators in $\fu(\cH)$ 
with respect to $(e_j)$ is isomorphic to $\ell^\infty(J,i\R)$, 
and since any permutation $\sigma$ of the basis vectors 
corresponds to an element $u_\sigma \in \U(\cH)$, the 
invariance of $W$ implies that the cone 
$W \cap \ft \subeq \ft \cong \ell^\infty(J,i\R)$ is invariant 
under the action of $S_J$. It is non-empty because it contains 
$X$. Therefore Lemma~\ref{lem:5.1a} implies that 
$W\cap i\R \1\not=\eset$. 
\end{prf}

\begin{thm} \mlabel{thm:coneinoh} 
If $\cH_\R$ is a real Hilbert space with $\dim \cH_\R > 2$, then 
the Banach--Lie algebra $\fo(\cH_\R)$ has no non-trivial open 
invariant cones. 
\end{thm}

\begin{prf} If $2 < n := \dim \cH_\R < \infty$, then 
$\fo(\cH_\R) \cong \so_n(\R)$ is a compact semisimple Lie algebra. 
Hence every open invariant convex cone in $\fo(\cH_\R)$ contains the only 
fixed point $0$ of the adjoint action, hence cannot be proper. 

We now assume that $\dim \cH_\R = \infty$ 
and choose an orthogonal complex structure~$I$ on $\cH_\R$ 
(cf.\ Example~\ref{ex:d.1c}). 
This complex structure defines on $\cH_\R$ the structure of a complex 
Hilbert space $\cH := (\cH_\R,I)$. 

Let $W \subeq \fo(\cH_\R)$ be a non-empty open convex invariant 
cone and $p \: \fo(\cH_\R) \to \fu(\cH)$ be the fixed point 
projection for the circle group $e^{\R \ad I}$. 
Then $p(W) = W \cap \fu(\cH)$  is a non-empty open 
invariant cone in $\fu(\cH)$ (\cite[Prop.~2.11]{Ne10c}), so that 
Theorem~\ref{thm:coneinuh} implies that 
$\lambda I \in W$ holds for some $\lambda \in \R$. 

Let $J \: \cH \to \cH$ be an anticonjugation 
(cf.\ Appendix~\ref{app:d.1}). Then $J\in \fo(\cH_\R)$ is a second 
complex structure anticommuting with $I$. This leads to 
$[J,I] = JI-IJ = 2 JI$, and hence to 
$e^{\frac{\pi}{2}\ad J}I = e^{\pi J}I = -I$. 
This shows that $\pm I \in W$, which leads to 
$0 \in W$, and finally to $W = \fo(\cH_\R)$. 
\end{prf}

\begin{thm} \mlabel{thm:coneinquat} 
If $\cH$ is a quaternionic Hilbert space, then 
the Banach--Lie algebra $\sp(\cH) = \fu_\H(\cH)$ has no non-trivial open 
invariant cones. 
\end{thm}

\begin{prf} If $\dim \cH < \infty$, then 
$\sp(\cH)$ is a compact simple Lie algebra, so that 
every open invariant cone contains $0$, hence cannot be proper. 

Let $\cH^\C$ denote the complex Hilbert space underlying 
$\cH$ and note that $\cH^\C$ can be written as 
$\cK^2$, where $\cK$ is a complex Hilbert space and 
the quaternionic structure on 
$\cH^\C$ is given in terms of a conjugation 
$\tau$ on $\cK$ by the anticonjugation 
$\sigma(v,w) := (-\tau w, \tau v)$ on $\cK^2$. 
Then 
\[ \sp(\cH) = \{ x \in \fu(\cK^2) \: \sigma x = x \sigma \}\] 
contains the operator $d(v,w) := (iv,-iw)$ whose centralizer 
is 
\[ \{ x \in \sp(\cH) \subeq \gl(\cK^2) \: x_{12} = x_{21} = 0\} 
= \{ \diag(x,-x^\top) \: x \in \fu(\cK)\} \cong \fu(\cK).\] 

Let $W \subeq \sp(\cH)$ be a non-empty open convex invariant 
cone and $p \: \sp(\cH) \to \fu(\cK)$ be the fixed point 
projection for the circle group $e^{\R \ad d}$. 
Then $p(W) = W \cap \fu(\cH,I)$ is a non-empty open 
invariant cone in $\fu(\cK)$ (\cite[Prop.~2.11]{Ne10c}), so that 
Theorem~\ref{thm:coneinuh} implies that 
$\lambda d \in W$ holds for some $\lambda \in \R$. 

The Lie algebra $\sp(\cH)$ also contains 
the operator $K(v,w) := (-w,v)$ because it commutes with~$\sigma$. 
Hence $W$ is also invariant under $e^{\R \ad K}$. 
From $K d = -dK$ it follows that 
\[ [K,[K,d]] = - 2[K,dK]= -2[K,d]K = 4 dK^2 = - 4 d,\] 
so that $e^{\frac{\pi}{2} \ad K}d = -d$ leads to $0 \in W$, 
and finally to $W = \sp(\cH)$. 
\end{prf}

\subsection{Applications to hermitian Lie algebras} 

\begin{lem} \mlabel{lem:12.1} 
Let $\fk$ be a Banach--Lie algebra for which 
$\fk/\z(\fk) = \oplus_{j = 1}^n \fk_j$, where 
each $\fk_j$ is either compact semisimple or isomorphic to 
$\fu(\cH)/\z(\fu(\cH))$ for a real, complex or quaternionic 
Hilbert space $\cH$. Then every non-empty open invariant convex cone 
$W \subeq \fk$ intersects $\z(\fk)$. 
\end{lem}

\begin{prf} In view of Theorems~\ref{thm:coneinuh} and 
\ref{thm:coneinoh}, \ref{thm:coneinquat} and Example~\ref{ex:cpt}, 
Lemma~\ref{lem:conext}(c) implies that $\fk/\z(\fk)$ contains 
no open invariant cones. Hence the assertion follows from 
Lemma~\ref{lem:conext}(a). 
\end{prf} 

\begin{lem} \mlabel{lem:12.2} 
If $(\g,\theta,d)$ is a full irreducible hermitian Lie algebra, 
then $\fk/\z(\fk)$ contains no open invariant cones. 
\end{lem} 

\begin{prf} In view of Lemma~\ref{lem:12.1}, we only have to show that 
$\fk/\z(\fk) \cong \aut(\fp)/\z(\aut(\fp))$ (cf.\ \eqref{eq:transfer}) 
has the required structure. This is an easy consequence of 
Theorem~\ref{thm:classif}: 

{\bf Type I:} $\aut(\fp)  \cong (\fu(\cH_+) \oplus 
\fu(\cH_-))/\R i(\1,\1)$ implies 
$\fk/\z(\fk) \cong  \pu(\cH_+) \oplus \pu(\cH_-).$

{\bf Type II/III:} $\fk/\z(\fk) \cong \pu(\cH)$. 

{\bf Type IV:} $\fk/\z(\fk) \cong \fo(\cH_\R)$. 
\end{prf}

\begin{prop} \mlabel{prop:12.2} 
If $(\g,\theta,d)$ is a hermitian Lie algebra for which 
$\fk/\z(\fk)$ contains no open invariant cones, then 
every open invariant cone $\eset \not=W \subeq \g$ intersects $\z(\fk)$. 
\end{prop} 

\begin{prf} Let $p_\fk \: \g \to \fk$ denote the projection 
along $\fp$. This is the fixed point projection with respect to the 
action of the circle group $e^{\R \ad d}\subeq \Aut(\g)$, so that 
\[ p_\fk(x) = \frac{1}{2\pi} \int_0^{2\pi} e^{t \ad d}x\, dt \] 
implies that $p_\fk(W) = W \cap \fk$ (cf.~\cite[Prop.~2.11]{Ne10c}). 
We conclude that $W \cap \fk$ is a non-empty open invariant 
cone in $\fk$ which intersects $\z(\fk)$ by Lemma~\ref{lem:conext}(a). 
\end{prf}

\section{Semiboundedness and positive energy} \mlabel{sec:5} 

In this section we start with our analysis of semibounded unitary 
representations of hermitian Lie groups. 
If $(G,\theta, d)$ is irreducible, we first show that 
each irreducible semibounded representation $(\pi, \cH)$ satisfies 
$d \in W_\pi \cup - W_\pi$. If $d \in W_\pi$, then we call 
$(\pi, \cH)$ a {\it positive energy representation}. In this 
case the maximal spectral value of the essentially selfadjoint operator 
$i\dd\pi(d)$ is an eigenvalue and the $K$-representation 
$(\rho,V)$ on the corresponding eigenspace is bounded and 
irreducible. Since $(\pi, \cH)$ is uniquely determined by $(\rho,V)$, 
we call a bounded representation $(\rho,V)$ of $K$ 
{\it (holomorphically) inducible} 
if it corresponds in this way 
to a unitary representation $(\pi, \cH)$ of~$G$. 
The classification of the holomorphically inducible representations 
$(\rho,V)$ is carried out for flat, negatively and positively 
curved spaces $G/K$ in Sections~\ref{sec:6}-\ref{sec:8}.

\subsection{From semiboundedness to positive energy} 

From the introduction we recall the concept 
of a semibounded unitary representation and the 
corresponding open invariant cone $W_\pi$.

\begin{defn} Let $(G,\theta, d)$ be a connected hermitian 
Banach--Lie group.  

(a) A unitary representation $(\pi, \cH)$ of 
a hermitian Lie group $G$ is called  a {\it positive energy representation} 
if $\sup \Spec(i\dd\pi(d)) < \infty$, i.e., 
if the infinitesimal generator $-i\dd\pi(d)$ of the one-parameter group 
$\pi(\exp td)$ is bounded from below. Note that this requirement 
is weaker than the condition $d \in W_\pi$. 

(b) We recall that $\fk_\C + \fp^-$ is a closed subalgebra of $\g_\C$ 
defining a complex structure on the homogeneous space 
$M := G/K$. Let $(\rho,V)$ be a bounded representation of $K$ and 
define a bounded Lie algebra representation 
$\beta \: \fp^+ \rtimes \fk_\C \to B(V)$ by 
$\beta(\fp^+) = \{0\}$ and $\beta\res_\fk = \dd\rho$ 
(cf.\ \eqref{eq:tridec}). Then the associated Hilbert bundle 
$\bV = G \times_K V$ over $G/K$ 
carries the structure of a holomorphic bundle on which $G$ 
acts by holomorphic bundle automorphisms (\cite[Thm.~I.6]{Ne10c}). 
A unitary representation $(\pi, \cH)$ of $G$ is said to be 
{\it holomorphically induced from $(\rho,V)$} 
if there exists a $G$-equivariant realization 
$\Psi \: \cH \to \Gamma(\bV)$ as a Hilbert 
space of holomorphic sections such that the 
evaluation $\ev_{\1 K} \: \cH \to V = \V_{\1 K}$ 
defines an isometric embedding $\ev_{\1 K}^* \: V \into \cH$. 
If a unitary representation $(\pi, \cH)$ holomorphically induced 
from $(\rho,V)$ exists, then it is uniquely determined 
(\cite[Def.~2.10]{Ne10d}) and we call $(\rho,V)$ {\it (holomorphically)  
inducible}. 

This concept of inducibility involves a choice of sign. 
Replacing $d$ by $-d$ changes the complex structure on $G/K$ and 
exchanges $\fp^+$ with $\fp^-$. 
\end{defn}

\begin{thm} \mlabel{thm:posen} 
Let $(\pi, \cH)$ be a semibounded representation 
of the hermitian Lie group $G$ for which $\fk/\z(\fk)$ contains no 
open invariant cones. Then the following assertions 
hold: 
\begin{description}
\item[\rm(i)] $\pi\res_{Z(K)_0}$ is semibounded and $W_\pi \cap \fz(\fk) \not=\eset$. 
\item[\rm(ii)] $\pi$ is a direct sum of representations which are bounded 
on $Z(G)$. 
\item[\rm(iii)] If $\pi$ is bounded on $Z(G)_0$ and $(G,\theta,d)$ is 
irreducible, then 
$d \in W_\pi \cup -W_\pi$. If $d \in W_\pi$, then 
$\pi$ is a positive energy representation, and if $d \in - W_\pi$, 
then  the dual representation $\pi^*$ is. 
\item[\rm(iv)] If $(G,\theta,d)$ is 
irreducible and $\pi$ is irreducible, then 
$d \in W_\pi \cup - W_\pi$. 
\end{description}
\end{thm}

\begin{prf} (i) The representation $\zeta := \pi\res_{Z(K)_0}$ satisfies 
$W_\zeta \supeq W_\pi \cap \z(\fk)$, and 
Proposition~\ref{prop:12.2} implies that the latter open cone 
is non-empty. 

(ii) From (i) it follows that $\pi\res_{Z(K)_0}$ is semibounded, 
hence given by a regular spectral measure on 
the dual space $\z(\fk)'$ (\cite[Thm.~4.1]{Ne09}). 
Now the regularity of the spectral measure implies (ii). 

(iii) If $\pi$ is bounded on $Z(G)_0$, then 
$\z(\g) + W_\pi = W_\pi$, so that 
\[ \eset\not = W_\pi \cap \z(\fk) = \z(\g) \oplus (W_\pi \cap \R d)\] 
follows from $\z(\fk) =\R d \oplus \z(\g)$ (Remark~\ref{rem:1.2}(d)). 
This proves (iii). 
If $d \in W_\pi$, then $\pi$ is a positive energy representation 
and otherwise $\pi^*$ with $W_{\pi^*} = - W_\pi$ has this property.

(iv) If $\pi$ is irreducible, then 
$\pi(Z(G)) \subeq \T \1$ follows from Schur's Lemma. 
In particular, $\pi\res_{Z(G)}$ is bounded, so that (iii) applies. 
\end{prf}


\begin{lem} \mlabel{lem:6.2b}
  Suppose that 
$(G,\theta, d)$ is a hermitian Lie group 
and that $(\pi, \cH)$ is a smooth unitary representation of~$G$. 
Then the following assertions hold: 
\begin{description}
\item[\rm(a)] $\cH$ decomposes into $Z(K)_0$-eigenspaces, which coincide with 
the $\exp(\R d)$-eigenspaces 
\[ \cH_\lambda = \{ v \in \cH \: (\forall t \in \R)\, \pi(\exp td)v 
= e^{it\lambda}v\},\quad \lambda \in \R.\] 
\item[\rm(b)] The corresponding 
projections $P_\lambda\: \cH \to \cH_\lambda$ preserve the space
$\cH^\infty$ of smooth vectors. In particular,  
$\cH^\infty \cap \cH_\lambda\supeq P_\lambda(\cH^\infty)$ is 
dense in $\cH_\lambda$ for every~$\lambda$.  
\item[\rm(c)] If $\lambda := \sup\Spec(i\dd\pi(d)) < \infty$, i.e., 
$(\pi, \cH)$ is a positive energy representation, 
then $\{0\} \not= \cH^\infty \cap \cH_\lambda\subeq 
(\cH^\infty)^{\fp^-}$. 
\end{description}
\end{lem}

\begin{prf} (a) Since $\pi$ is irreducible, $\pi(Z(G)) \subeq \T \1$ follows 
from Schur's Lemma. We consider the representation of the abelian group 
\[ Z(K)_0 = \exp(\z(\fk)) = Z(G)_0 \exp(\R d)\] 
on $\cH$ (Remark~\ref{rem:zent}). 
From $\pi(Z(G)) \subeq \T\1$ and 
$\exp(2\pi d) \in Z(G)$, we derive that 
\[ \pi(Z(K)_0) 
= \pi(Z(G)_0) \pi(\exp([0,2\pi]d))
\subeq \T \pi(\exp([0,2\pi]d)) \] 
is relatively compact in the strong operator topology, so that 
$\cH$ is an orthogonal direct sum of $Z(K)_0$-eigenspaces 
which coincide with the $\exp(\R d)$-eigenspaces. 
Suppose that $\pi(\exp 2 \pi d) = e^{2\pi i\mu_0} \1$, so that 
$\mu - \mu_0 \in \Z$ whenever $\cH_\mu\not=\{0\}$. 
Then 
\[ P_\mu v := \frac{1}{2\pi} \int_0^{2\pi} 
e^{-it \mu} \pi(\exp t d)v\, dt \] 
is the orthogonal projection onto $\cH_\mu$. 

(b) All these projections preserve $\cH^\infty$ 
(\cite[Prop.~3.4(i)]{Ne10d}), so that 
$P_\mu(\cH^\infty) = \cH^\infty \cap \cH_\mu$ 
is dense in $\cH_\mu = P_\mu(\cH)$. 

(c) If $\lambda := \sup(\Spec(i\dd\pi(d)))$ is finite 
and $v \in \cH_\lambda$ any smooth vector, then 
$\dd\pi(\fp^-)v \in  \cH_{\lambda+1} = \{0\}$. 
\end{prf}

\begin{thm} \mlabel{thm:6.2b} Suppose that 
$(G,\theta, d)$ is a hermitian Lie group 
for which $\fk/\z(\fk)$ contains no open invariant cones. 
For any irreducible 
semibounded positive energy representation $(\pi, \cH)$ of $G$,  
the $K$-representation $\rho$ on $V := \oline{(\cH^\infty)^{\fp^-}}$ 
is bounded and irreducible, and $(\pi,\cH)$ is holomorphically induced from 
$(\rho,V)$. 
\end{thm} 

\begin{prf} The smoothness of $\rho$ follows from 
the density of $V \cap \cH^\infty$ in~$V$. 
Since $\pi$ is semibounded, the corresponding cone $W_\pi$ 
is non-empty, and we know from 
Proposition~\ref{prop:12.2} that $W_\pi \cap \z(\fk) \not=\eset$. 
We conclude that the open cone $W_\rho \supeq W_\pi \cap \fk$ also 
intersects $\z(\fk)$. 

From Lemma~\ref{lem:6.2b}(c) we know that 
$V_\lambda := V \cap \cH_\lambda\not=\{0\}$. 
Since $V^\infty := V \cap \cH^\infty$ 
is dense in $V$ and $P_\lambda$-invariant, the subspace 
$V^\infty_\lambda = P_\lambda(V^\infty)$ is dense in $V_\lambda$. 
The representation $\rho_\lambda$ of $K$ on $V_\lambda$ 
satisfies  $\rho_\lambda(Z(K)_0) \subeq \T \1$, 
which implies that $\z(\fk) + W_{\rho_\lambda} = W_{\rho_\lambda}$. 
As the open invariant cone $W_{\rho_\lambda}$ in $\fk$ also intersects 
$\z(\fk)$ (Lemma~\ref{lem:12.1}), we obtain 
$W_{\rho_\lambda} = \fk$, i.e., $\rho_\lambda$ is bounded. 

Now we apply \cite[Thm.~2.17]{Ne10d} to the closed subspace 
$V_\lambda \subeq \cH$ and the representation of 
$\fp^+ \rtimes \fk_\C$ 
on $V_\lambda$ defined by $\beta(\fp^+) = \{0\}$. 
We derive in particular that $(\pi,\cH)$ is holomorphically 
induced from $(\rho_\lambda, V_\lambda)$ and hence that 
the irreducibility of $\pi$ implies the irreducibility 
of $\rho_\lambda$ (\cite[Cor.~2.14]{Ne10d}). 

We obtain in particular $V_\lambda \subeq \cH^\omega$, so that 
$U(\g)_\C V_\lambda^\infty$ spans a dense subspace of $\cH$. 
Next we use the PBW Theorem to see that 
\[ U(\g_\C)V_\lambda^\infty 
= U(\fp^+)U(\fk_\C) U(\fp^-) V_\lambda^\infty 
= U(\fp^+) V_\lambda^\infty,\] 
and the eigenvalues of $-i\dd\pi(d)$ on this space 
are contained in $\lambda +\N_0$. In particular, 
$\lambda$ is the minimal eigenvalue of $-i\dd\pi(d)$ on $\cH$. 
Since the same argument applies to all other eigenvalues of 
$-i\dd\pi(d)$ in $V$, they are equal, so that $V = V_\lambda$. 
\end{prf}

\begin{rem}
  \mlabel{rem:5.5} 
In view of Theorem~\ref{thm:6.2b}, every irreducible semibounded 
representation $(\pi, \cH)$ of $G$ is determined by the bounded 
irreducible $K$-representation $(\rho, V)$. 
According to \cite[Cor.~2.16]{Ne10d}, 
two such $G$-representations $(\pi_j, \cH_j)$, $j =1,2$, are equivalent if and 
only if the corresponding $K$-representations $(\rho_j, V_j)$ 
are equivalent. Therefore the classification of the separable 
semibounded irreducible $G$-representations means to 
determine the inducible bounded irreducible 
representations $(\rho, V)$ of~$K$. 
\end{rem}

The following theorem follows from \cite[Thm.~3.7, Thm.~3.14]{Ne10d}. 
In particular, it provides a criterion for semiboundedness of a 
$G$-representation in terms of boundedness of a $K$-representation 
plus a positive energy condition. 

\begin{thm} \mlabel{thm:3.15} 
If $(\pi, \cH)$ is a smooth positive energy representation 
for which the $K$-representation $\rho$ 
on $V := \oline{(\cH^\infty)^{\fp^-}}$ 
is bounded, then it is holomorphically induced from $(\rho,V)$ 
and it is semibounded with $d \in W_\pi$. 
\end{thm}

\begin{prop} If the irreducible hermitian Lie group 
$(G,\theta, d)$ has a bounded unitary representation  
with $\dd\pi(\fp) \not= \{0\}$, then the $JH^*$-triple $\fp$ is negative. 
\end{prop}

\begin{prf} Let $(\pi, \cH)$ be a bounded unitary representation  
with $\dd\pi(\fp) \not= \{0\}$. Then 
$\ker \dd\pi \trile \g$ is a closed ideal, so that the 
simplicity of the $JH^*$-triple $\fp$ implies that 
$\fp \cap \ker \dd\pi = \{0\}$, hence that 
$\ker \dd\pi \subeq \ker \ad_\fp = \z(\g)$. We may therefore 
assume that $\dd\pi$ is faithful, so that 
$\|x\| := \|\dd\pi(x)\|$ defines an invariant norm on~$\g$. 

Since $(G,\theta, d)$ is irreducible, 
either $\fp$ is flat, positive or negative. 
If $\fp$ is flat, then $\fp + [\fp,\fp]$ is a $2$-step nilpotent 
Lie subalgebra and $\fp + [\fp,\fp] + \R d$ is solvable 
and locally finite, i.e., every finite subset generates a 
finite dimensional subalgebra. For a finite dimensional 
Lie algebra, the existence of an invariant norm implies that 
it is compact. If it is also solvable, it must be abelian, 
but $\fp + [\fp,\fp] + \R d$ is non-abelian. Therefore 
$\fp$ is not flat. 

If $\fp$ is positive, then \cite[Lemma~3.3]{Ka83} implies the 
existence of an element $z \in \fp$ which is a tripotent for the 
Jordan triple structure, i.e., $\{z,z,z\} = z$. 
Then $\g_z := \R z + \R [d,z] + \R[z,[d,z]]$ 
is a $\theta$-invariant subalgebra of $\g$ isomorphic to 
$\fsl_2(\R)$. This contradicts the existence of an invariant 
norm. Therefore $\fp$ is negative. 
\end{prf}

\begin{rem} If $(\pi, \cH)$ is a semibounded representation 
with $d \in W_\pi$, then the boundedness of $\pi$ is equivalent to the 
boundedness of the operator $\dd\pi(d)$. In fact, 
the boundedness of $\dd\pi(d)$ implies that 
$d \in W_\pi \cap -\oline{W_\pi}$, hence that 
$W_\pi \cap -W_\pi \not=\eset$. But this implies that 
$0 \in W_\pi$, i.e., $W_\pi = \g$. 
\end{rem}

\subsection{Necessary conditions for inducibility} 

In this subsection we discuss some necessary conditions for 
inducibility of a bounded representation $(\rho,V)$ of $K$. 
To evaluate these conditions in concrete cases, we also express them 
in terms of root decompositions.

\begin{defn} (a) Let $C_\fp \subeq \fk$ be the closed convex cone generated by 
the elements of the form $i[z^*,z]$, $z \in \fp^+$. 

(b) For a convex cone $C \subeq \fk$, a smooth 
unitary representation $(\rho, V)$ of $K$ is said to be 
{\it $C$-dissipative} if 
$i\dd\rho(y) \leq 0$ for $y \in C$. 
\end{defn}

\begin{rem} \mlabel{rem:6.4} 
(a) Since its set of generators is $\Ad(K)$-invariant, the cone 
$C_\fp \subeq \fk$ is $\Ad(K)$-invariant. In general it does not have 
any interior points. For the hermitian Lie algebras 
$\g = \fu_{\rm res}(\cH_+, \cH_-)$, $\fsp_{\rm res}(\cH)$ and $\fo^*_{\rm res}(\cH)$, 
the cone $C_\fp$ lies in the subspace of compact operators 
in $\fk$, so that it never has interior points with respect to the 
operator norm. 

(b) Passing from $\g = \fk + \fp$ to  the $c$-dual 
symmetric Lie algebra $\g^c = \fk + i \fp \subeq \g_\C$ does not 
change the subspaces $\fp^\pm \subeq \g_\C= \g^c_\C$, 
but for $z \in \fp_\C$ the corresponding element $z^*$ changes sign. 
This leads to $C_{i\fp} = -C_\fp$.   

(c) The elements in $\fp^+$ can also be written 
as $x - i I x$, $x \in \fp$, where $I = \ad_\fp d$ is the canonical 
complex structure on $\fp$. We then have 
\[ i [(x - i I x)^*, x - i I x] 
=  i [-x - i I x,x - i I x] 
=   2[Ix,x].\] 
\end{rem}

\begin{lem} \mlabel{lem:neccond} If the $K$-representation 
$(\rho, V)$ is inducible, then 
\begin{equation}
  \label{eq:inducness}
\dd\rho([z^*,z]) \geq 0 \quad \mbox{ for } \quad z \in \fp^+,
\end{equation}
i.e., $(\rho, V)$ is $C_\fp$-dissipative. 
\end{lem} 

\begin{prf} If the $G$-representation $(\pi, \cH_V)$ is holomorphically 
induced from $(\rho, V)$, then we have for $v \in (\cH^\infty)^{\fp^-}$ 
and $z \in \fp^+$ the relation 
\[ \la \dd\rho([z^*,z])v,v\ra 
=  \la [\dd\pi(z^*), \dd\pi(z)]v,v\ra 
=  \la \dd\pi(z^*)\dd\pi(z)v,v\ra 
= \|\dd\pi(z)v\|^2 \geq 0.\] 
Since $(\cH^\infty)^{\fp^-}$ is dense in $V$ and 
$\dd\rho([z^*,z])$ is bounded, this proves \eqref{eq:inducness}. 
\end{prf}

\begin{ex} \mlabel{ex:5.12} 
Consider the full hermitian Lie algebra 
$\g = \fu_{\rm res}(\cH)$, where the involution is obtained from a 
decomposition $\cH = \cH_+\oplus \cH_-$ and the element 
\[ d := \frac{i}{2}\pmat{\1 & 0 \\ 0 & -\1}.\] 
Then $\fk = \fu(\cH_+) \oplus \fu(\cH_-)$ and  
$\fp^+ = \pmat{0 & B_2(\cH_-, \cH_+) \\ 0 & 0}.$ 
For $z \in B_2(\cH_-, \cH_+)$ we have 
\[ \Big[ \pmat{0 & z \\ 0 & 0}, \pmat{0 & 0 \\ z^* & 0}\Big] 
=   \pmat{zz^* & 0 \\ 0 & -z^*z}.\] 
Therefore 
\[ C_\fp \subeq \{ (a,d) \in \fk \: 
i\cdot a \geq 0, -i\cdot d \geq 0\}.\] 

(b) For the $c$-dual Lie algebra 
$\g^c= \fu_{\rm res}(\cH_+,\cH_-)$, we find that 
\[ C_\fp \subeq \{ (a,d) \in \fk \: 
-ia \geq 0, id \geq 0\}.\] 
\end{ex}

The most effective way to draw further consequences from 
the necessary condition Lemma~\ref{lem:neccond} 
is to use root decompositions (see Appendix~\ref{app:d.1a} 
for detailed definitions).

\begin{lem} \mlabel{lem:5.16} Suppose that $\g = \fk + \fp$ is hermitian, 
 $\ft \subeq \fk$ is maximal abelian in $\g$, 
and $\g_\C^\alpha = \C x_\alpha$ with 
$\g_\C(x_\alpha) \cong \fsl_2(\C)$. 
Then the following assertions hold: 
 \begin{description}
 \item[\rm(i)] If $x_\alpha \in  \fk_\C$, then $\alpha$ is compact, i.e., 
$\alpha \in \Delta_c$. 
 \item[\rm(ii)]  If $\fp$ is a negative/positive  $JH^*$-triple and 
$x_\alpha \in \fp_\C$, then $\alpha$ is compact/non-compact. 
 \item[\rm(iii)] If $(\beta,V)$ is a $*$-representation 
of $\g_\C$ on a pre-Hilbert space and 
$v_\lambda \in V$ an eigenvector of $\g_\C^{-\alpha} + \ft_\C$, 
then 
\[ \lambda(\check \alpha) 
\begin{cases}
\leq 0  & \text{ for } \alpha \in \Delta_c \\ 
\geq 0  & \text{ for } \alpha \in \Delta_{nc}.
\end{cases}\]
 \end{description}
\end{lem} 

\begin{prf} (i) If $x_\alpha \in \g_\C^\alpha \subeq \fk_\C$, 
then the real Lie algebra $\g(\alpha) := \g_\C(x_\alpha) \cap \g 
\subeq \fk$ 
has a non-trivial bounded representation $\ad_\fp$, hence is compact. 
Now Lemma~\ref{lem:e.1} implies that $\alpha \in \Delta_c$. 

(ii) The $3$-dimensional Lie algebra 
$\g(\alpha) = \g_\C(x_\alpha) \cap \g$ is hermitian with 
$d_\alpha := -\frac{i}{2} \check \alpha$ and 
$\fp^+ = \C x_\alpha$. Now 
\[ \{x_\alpha, x_\alpha, x_\alpha\} 
= [[x_\alpha, \oline{x_\alpha}], x_\alpha]
= -[[x_\alpha, x_\alpha^*], x_\alpha]
= \alpha([x_\alpha^*, x_\alpha]) x_\alpha.\] 
Therefore the $JH^*$-triple $\C x_\alpha$ is positive if 
$\g(\alpha) \cong \su_{1,1}(\C)$ and 
negative if $\g(\alpha) \cong \su_2(\C)$ (Lemma~\ref{lem:e.1}). 

(iii) As in the proof of Lemma~\ref{lem:neccond}, we obtain for a 
unit vector $v_\lambda$ of weight $\lambda$ and $x_\alpha \in \g_\C^\alpha$: 
\[ \lambda([x_\alpha^*,x_\alpha])
= \la \beta([x_\alpha^*,x_\alpha])v,v\ra 
= \la \beta(x_\alpha)^*\beta(x_\alpha)v,v\ra 
= \|\beta(x_\alpha)v\|^2 \geq 0.\]  
Now (iii) follows from \eqref{eq:corootsign} in 
Definition~\ref{def:coroot}. 
\end{prf}

\begin{ex}\mlabel{ex:5.12b} (Example~\ref{ex:5.12} continued) 
For the Hilbert space $\cH = \cH_+ \oplus \cH_-$,  
we choose an orthonormal basis $(e_j)_{j \in J}$ such that 
$J= J_+ \cup J_-$, where $(e_j)_{j \in J_\pm}$ is an 
orthonormal basis of $\cH_\pm$. 

(a) For the Lie algebra 
$\g = \fu_{\rm res}(\cH)$, the subalgebra 
\[ \ft := \{ X \in \g \: (\forall j\in J)\, 
X e_j \in \C e_j \} \cong \ell^\infty(J,i\R) \] 
is elliptic and maximal abelian with root system 
\[ \Delta = \Delta_c = \{ \eps_j - \eps_k \: j\not=k \in J\}\] 
(cf.\ Example~\ref{ex:d.1a}). 
From $\Delta_p^+ = \{ \eps_j - \eps_k \: j\in J_+, k \in J_-\}$ 
(Example~\ref{ex:d1}),  we therefore get 
$E_{jj} - E_{kk} \in i C_{\fp}$ for $j \in J_+, k \in J_-.$ 
For any weight $\mu = (\mu_j)_{j \in J}$ of a representation 
$(\rho,V)$ of $K$, the relation $\dd\rho([z^*,z]) \geq 0$ for 
each $z \in \fp^+$ now implies 
\begin{equation}
  \label{eq:weight1}
\mu_j \leq\mu_k \quad \mbox{ for } \quad j \in J_+, k \in J_-,
\quad \mbox{ i.e., } \quad 
\sup(\mu_+) \leq \inf(\mu_-) 
\ \ \  \mbox{ for } \ \ \ 
\mu_\pm := \mu\res_{J_\pm}. 
\end{equation}

(b)  For the $c$-dual Lie algebra 
$\g^c= \fu_{\rm res}(\cH_+,\cH_-)$,  
the sets $\Delta_k$ and $\Delta_p^\pm$ are the same, but 
we have $E_{jk}^* = - E_{kj}$ for $j \in J_+, k \in J_-$, 
so that $\Delta_c = \Delta_k$ and $\Delta_{nc} =\Delta_p$. 
For any weight $\mu = (\mu_j)_{j \in J}$ of a representation 
$(\rho,V)$ of $K$, the relation $\dd\rho([z^*,z]) \geq 0$ therefore implies 
\begin{equation}
  \label{eq:weight2}
\mu_j \geq\mu_k \quad \mbox{ for } \quad j \in J_+, k \in J_-,
\quad \mbox{ i.e.,} \quad 
\inf(\mu_+) \geq \sup(\mu_-).
\end{equation}
\end{ex}

\subsection{The Reduction Theorem} 

Suppose that $(G,\theta,d)$ is a simply connected 
full hermitian Lie group for which $K$ is a direct product 
of a finite dimensional Lie group with compact Lie algebra 
and groups isomorphic to full unitary groups 
$\U(\cH)$, where $\cH$ is an infinite dimensional real, complex 
or quaternionic Hilbert space. 
We write $K_\infty \trile K$ for the normal integral subgroup 
obtained by replacing each $\U(\cH)$-factor 
by $\U_\infty(\cH)$. In Theorem~\ref{thm:kredux} below we show that 
bounded irreducible representations $(\rho, V)$ of 
$K$ are tensor products of two bounded irreducible 
representations $(\rho_0, V_0)$ and $(\rho_1, V_1)$, 
where $\rho_0\res_{K_\infty}$ is irreducible and 
$K_\infty \subeq \ker \rho_1$. In this section we see 
how this can be used to derive 
a factorization of irreducible semibounded representations, resp., 
to reduce the inducibility problem to the case~$\rho = \rho_0$. 

We start with a  general lemma on unitary representations 
of discrete groups. 

\begin{lem} \mlabel{lem:factorize} {\rm(Factorization Lemma)}
Let $G$ be a group, $N \trile G$ be a normal 
subgroup and  $(\pi, \cH)$ be an irreducible unitary representation 
of $G$ such that the restriction $\pi\res_N$ contains an 
irreducible representation $(\pi_0, \cH_0)$ 
whose equivalence class is $G$-invariant. 
Let 
\[ \hat G := \{ (g,u) \in G \times \U(\cH_0) \: 
(\forall n \in N)\ 
\pi_0(gng^{-1}) = u \pi_0(n) u^* \} \] 
denote the corresponding central $\T$-extension 
of $G$ with kernel $Z := \{ \1\} \times \T\1$ and 
note that $\sigma \: N \to \hat G, n \mapsto (n,\pi_0(n))$ defines a 
group homomorphism. Then 
\begin{description}
\item[\rm(a)] $\hat\pi_0(g,u) := u$ 
defines an irreducible representation 
$(\hat\pi_0,\cH_0)$ of $\hat G$ for which 
$\hat\pi_0 \circ \sigma = \pi_0$ is irreducible, and there exists  
\item[\rm(b)] an irreducible representation 
$(\hat\pi_1,\cH_1)$ of $\hat G$ with 
$\sigma(N)  \subeq \ker \hat\pi_1$ and 
$\hat\pi_1(\1,t) = t^{-1}\1$ for $t \in \T$,  
\end{description}
 such that 
$q^*\pi \cong \hat\pi_0 \otimes \hat\pi_1$ holds for the 
projection map $q \: \hat G \to G$. 
\end{lem}

\begin{prf} Our assumption implies that $G$ 
preserves the $\pi_0$-isotypic subspace for $\pi\res_N$, 
so that the irreducibility of $\pi$ implies that 
$\pi\res_N$ is of the form 
$\pi_0 \otimes \pi_1$, where 
$(\pi_0, \cH_0)$ is irreducible and $(\pi_1, \cH_1)$ is trivial.
Since the class $[\pi_0] \in \hat N$ is $G$-invariant, 
the map $\hat G \to G, (g,u) \mapsto g$ is surjective 
with kernel $\{\1\} \times \T\1$, 
hence a central extension. 
Clearly, $(\hat\pi_0, \cH_0)$ is a unitary representation 
of $\hat G$ whose restriction to $\sigma(N)$ is irreducible. 
Since the operators 
$(\hat\pi_0(g,u) \otimes \1)^{-1}\pi(g)$  
commute with $\pi(N)$, they are of the form 
$\1 \otimes \hat\pi_1(g,u)$ for some unitary 
operator $\hat\pi_1(g,u) \in \U(\cH_1)$. 
By definition, $\hat\pi_1(n,\pi_0(n)) = \1$ for $n \in N$. 
That $\hat\pi_1$ is a representation follows from 
\begin{align*}
&\ \ \ \1 \otimes \hat\pi_1((g_1,u_1)(g_2,u_2)) 
= \big(\hat\pi_0(g_1g_2,u_1u_2) \otimes \1\big)^{-1}\pi(g_1g_2) \\
&= (u_2^{-1} u_1^{-1} \otimes \1) \pi(g_1) \pi(g_2) 
= (u_2^{-1} \otimes \1)(\1 \otimes \hat\pi_1(g_1,u_1)) \pi(g_2) \\
&= (\1 \otimes \hat\pi_1(g_1,u_1))(u_2^{-1} \otimes \1) \pi(g_2) 
= (\1 \otimes \hat\pi_1(g_1,u_1))(\1  \otimes \hat\pi_1(g_2,u_2)).
\qedhere\end{align*}
\end{prf}

\begin{lem} \mlabel{lem:product} 
Let $K_j$, $j= 1,\ldots, n$, be Lie groups for which 
all bounded unitary representations are direct sums of irreducible 
ones. Then the same holds for the product group 
$K := \prod_{j = 1}^n K_j$ and each bounded irreducible representation 
of $K$ is a tensor product of bounded irreducible representations 
of the~$K_j$. 
\end{lem}

\begin{prf} We argue by induction, so that it suffices to consider the 
case $n=2$. Let $(\pi, \cH)$ be a bounded unitary representation 
of $K =K_1 \times K_2$. Since $\pi\res_{K_1}$ is a direct sum of 
irreducible representations, we obtain a decomposition 
of $\cH$ into isotypical subspaces $\cH_j$ which are also invariant 
under $K_2$. Hence we may assume that $\pi\res_{K_1}$ is isotypical, so 
that $\cH \cong \cH_1 \otimes \cH_2$ and 
$\pi(k_1, k_2) = \pi_1(k_1) \otimes \pi_2(k_2)$, where 
$(\pi_1, \cH_1)$ is an irreducible representation of $K_1$. 
As $(\pi_2, \cH_2)$ also is a direct sum of irreducible representations, 
the assertion follows from the fact that tensor products of  irreducible 
representations of $K_1$, resp., $K_2$, define irreducible 
representations of~$K_1 \times K_2$. 
\end{prf}

\begin{thm} \mlabel{thm:kredux} 
Suppose that $K = \prod_{j = 0}^n K_j$ is a product of connected Lie groups, 
where $K_0$ has a compact Lie algebra and 
the groups $K_j$, $j > 0$, are isomorphic to a full unitary 
group $\U(\cH_j)$ of an infinite dimensional 
Hilbert space $\cH_j$ over $\R, \C$ or $\H$. 
Accordingly, we define 
$K_{j,\infty}$ as $\U_\infty(\cH_j)$. 
Then 
\[ K_\infty := K_0 \times \prod_{j = 1}^n K_{j,\infty} \] 
is a normal Lie subgroup and each irreducible 
bounded representation $(\rho,V)$ of $K$ 
is a tensor product of two bounded irreducible 
representations $(\rho_0, V_0)$ and $(\rho_1, V_1)$, 
where $\rho_0\res_{K_\infty}$ is irreducible and 
$K_\infty \subeq \ker \rho_1$. 
\end{thm}

\begin{prf} Each bounded representation of $K_j$ 
extends to a holomorphic representation of its universal 
complexification $K_{j,\C}$. For the groups of the form 
$\U_\infty(\cH_j)$ where $\cH_j$ is complex, 
it therefore follows from \cite[Thm.~III.14]{Ne98} 
that all bounded representations are direct sums of irreducible 
ones, and if $\cH_j$ is a Hilbert spaces over $\K = \R,\H$, 
this is a consequence of 
Theorems~\ref{thm:oprep} and \ref{thm:sprep}. 
Since $K_0$ has a compact Lie algebra, it is a product 
$K_0 = Z_0 \times K_0'$, where $K_0'$ is compact semisimple 
and $Z_0$ is abelian. Since every continuous unitary representation 
of $K_0'$ is a direct sum of irreducible ones 
and $\rho(Z_0) \subeq \T\1$ by Schur's Lemma, 
Lemma~\ref{lem:product} implies that the restriction $\rho\res_{K_\infty}$ 
is a direct sum of irreducible representations which 
in turn are tensor products of irreducible ones. 

Since each bounded irreducible representation 
of $\U_\infty(\cH)$ extends to $\U(\cH)$ 
(cf.\ Definition~\ref{def:tensrep} and 
Theorems~\ref{thm:oprep} and \ref{thm:sprep}), 
the same holds for the 
bounded irreducible representations of $K_\infty$. 
Hence the assertion follows  from Lemma~\ref{lem:factorize}, which 
applies with the trivial central extension $\hat K \cong \T \times K$ 
of $K$. 
\end{prf}

Let $(G,\theta,d)$ be a simply connected full irreducible 
hermitian Lie group for which 
$K = (G^\theta)_0$ satisfies the assumption of Theorem~\ref{thm:kredux}. 
Then  $[\fp,\fp] \subeq \fk_\infty = \L(K_\infty)$ 
follows from the concrete description of $\fp$ in all types 
(cf.\ Remark~\ref{rem:2.7}). 
Therefore each bounded irreducible representation 
$(\rho, V)$ of $K$ is a tensor product 
of two bounded irreducible 
representations $(\rho_0, V_0)$ and $(\rho_1, V_1)$, 
where $\rho_0\res_{K_\infty}$ is irreducible and 
$K_\infty \subeq \ker \rho_1$. 

From $[\fp,\fp] \subeq \fk_\infty$ we derive the existence 
of a homomorphism $p \: G \to K/K_\infty$, 
so that, if $\oline\rho_1 \: K/K_\infty \to \U(\cH_1)$ denotes the 
induced representation of the quotient group, then 
 $\pi_1 := \oline\rho_1  \circ p$ is a bounded 
irreducible representation of $G$. 
If $(\rho_0, V_0)$ is extendable to a semibounded 
representation $(\pi_0, \cH_0)$ of~$G$, 
then 
\[ \pi := \pi_0 \otimes \pi_1, \qquad \cH := \cH_0 \otimes V_1\] 
defines a semibounded unitary representation of 
$G$ on $\cH$ which is holomorphically induced from 
the irreducible representation 
$(\rho, V)$, hence irreducible. 

\begin{thm} \mlabel{thm:redux} 
If $(G,\theta,d)$ is a simply connected full irreducible 
hermitian Lie group for which 
$K = (G^\theta)_0$ satisfies the assumptions of Theorem~\ref{thm:kredux}, 
then the representation $(\rho,V)$ is holomorphically inducible if and only 
if $(\rho_0, V_0)$ has this property. 
\end{thm}

\begin{prf} We know from Theorem~\ref{thm:c.2} below that 
$(\rho,V)$ is holomorphically inducible if and only if 
the corresponding function 
$f_\rho$ on a $\1$-neighborhood of $G$, defined by 
\[ f_\rho(\exp x) = F_\rho(x) \quad \mbox{ and } \quad 
F_\rho(x_+ * x_0 * x_-) = e^{\dd\rho(x_0)}, \quad 
x_\pm \in \fp^\pm, x_0 \in \fk_\C, \]  
is positive definite on a $\1$-neighborhood in~$G$.
From the factorization of $\rho$, we obtain a factorization 
\[ f_\rho(g) 
= f_{\rho_0}(g) \otimes f_{\rho_1}(g)
= f_{\rho_0}(g) \otimes \pi_1(g), \] 
which leads to 
\[ f_\rho(gh^{-1}) 
= f_{\rho_0}(gh^{-1}) \otimes \pi_1(gh^{-1}) 
= (\1 \otimes \pi_1(g))(f_{\rho_0}(gh^{-1}) \otimes \1)
(\1 \otimes \pi_1(h))^*.\] 
Therefore the function $f_\rho$ is positive definite if and only if 
$f_{\rho_0}$ is positive definite (cf.\ Remark~\ref{rem:a.2}), 
and this 
shows that $(\rho, V)$ is inducible if and only if 
$(\rho_0, V_0)$ has this property. 
\end{prf}

\begin{rem}
  \mlabel{rem:5.14}
In view of Theorem~\ref{thm:redux}, the classification 
problem for irreducible semibounded 
representations of $G$ reduces to the classification 
of the inducible irreducible bounded representations 
for which $\rho\res_{K_\infty}$ is irreducible and the disjoint 
problem of the classification of irreducible bounded representations 
of the group $K/K_\infty$. The latter contains in particular the 
classification problem 
for irreducible bounded representations of groups 
like $\U(\cH)/\U_\infty(\cH)$, which is the identity component 
of the unitary group of the $C^*$-algebra $B(\cH)/K(\cH)$. 
For a construction procedure for 
irreducible unitary representations 
of unitary groups of $C^*$-algebras, we refer to \cite{BN11}. 
\end{rem}

\subsection{Triple decompositions in complex groups} 

This brief subsection sets the stage for a uniform treatment of the 
inducibility problem for bounded $K$-representations $(\rho,V)$ 
for groups related to positive and negative $JH^*$-triples. 

\begin{lem} \mlabel{lem:tridec} 
Let $(G,\theta, d)$ be a hermitian Lie group,  
$G_c$ be a complex Lie group with Lie algebra $\g_\C$, 
$P^\pm := \exp \fp^\pm$, and 
$K_c := \la \exp \fk_\C \ra$ be the Lie subgroup 
of $G_c$ corresponding to $\fk_\C$. Then the  multiplication map 
\[ \mu \:  P^+ \times K_c \times P^- \to G_c, \quad (p_+,k,p_-) 
\mapsto p_+ k p_- \] 
is biholomorphic onto an open subset of $G_c$ and 
the subgroups $P^\pm$ are simply connected. 
\end{lem} 

\begin{prf} The existence of the Lie subgroups 
$P^\pm$ and $K_c$ of $G_c$ follows from the fact the 
Lie algebras $\fp^\pm$ and $\fk_\C$ are closed 
(\cite{Mais62}). 
That the subgroups $P^\pm$ are simply connected follows from 
\[ \Ad(\exp z)d = e^{\ad z} d = d + [z,d] 
= d \mp iz\quad \mbox{ for } \quad 
z \in \fp^\pm, \] 
which implies that $\exp\res_{\fp^\pm}$ is injective. 
It also shows that, for $p_\pm \in P^\pm$, the relation 
$\Ad(p_\pm) d = d$ implies $p_\pm = \1$. 

The map $\mu$ is an orbit map for the holomorphic 
action of the direct product group $P^+ \times (K_c \rtimes P^-)^{\rm op}$ 
on $G_c$ by $(g,h).x := gxh$. Since its differential 
in $(\1,\1,\1)$ is the summation map, hence invertible, 
it follows from the complex Inverse Function Theorem 
that $\mu$ is a covering map onto an open subset 
of $G_c$. Therefore it remains to show that the fiber of 
$\1$ is trivial. If $p_+ k p_- = \1$, then 
$p_+ = p_-^{-1} k^{-1}$ implies that 
\[ \Ad(p_+)d  -d 
= \Ad(p_-^{-1} k^{-1})d -d 
= \Ad(p_-^{-1})d -d \in \fp^+ \cap \fp^- = \{0\},\] 
so that $p_+ = p_- = \1$. 
It follows that $\mu$ is injective, hence biholomorphic onto its image. 
\end{prf}

Now we consider a simply connected hermitian 
Lie group $(G,\theta, d)$ for which 
there exists a simply connected Lie group $G_\C$ with Lie algebra 
$\g_\C$. This assumption is satisfied if $G$ is the simply connected 
Lie group with Lie algebra $\g = \hat\g(\fp)$ 
(Theorem~\ref{thm:exist}). 
Then there exists a homomorphism $\eta_G \: G \to G_\C$ 
integrating the inclusion $\g \into \g_\C$. 
Since $\L(\eta_G)$ is injective, the kernel of $\eta_G$ is 
discrete. This in turn implies that $\eta_G\res_K \: K \to G_\C$ 
has a discrete kernel, so that the Complexification Theorem 
in \cite[Thm.~IV.7]{GN03} implies the existence of a 
universal complexification 
$\eta_K \: K \to K_\C$ of $K$. 
Let $K_c := \la \exp \fk_\C \ra \subeq G_\C$ denote the 
integral subgroup corresponding to $\fk_\C \subeq \g_\C$ 
and $q_c \: K_\C \to K_c$ be the natural map whose existence follows 
from the universal property of $K_\C$. 
We then obtain a covering map 
\[ \tilde\mu \: P^+ \times K_\C \times P^- \to 
\Omega := P^+ K_c P^- \subeq G_\C, \quad 
(p^+, k, p^-) \mapsto p^+ q_c(k) p^-. \] 
Let $U \subeq \fp$ be a $K$-invariant open convex symmetric 
$0$-neighborhood 
for which the map $U \times K \to \exp(U)K, (x,k) \mapsto \exp x \cdot k$ 
is a diffeomorphism onto an open subset $U_G = \exp(U)K$ 
of $G$ contained in $\eta_G^{-1}(\Omega)$. The $K$-invariance and 
the symmetry of $U$ imply that 
\[ U_G = \exp(U) K = K \exp U = K \exp(-U) = U_G^{-1}.\]
Let 
\[ \tilde\eta_G \: U_G \to P^+ \times K_\C \times P^- \] 
be the unique continuous lift of 
$\eta_G\res_{U_G}$ with $\tilde\eta_G(\1) = (\1,\1,\1)$ and
$\tilde\mu \circ \tilde\eta_G = \eta_G\res_{U_G}$. 
Writing $x = (x_+, x_0, x_-)$ for the elements of 
$P^+ \times K_\C \times P^-$, we thus obtain an analytic map 
\[ \kappa \: U_G \to K_\C, \quad g \mapsto \tilde\eta_G(g)_0.\]
From the uniqueness of lifts and the fact that $K$ normalizes 
$P^\pm$, we derive that $\kappa$ is $K$-biequivariant, i.e., 
\[ \kappa(k_1 gk_2) = k_1 \kappa(g) k_2 \quad \mbox{ for } \quad 
g \in G, k_1, k_2 \in K.\] 

For a bounded representation $(\rho,V)$ of $K$ and its holomorphic 
extension 
\[ \rho_\C \: K_\C \to \GL(V)\quad \mbox{ with } \quad 
\rho_\C \circ \eta_K = \rho, \] 
 we now obtain the analytic function 
\begin{equation}
  \label{eq:frho-1}
f_\rho := \rho_\C \circ \kappa \: U_G \to \GL(V)
\end{equation}
which coincides with the corresponding function 
in Theorem~\ref{thm:c.2} in an identity neighborhood. 
Now let $W \subeq U$ be an open symmetric $K$-invariant $0$-neighborhood 
for which $W_G \subeq U_G$ satisfies $W_G W_G^{-1} = W_G W_G\subeq U_G$. 
We write $M_W := W_G/K$ for the corresponding open subset of 
$M  = G/K$. From the equivariance property 
\[ f_\rho(k_1 g k_2) = \rho(k_1) f_\rho(g) \rho(k_2) \quad \mbox{ for } 
\quad g \in G, k_1, k_2 \in K,\] 
it follows that $f_\rho$ defines a kernel function 
\[ F_\rho \: M_W \times M_W \to B(V), \quad 
(gK,hK) \mapsto (f_\rho(h)^{-1})^* f_\rho(h^{-1}g) f_\rho(g)^{-1}. \] 
The kernel $F_\rho$ is hermitian in the sense that 
\[ F_\rho(z,w)^* = F_\rho(w,z) \quad \mbox{ for } \quad 
z,w \in M_W.\] 
Moreover, for each $h \in G$, the function 
\[ G \to B(V), \quad g \mapsto 
(f_\rho(h)^{-1})^* f_\rho(h^{-1}g) f_\rho(g)^{-1} \] 
is annihilated by the left invariant differential operators 
$L_w$, $w \in \fp^-$, because $f_\rho$ has this property, 
which in turn is a consequence of the construction of~$\kappa$. 
Therefore $F_\rho$ is holomorphic in the first argument 
(\cite[Def.~1.3]{Ne10d}), 
and since it is hermitian, it is antiholomorphic in the second 
argument, i.e., holomorphic as a function on the complex manifold 
$M_W \times \oline{M_W}$
because we know already that it is smooth. 
Here $\oline{M_W}$ denotes the complex 
manifold $M_W$, endowed with the opposite complex structure.

\begin{prop}\mlabel{prop:kernel-1} 
Let $(G,\theta,d)$ be a hermitian Lie group 
for which $\g_\C$ is integrable and 
$M_W = W_G/K$ be as above. 
Then the following are equivalent for a bounded unitary representation 
$(\rho,V)$ of $K$: 
\begin{description}
\item[\rm(i)] The representation $(\rho,V)$ is holomorphically 
inducible. 
\item[\rm(ii)] The function $f_\rho$ is positive definite on 
$W_G W_G^{-1}$. 
\item[\rm(iii)] The holomorphic kernel $F_\rho$ on $M_W \times \oline{M_W}$ 
is positive definite. 
\end{description}
\end{prop}

\begin{prf} Since $f_\rho$ is analytic, 
the Extension Theorem~\ref{thm:extension} 
implies that it is positive definite if and only 
if it is positive definite on some identity neighborhood. 
Therefore the equivalence of (i) and (ii) follows from  
Theorem~\ref{thm:c.2}. 

Further, Remark~\ref{rem:a.2} implies that the kernel 
$F_\rho$ is positive definite if and only if 
$f_\rho$ is a positive definite function, so that 
(ii) and (iii) are also equivalent. 
\end{prf}

\begin{rem}\mlabel{rem:5last}
Suppose that $\fp$ is a simple $JH^*$-triple and  
$\g = \hat\g(\fp)$, so that $\g_\C$ is integrable 
(Theorem~\ref{thm:exist}). 
As we have seen in Remark~\ref{rem:3.9}, 
the simply connected group $G$ with Lie algebra 
$\hat\g(\fp)$ is a quotient 
of a semidirect product $G_1 \rtimes K \to G$, 
where $G_1$ is the simply connected covering of a group of operators. 
Then $f_\rho$ is positive definite on $W_GW_G^{-1}$ 
if and only if its pullback 
$\tilde f_\rho$ to $G_1 \rtimes K$ is positive definite, 
but this function satisfies 
\begin{align*}
\tilde f_\rho((g_2,k_2)^{-1}(g_1, k_1))
&= \rho(k_2)^{-1} \tilde f_\rho((g_2^{-1},\1)(g_1,\1)) \rho(k_1)\\
&= \rho(k_2)^* \tilde f_\rho((g_2^{-1}g_1,\1)) \rho(k_1),
\end{align*}
so that $\tilde f_\rho$ is positive definite if and only 
if the restriction $\tilde f_\rho\res_{G_1}$ is positive definite 
(Remark~\ref{rem:a.2}). 
For $K_1  := (G_1^\theta)_0$, the function 
$\tilde f_\rho$ only depends on the representation 
$\tilde\rho  \: K_1 \to \U(V)$ obtained by pulling $\rho$ back 
by  the canonical map $K_1 \to K$. 
\end{rem}

\section{Motion groups and their central extensions} 
\mlabel{sec:6}

In this section we obtain a complete classification of 
the irreducible semibounded representations of 
the hermitian motion group $G = \fp \rtimes_\alpha K$ 
and of its canonical 
central extension $\hat G = \Heis(\fp) \rtimes_\alpha K$ 
(Example~\ref{ex:1.3}). 
For $G = \fp \rtimes K$, Theorem~\ref{thm:6.1} asserts 
that all semibounded unitary representation are trivial on $\fp$, 
hence factor through (semibounded) representations of $K$. 
If, in addition, $\fk/\z(\fk)$ contains no invariant cones, 
then Proposition~\ref{prop:4.5} shows that 
all semibounded irreducible representations of $K$ are bounded. 
For the central extension 
$\hat G$, Theorem~\ref{thm:classi-mot} asserts that a bounded 
irreducible representation $(\rho,V)$ of $\hat K = \R \times K$ 
with $\rho(t,\1) = e^{ict}\1$ on the center is inducible if and only 
if $c \geq 0$.

\begin{thm} \mlabel{thm:6.1} If $[\fp,\fp] = \{0\}$ and $G = \fp \rtimes K$, 
then all semibounded unitary representations $(\pi, \cH)$ 
factor through $K$, i.e., $\fp \subeq \ker \dd\pi$. 
\end{thm}

\begin{prf} Let $(\pi, \cH)$ be a semibounded unitary 
representation of $G$ and 
$\hat G := G \times \R$. We consider the representation 
$\hat\pi$ of $\hat G$ defined by $\hat\pi(g,t) := e^{it} \pi(g)$. 
Since $W_\pi \not=\eset$, the open cone 
$C_{\hat\pi}$ of all elements $x$ of 
$\hat\g = \g \oplus \R$ for which $s_{\hat\pi} < 0$ holds 
on a neighborhood of $x$ is non-empty. 
For the abelian ideal $\fp$ of $\hat\g$, we thus obtain 
from \cite[Lemma~7.6]{Ne10c} that 
\[ C_{\hat\pi} + \fp = C_{\hat\pi} + [d,\fp] \subeq 
C_{\hat\pi} + [\fp,\hat\g] = C_{\hat\pi},\]  
and hence that $\fp \subeq H(C_{\hat\pi}) = \ker \dd\pi$. 
\end{prf}

Now we turn to the central extension 
\[ \hat G = \Heis(\fp) \rtimes K 
\cong \fp \oplus (\R \oplus K)\]  
with Lie algebra $\hat\g \cong \R \oplus_{\omega_\fp} \g$, 
where $\omega_\fp(x,y) = 2\Im \la x,y\ra$. We assume that 
$\fk/\z(\fk)$ contains no open invariant cones. 
From Theorem~\ref{thm:6.2b} we know that 
every irreducible semibounded representation $(\pi, \cH)$ of positive energy 
is holomorphically induced from the bounded $\hat K$-representation 
on $V = \oline{(\cH^\infty)^{\fp^-}}$. We now turn to the question 
which bounded $\hat K$-representations $(\rho, V)$ are holomorphically 
inducible. 

In view of Remark~\ref{rem:6.4}(c), the 
cone $C_\fp$ is generated by the elements 
of the form $[Ix,x]$, $x \in \fp$. 
From 
\[ [Ix,x] = 2\Im \la Ix,x\ra = 2 \|x\|^2\] 
it now follows that 
\[ C_\fp = \R^+ \times \{0\} \subeq \R \times \fk,\] 
which leads to the necessary condition 
\begin{equation}
  \label{eq:neccond3}
 -i\dd\rho(1,0)\geq 0.
\end{equation}
Any irreducible representation 
$(\rho,V)$ of $\hat K \cong \R \times K$ is of the form 
$\rho(t,k) = e^{ic t} \rho_1(k)$, where 
$(\rho_1, V)$ is an irreducible $K$-representation. 
Now condition \eqref{eq:neccond3} means that 
$c \geq 0$. 

For $c = 0$, we obtain the unique extension 
to $G$ by $\pi(t,v,k) := \rho(t,k) = \rho_1(k)$ satisfying 
$\fp \subeq \ker \pi$. To deal with the cases $c > 0$, 
we have to introduce the canonical representation 
of the Heisenberg group on the Fock space 
$S(\fp) = \hat\bigoplus_{n \geq 0} S^n(\fp)$ 
(cf.\ \cite{Ot95}). 
On the dense subspace $S(\fp)_0 = \sum_{n = 0}^\infty S^n(\fp)$ of 
$S(\fp)$, 
we have for each $v \in \fp$ the {\it creation operator} 
$$ a^*(v)(v_1 \vee \cdots \vee v_n) := v \vee v_1 \vee \cdots \vee v_n. $$
This operator has an adjoint $a(v)$ on $S(\fp)_0$, given by 
$$ a(v)\Omega = 0, \quad 
a(v)(v_1 \vee \cdots \vee v_n) 
= \sum_{j = 0}^n \la v_j, v \ra v_1 \vee \cdots \vee \hat{v_j} \vee 
\cdots \vee v_n, $$
where $\hat{v_j}$ means omitting the factor $v_j$. 
One easily verifies 
that these operators satisfy the 
canonical commutation relations (CCR): 
\begin{equation}
  \label{eq:ccr}
[a(v), a(w)] = 0, \quad [a(v), a^*(w)] = \la w, v\ra \1. 
\end{equation}
For each $v \in \fp$, the operator $a(v) + a^*(v)$ on $S(\fp)_0$ 
is essentially self-adjoint and 
$$ W(v) := e^{\frac{i}{\sqrt{2}}\oline{a(v) + a^*(v)}} \in \U(S(\fp)) $$
is a unitary operator. These operators satisfy the {\it Weyl 
relations} 
$$ W(v) W(w) = W(v+w) e^{\frac{i}{2}\Im \la v, w\ra} \quad \mbox{ for } 
\quad v,w \in \fp $$
(cf.\ \cite{Ne10c}). 
For the Heisenberg group $\Heis(\fp) := \R \times \fp$ 
with the multiplication 
\[  (t,v) (t',v') := (t + t' + \Im \la v,v'\ra, v + v') \] 
we thus obtain by $W(t, v) := e^{it/2} W(v)$
a unitary representation on $S(\fp)$, called the {\it 
Fock representation} which is actually 
smooth (cf.\ \cite[Sect.~9.1]{Ne10c}). Using the natural representations 
\[ S^n(U)(v_1 \vee \cdots \vee v_n) 
:= Uv_1 \vee \cdots\vee Uv_n \] 
of the unitary group $\U(\fp)$ on the symmetric powers $S^n(\fp)$, 
and combining them to a unitary representation 
$(S, S(\fp))$, we obtain a 
smooth representation $(\pi_s, S(\fp))$, defined by 
\[ \pi_s(t,v,k) := W(t,v)S(k) = e^{it/2}W(v)S(k) \] 
of $\hat G = \Heis(\fp) \rtimes K$ on $S(\fp)$ with 
$\Omega \in S(\fp)^\infty$. 
Then $\dd\pi_s(\fp^+) = a^*(\fp)$ and $\dd\pi_s(\fp^-) = a(\fp)$ 
lead to $(S(\fp)^\infty)^{\fp^-} = \C \Omega$, 
and it follows from \cite[Thm.~2.17]{Ne10d} 
that $(\pi_s, S(\fp))$ is holomorphically induced 
from the character 
$\rho_s \: \hat K = \R \times K \to \T$, given by 
$\rho_s(t,k) = e^{it/2}$. 
Composing with the automorphisms of 
$\Heis(\cH)$ given by $\mu_h(t,v) := (h^2 t, h v)$, $h > 0$,  
we see that all characters $\rho_c(t,k) := e^{ict}\1$ of $\hat K$ 
are holomorphically inducible. 
We write $(\pi_c, S(\fp))$ for the corresponding 
unitary representation of $\hat G$.

With $\pi_1(t,v,k) := \rho_1(k)$ we now put 
\[ \pi := \pi_c \otimes \pi_1 \quad \mbox{ on } \quad 
\cH := S(\fp) \otimes V_1\]  
an observe that this unitary representation is holomorphically 
induced from \break $\rho = \rho_c \otimes \rho_1$. Its semiboundedness 
follows from Theorem~\ref{thm:3.15}. 
This completes the proof of the following theorem: 

\begin{thm} \mlabel{thm:classi-mot} 
A bounded irreducible unitary representation 
$(\rho, V)$ of $\hat K = \R \times K$ with 
$\rho(t,\1) = e^{ict}\1$  
is holomorphically inducible to a semibounded positive energy 
representation of  $\hat G= \Heis(\fp) \rtimes K$ if and only if $c \geq 0$. 
\end{thm}

\begin{rem} If $\fk/\z(\fk)$ contains no open invariant cones, 
then Theorem~\ref{thm:posen}(iii) implies that 
every semibounded unitary representation 
$\pi$ of $\hat G$ satisfies $d \in W_\pi \cup - W_\pi$. 
If $d \in W_\pi$, then it is of positive energy and if, in addition, 
$\pi$ is irreducible, 
Theorem~\ref{thm:6.2b} implies that it is holomorphically induced 
from a bounded representation of $\hat K$, hence covered by 
Theorem~\ref{thm:classi-mot}. 
\end{rem}

\begin{ex} (a) Theorem~\ref{thm:classi-mot} applies in particular 
to the oscillator group 
$\hat G = \Osc(\cH) = \Heis(\cH) \rtimes \T$ of a complex 
Hilbert space $\cH$ from Example~\ref{ex:1.3}(b). 
In this case $G = \cH \rtimes \T$ is the minimal hermitian 
Lie group with $\fp = \cH$. 

In the physics literature the unitary representations 
of $\Osc(\cH)$ are known as the representations of the 
CCR with a number operator. It is known 
that the positive energy representations 
of $\Osc(\cH)$ are direct sums of Fock representations 
(cf.\ \cite{Ch68}). Therefore our Theorem~\ref{thm:classi-mot} 
generalizes Chaiken's uniqueness result on the 
irreducible positive energy representations of $\Osc(\cH)$. 

(b) The maximal hermitian group 
$(G,\theta, d)$ with $\fp = \cH$ is 
$G = \cH \rtimes \U(\cH)$. For this group,  
Theorem~\ref{thm:classi-mot} 
provides a complete classification of the 
semibounded positive energy representations in terms 
of bounded irreducible representations $(\rho, V)$ of 
$\U(\cH)$. Those representations for which 
$\rho\res_{\U_\infty(\cH)}$ is irreducible 
are described in 
Appendix~\ref{app:e} (cf.\ Definition~\ref{def:tensrep}). 
According to Theorem~\ref{thm:8.1b}, all separable continuous 
irreducible representations of $\U(\cH)$ have this property. 
\end{ex}

\section{Hermitian groups with positive $JH^*$-triples} 
\mlabel{sec:7}

In this section we consider the case where 
$M = G/K$ is an infinite dimensional 
irreducible symmetric Hilbert domain, i.e., 
$\fp$ is a positive simple infinite dimensional 
$JH^*$-triple of type I-IV,  
and $G$ is a simply connected Lie group with Lie algebra 
$\g = \hat\g(\fp)$. In view of \cite[Prop.~3.15, Thm.~5.1]{Ne02c}, 
the polar map 
\begin{equation}
  \label{eq:polar}
K \times \fp \to G, \quad 
(k,x) \mapsto k \exp x 
\end{equation}
is a diffeomorphism because the positivity of the $JH^*$-triple 
$\fp$ implies that $G/K$ is a symmetric space with seminegative 
curvature. 
Let $G_\C$ denote the simply connected Lie group 
with Lie algebra $\g_\C$ whose existence follows from Theorem~\ref{thm:exist}. 

\begin{prop} \mlabel{prop:7.2} 
Let $(G,\theta, d)$ be a simply connected 
hermitian Lie group for which $\fp$ is a positive $JH^*$-triple 
 and $G_\C$ be a simply connected group 
with Lie algebra $\g_\C$. Then the canonical morphism 
$\eta_G \: G \to G_\C$ is a covering of a closed 
subgroup $\eta_G(G)$ of $G_\C$ which is contained 
in the open subset $P^+K_c P^-$, where 
$P^\pm = \exp \fp^\pm$ and $K_c := \la \exp \fk_\C \ra$. 
\end{prop} 

\begin{prf} In view of $\z(\g) \subeq \fk$, $\eta_G(G) \subeq P^+ K_c P^-$ 
follows from the corresponding assertion for 
the group $\Ad(G) \subeq \Aut(\fp)_0$, 
for which it follows from the fact that the 
$G$-orbit of the base point in the complex homogeneous space 
$G_\C/K_c P^-$ corresponds to the domain 
\[ \cD = \{ z \in \fp^+ \: \|z\|_\infty <1 \} \subeq \fp^+ \cong P^+ \] 
(via the exponential map of $P^+$) 
and the open embedding $P^+ \into G_\C/(K_c P^-)$ 
(cf.\ \cite{Ka83}). For domains of type I-III, this also follows 
from the discussion in \cite[Sect.~III]{NO98}. 
\end{prf}

We can now refine some of the results developed in the context of 
Proposition~\ref{prop:kernel-1}. Since $G$ has a diffeomorphic 
polar map (see \eqref{eq:polar}), Proposition~\ref{prop:7.2} shows that we may 
put $U_G = W_G= G$. 
We thus obtain an analytic map 
\[ \kappa \: G \to K_\C, \quad g \mapsto \tilde\eta_G(g)_0\] 
defined on the whole group $G$. 
Let $\oline\kappa(g) := \oline{\kappa(g)}$, where 
$k \mapsto \oline k$ denotes complex conjugation of $K_\C$ with respect 
to $K$. This function plays the role of $\kappa$ if we exchange 
$P^+$ and $P^-$, which is the context of \cite[Sect.~II]{NO98}. 
We now obtain on $M  =  G/K$ a well-defined function 
\[ gK \mapsto 
(\kappa(g)^{-1})^*\kappa(g)^{-1} = \oline\kappa(g)\oline\kappa(g)^*.\] 
According to \cite[Lemma~II.3]{NO98}, the canonical 
holomorphic kernel function 
$Q_M^c \: M\times \oline M \to K_c$ satisfies 
\[Q_c(gK, gK) = \eta_K(\oline\kappa(g)\oline\kappa(g)^*)
 \quad \mbox{ for } \quad g \in G.\] 
Since $M = G/K \cong \fp$ according to the polar decomposition, 
the space $M$ is simply connected. Hence the canonical lift 
\[ Q_M \: M\times \oline M \to K_\C \] 
of $Q_c$ to the covering group  $K_\C$ of $K_c$ 
determined by $Q_M(\1K, \1 K) =\1$ satisfies
\[Q_M(gK, gK) = \oline\kappa(g)\oline\kappa(g)^* \quad \mbox{ for } 
\quad g \in G.\] 

For a bounded representation $(\rho,V)$ 
of $K$ and its holomorphic 
extension $\rho_\C \: K_\C \to \GL(V)$, we recall from 
\eqref{eq:frho-1} in Section~\ref{sec:5} the analytic function 
$f_\rho := \rho_\C \circ \kappa \: G \to \GL(V)$ 
and the corresponding holomorphic kernel function 
\[ F_\rho \: M \times \oline M \to B(V), \quad 
(gK,hK) \mapsto (f_\rho(h)^{-1})^* f_\rho(h^{-1}g) f_\rho(g)^{-1}. \] 
On the diagonal $\Delta_M \subeq M \times \oline M$ we have 
\begin{align*}
F_\rho(gK,gK) 
&= (f_\rho(g)^{-1})^* f_\rho(g^{-1}g) f_\rho(g)^{-1}
= (f_\rho(g)^{-1})^* f_\rho(g)^{-1}\\
&= \rho_\C((\kappa(g)^{-1})^* \kappa(g)^{-1})
= \rho_\C(Q_M(gK,gK)). 
\end{align*}
Therefore $F_\rho$ coincides with the kernel 
$Q_\rho := \rho_\C \circ Q_M$, because
both kernels are holomorphic on $M \times \oline M$ 
and the diagonal $\Delta_M \subeq M \times M$ is totally real. 

Proposition~\ref{prop:kernel-1} now specializes as follows: 
\begin{prop}\mlabel{prop:kernel} 
Let $(G,\theta,d)$ be a hermitian Lie group 
for which $\fp$ is a positive simple $JH^*$-triple and 
$\g_\C$ is integrable. 
Then the following are equivalent for a bounded unitary representation 
$(\rho,V)$ of $K$: 
\begin{description}
\item[\rm(i)] The representation $(\rho,V)$ is holomorphically 
inducible. 
\item[\rm(ii)] The function $f_\rho \: G \to B(V)$ is positive definite. 
\item[\rm(iii)] The kernel $Q_\rho \: M \times M \to B(V)$ is positive 
definite. 
\end{description}
\end{prop}

From now on we assume that $\g = \hat\g(\fp)$ 
is the universal central extension of 
$\g(\fp)$ (Remark~\ref{rem:3.5b}) 
and recall from the discussion in the proof of 
Theorem~\ref{thm:exist} that in all cases the corresponding 
simply connected Lie group $G$ has the property that 
$K := (G^\theta)_0$ has the product structure required in 
Theorem~\ref{thm:kredux}. 
In Remark~\ref{rem:5last} we have seen that 
$f_\rho$ is positive definite if and only if 
the corresponding function $\tilde f_\rho \: G_1 \to B(V)$ 
is positive definite, which only depends on the 
representation $\tilde\rho  \: K_1 \to \U(V)$ 
obtained by pulling $\rho$ back via the canonical map $K_1 \to K$. 

As $\fp^\pm \subeq (\g_1)_\C$, the kernel $Q_\rho$ coincides with 
the corresponding kernel that we obtain from the 
triple decomposition of the open subset 
$P^+ K_{1,c} P^-$ of $(G_1)_\C$. Therefore 
$\tilde f_\rho$ is positive definite if and only if the kernel 
$Q_\rho = Q_{\tilde \rho}$ is positive definite on $M$. 
Since the positive definite kernels $Q_{\tilde\rho}$ 
have been classified in \cite[Thm.~IV.1]{NO98}  
in terms of unitarity of related highest weight modules, 
and all bounded irreducible representations 
$(\tilde\rho,V)$ of $K_1$ can be described as highest weight 
representations, we can now derive a classification 
of the inducible irreducible bounded representations 
in terms of conditions 
on the highest weight of $\rho_0$ (cf.\ Theorem~\ref{thm:redux}). 
In the following theorem, we 
evaluate the characterizations from \cite{NO98} 
case by case for the four types of hermitian groups. 
Here we use the classification of bounded 
irreducible representations of $\U(\cH)$ whose restriction to 
$\U_\infty(\cH)$ is irreducible in terms of 
highest weights (cf.\ Definition~\ref{def:tensrep}).

\begin{thm}{\rm(Classification Theorem)} \mlabel{thm:ind-classif} 
For an irreducible bounded representation 
$(\rho,V)$ of $K$ whose restriction to $K_\infty$ is irreducible, we 
obtain the following characterization 
of the inducible representations with respect to the element 
$d = \frac{i}{2}\diag(\1,-\1)$: 
\begin{description}
\item[\rm(I$_\infty$)] Let 
$\cH_\pm$ both be infinite dimensional and 
$\hat G = \hat\U_{\rm res}(\cH_+, \cH_-)$ 
with \break $\hat K \cong \R \times \U(\cH_+) \times \U(\cH_-)$. 
Then $\rho(t,k_1, k_2) = e^{ict} \pi_{\lambda_+}(k_1) \otimes \pi_{\lambda_-}(k_2)$ is inducible if and only if 
$\lambda_+ \geq 0 \geq \lambda_-$ and 
\[ c \in |\supp(\lambda_+)| + |\supp(\lambda_-)| + \N_0.\]
\item[\rm(I$_{\rm fin}$)] Let  $\dim \cH_- < \infty$, 
$G = \tilde\U(\cH_+, \cH_-)$ 
and $K \cong \U(\cH_+) \times \tilde\U(\cH_-)$. 
Then $\rho(k_1, k_2) = 
\pi_{\lambda_+}(k_1) \otimes \pi_{\lambda_-}(k_2)$ 
is inducible if and only if $\lambda_+ \geq 0 \geq \lambda_-$ and 
\[ c := -\max(\lambda_-)\in \{ a, a+1, \ldots, b\} \cup ]b,\infty[ \subeq \Z\] 
for $a := |\supp(\lambda_+)| + |\supp(\lambda_-+c)|$ and 
$b := a - 1+ |J_-\setminus \supp(\lambda_-+c)|.$
\item[\rm(II)] For $\hat G = \hat\OO^*_{\rm res}(\cH)$ 
and $\hat K = \R \times \U(\cH)$, the representation 
$\rho(t,k) = e^{ict} \pi_{\lambda}(k)$ 
is inducible if and only if 
\[ \lambda \geq 0 \quad \mbox{ and } \quad 
c \in |\supp(\lambda)| + \N_0.\] 
\item[\rm(III)]  For $\hat G = \hat\Sp_{\rm res}(\cH)$ 
and $\hat K \cong \R \times \U(\cH)$, the representation 
$\rho(t,k) = e^{ict} \pi_{\lambda}(k)$ 
is inducible if and only if 
\[ \lambda \geq 0 \quad \mbox{ and } \quad 
 2c \in |\supp(\lambda)| + |\{ j \in J \: \lambda_j > 1\}| + \N_0.\] 
\end{description}
\end{thm}

\begin{prf} (I$_\infty$) From Example~\ref{ex:3.5} we obtain 
the structure of $\hat K$ and that $\hat G$ is a quotient of a 
semidirect product via the homomorphism 
$\gamma_G \: \tilde{SG}_1 \rtimes K \to \hat G$ 
integrating 
\[ \gamma_\g \: \sg_1 \rtimes \fk \to \hat\g, \quad 
(x,y) \mapsto (-i\tr(x_{11}), x + y).\] 

Let $(\rho, V)$ be a bounded irreducible representation 
of $\hat K$ for which the restriction to 
$\hat K_\infty = \R \times \U_\infty(\cH_+) \times \U_\infty(\cH_-)$ 
is irreducible. Then 
\[ \rho(t,k_1, k_2) = e^{ict} 
\pi_{\lambda_+}(k_1) \otimes \pi_{\lambda_-}(k_2),\] 
where $\lambda_\pm \: J_\pm \to \Z$ are finitely supported functions, 
corresponding to the highest weights 
and $(e_j)_{j \in J_\pm}$ in $\cH_\pm$ 
are orthonormal bases (cf.\ Definition~\ref{def:tensrep}). 
On 
\[ \sk_1 \cong \{ (x,y) \in \fu_1(\cH_+) \times \fu_1(\cH_-) \: 
\tr(x) + \tr(y) = 0\}\] 
this leads to a representation with highest weight 
$\lambda 
= (\lambda_+ + c, \lambda_-) = (\lambda_+, \lambda_--c),$  
where we consider $\lambda$ as a function $J = J_+ \dot\cup J_- \to \R$. 
In fact, for a diagonal matrix 
$x = \diag((x_j)_{j \in J})$ and a $\lambda$-weight vector $v_\lambda$, 
we have
\begin{align*}
&\ \ \ \ \rho(\exp(-i\tr x_{11}, x))v_\lambda 
=  \rho(-i\tr x_{11}, \exp x)v_\lambda \\
&=  e^{c \sum_{j \in J_+} x_j} 
e^{\sum_{j \in J_+}(\lambda_+)_j x_j + \sum_{j \in J_-}(\lambda_-)_j x_j}v_\lambda 
=  e^{\sum_{j \in J}  \lambda_j x_j}v_\lambda.
\end{align*}
From Lemma~\ref{lem:5.16} and Example~\ref{ex:5.12b}(b) 
we obtain the necessary condition 
$c + \min(\lambda_+)  \geq \max(\lambda_-)$. 
Since $\lambda_\pm$ have finite support, this implies that $c \geq 0$. 

To derive the classification from \cite[Prop.~I.7]{NO98}, 
we put $J_1 := J_-$, $J_2 := J_+$ and 
\[ M := \min(\lambda_+ + c) = c + \min(\lambda_+) 
\geq m := \max(\lambda_-).\] 
According to loc.\ cit.,  a necessary condition for inducibility is that 
\[ q'' := |\supp(\lambda_+ + c- M)|\quad \mbox{ and } \quad 
p'' := |\supp(\lambda_--m)| \] 
are both finite. This implies that $M = c$ and $m = 0$, so that 
\[ a := p'' + q'' = |\supp(\lambda_+)| + |\supp(\lambda_-)|. \] 
Since $J_\pm$ are both infinite, 
\[ b := a - 1+ \min\{ |J_+ \setminus \supp(\lambda_+)|, 
|J_-\setminus \supp(\lambda_-)|\} = \infty, \] 
and loc.\ cit.\ shows that inducibility is equivalent to 
$c = M- m \in a + \N_0.$ 
 
(I$_{\rm fin}$) With the same group $SG_1$ as for the previous case, 
we obtain a surjective homomorphism 
$\gamma_G \: \tilde SG_1 \rtimes K \to G$ 
integrating the summation homomorphism 
$\gamma_\g \: \sg_1 \rtimes \fk \to \g, 
(x,y) \mapsto x + y.$  
Now $K_\infty = \U_\infty(\cH_+) \times \U(\cH_-)$ 
and $(\rho, V)$ has the form 
\[ \rho(k_1, k_2) = 
\pi_{\lambda_+}(k_1) \otimes \pi_{\lambda_-}(k_2).\] 
On $\sk_1$ this leads to a representation with highest weight 
$\lambda = (\lambda_+ , \lambda_-) \: J \to \R.$ 
From Lemma~\ref{lem:5.16} and Example~\ref{ex:5.12b}(b) 
we obtain the necessary condition 
\[ M := \min(\lambda_+) \geq m := \max(\lambda_-)\] 
(see also \cite[Prop.~I.5(ii)]{NO98}). 
To derive the classification from \cite[Prop.~I.7]{NO98}, 
we put $J_1 := J_-$ and $J_2 := J_+$. 
With the finiteness of 
\[ q'' := |\supp(\lambda_+ - M)|\quad \mbox{ and } \quad 
p'' := |\supp(\lambda_--m)|, \] 
we derive that $M = 0$ and hence that $\lambda_+ \geq 0\geq \lambda_-$, 
so that 
\[ a := p'' + q'' = |\supp(\lambda_+)| + |\supp(\lambda_--m)|. \] 
Now 
\[ b 
:= a - 1+ \min\{ |J_+\setminus \supp(\lambda_+)|, 
|J_-\setminus \supp(\lambda_-)|\}  
= a - 1+ |J_-\setminus \supp(\lambda_-)| < \infty \] 
and loc.\ cit.\ shows that inducibility is equivalent to 
\[ c := - m \in \{a, a + 1, \ldots, b\} \cup ]b,\infty[. \] 

(II) From Example~\ref{ex:3.6} we recall that 
$\hat K \cong \R \times \U(\cH) = \R \times K,$  
and that $\hat G$ is obtained as a quotient of a 
semidirect product via the homomorphism 
$\gamma_G \: \tilde{G}_1 \rtimes K \to \hat G$ 
integrating 
\[ \gamma_\g \: \g_1 \rtimes \fk \to \hat\g, \quad 
(x,y) \mapsto (-i\tr(x_{11}), x + y).\] 

Let $(\rho, V)$ be a bounded irreducible representation 
of $\hat K$ for which the restriction to 
$\hat K_\infty = \R \times \U_\infty(\cH)$ 
is irreducible. Then 
$\rho(t,k) = e^{ict} \pi_{\lambda}(k)$, 
where $\lambda \: J \to \Z$ is finitely supported. 
On $\fk_1 \cong \fu_1(\cH)$, 
this leads to a representation with highest weight 
$\lambda + c$. From Lemma~\ref{lem:5.16} and the description of 
the coroots for $\Delta_p^+ \subeq \Delta_{nc}$ in 
Example~\ref{ex:d3}, we derive the necessary condition 
\[ 2 c + \lambda_j + \lambda_k \geq 0 
\quad \mbox{ for } \quad j\not= k \in J \] 
(see also \cite{NO98}). 
As $\lambda$ has finite support, we obtain $c \geq 0$. 

Comparing our $\Delta_p^+$ with the positive system used in \cite{NO98}, 
it follows that we have to apply \cite[Prop.~I.11]{NO98} to 
$-\lambda -c$. 
For 
\[M := \max(-c - \lambda) = -\min(\lambda)-c 
\quad \mbox{ we put } \quad p' := |J\setminus \lambda^{-1}(-M-c)|.\] 
Then $p'$ has  to be finite, so that 
$M = -c  \leq 0$, $\lambda \geq 0$  
and $p' = |\supp(\lambda)|$. 
According to \cite[Prop.~I.11]{NO98}, inducibility is 
equivalent to 
$c \in |\supp(\lambda)| + \N_0.$ 

(III) Here we  have a similar situation as for type II. 
Now $\rho(t,k) = e^{ict} \pi_{\lambda}(k)$ 
leads to a representation with highest weight  
$\lambda + c$ of $\fk_1 \cong \fu_1(\cH)$ 
and Lemma~\ref{lem:5.16} with Example~\ref{ex:d2} 
 lead to the necessary condition 
\[ 2 c + \lambda_j + \lambda_k \geq 0 
\quad \mbox{ for } \quad j,k \in J, \] 
which is equivalent to $c + \lambda \geq 0$. 
As $\lambda$ has finite support, this implies $c \geq 0$. 

For $M := \max(-c - \lambda) = -c -\min(\lambda)$ we put 
\[ q' := |\{ j \in J \: \lambda_j > \min(\lambda)\}| \quad \mbox{ and } \quad 
r' := |\{ j \in J \: \lambda_j > \min(\lambda)+1\}|.\] 
Then $r'$ and $q'$ have to be finite, so that 
$M = -c$, $\min\lambda = 0$ and  
$q' = |\supp(\lambda)|$. 
According to \cite[Prop.~I.9]{NO98}, inducibility is 
equivalent to $2c \in q' + r' + \N_0.$ 
\end{prf}

\begin{rem} (a) For the types I$_\infty$, II and III, 
the number $c$ is called the {\it central charge} 
of the representation. The preceding theorem shows that $c = 0$ implies 
$\lambda = 0$, hence that the irreducible representation $\rho$ is trivial. 
In particular, the groups 
$\U_{\rm res}(\cH_+, \cH_-)$, 
$\OO^*_{\rm res}(\cH)$ and $\Sp_{\rm res}(\cH)$ have no 
non-trivial irreducible semibounded representation 
$(\pi,\cH)$ with $\fp \not\subeq \ker \dd\pi$. 
All irreducible semibounded representations of these groups 
are pullbacks of bounded representations of $K/K_\infty$ by 
the canonical homomorphism $G \to K/K_\infty$. 

For $K = \U(\cH)$, the quotient 
$K/K_\infty \cong \U(\cH)/\U_\infty(\cH)$ is the identity component 
of the unitary group of the $C^*$-algebra $B(\cH)/K(\cH)$, hence 
has enough bounded unitary representations to separate the points. 
However, these representations all live on non-separable spaces 
(see \cite{Pi88} and  Theorem~\ref{thm:8.1b} below). 

(b) The representations for which $V$ is one dimensional 
are called {\it of scalar type}. For the groups  
$\hat\U_{\rm res}(\cH_+, \cH_-)$ and 
$\hat\OO^*_{\rm res}(\cH)$ they 
are parameterized by $c \in \N_0$, and for 
$\hat\Sp_{\rm res}(\cH)$ we have $c \in \shalf\N_0$. 

For $c = \shalf$ we obtain in particular the 
metaplectic representation of 
$\hat\Sp_{\rm res}(\cH)$ on 
the even Fock space $S^{\rm even}(\cH)$ 
and for $\rho(t,k) = e^{it/2}k$ on $V = \cH$ we 
obtain the metaplectic representation 
on the odd Fock space $S^{\rm odd}(\cH)$ 
(cf.\ \cite{Ot95}, \cite[Sect.~9.1]{Ne10c}, \cite[Sect.~IV]{NO98}). 

(c) For type I$_{\rm fin}$ ($\dim \cH_- < \infty$), 
the scalar type condition means that 
$\lambda_+ = 0$ and $\lambda_- = -c$, 
so that $a = 0$ and  $b = |J_-|-1$ 
lead to the condition 
\[ c \in \{0,1,2, |J_-|-1\} \quad \mbox{ or } \quad 
c > |J_-| -1.\] 

For type $I_{\rm fin}$ and $\sum_j (\lambda_-)_j = - \sum_j (\lambda_+)_j$,  
the representation $\rho  =\pi_{\lambda_+} \otimes \pi_{\lambda_-}$ 
vanishes on 
\[ Z(\tilde G)_0 = \{ (e^{-it}, t) \: t\in \R\} \subeq 
K \cong \U(\cH_+) \times (\R \ltimes \SU(\cH_-)),\] 
so that we also obtain non-trivial inducible representations 
vanishing on the center. These lead to semibounded representations 
of $\U(\cH_+, \cH_-)/\T \1$. 
\end{rem}

For the domains of type IV, the situation is quite trivial, as the 
following theorem shows. 

\begin{thm} The universal covering group 
$\tilde\OO(\R^2, \cH_\R)$ of $\OO(\R^2, \cH_\R)$ has no 
semibounded unitary representation $(\pi,\cH_0)$ 
with $\fp \not\subeq \ker(\dd\pi)$. 
\end{thm}

\begin{prf} Put $G := \OO(\R^2, \cH_\R)$ and write 
$\tilde G$ for its universal covering group so that 
$K \cong \SO_2(\R) \times \OO(\cH_\R)$ and 
$\tilde K \cong \R \times \OO(\cH_\R) \subeq \tilde G$. 

Let $(\rho, V)$ be an irreducible holomorphically inducible bounded 
representation of $\tilde K$. 
Then $\rho$ is of the form 
$\rho(t,k) = e^{ict} \rho_0(k)$, where 
$(\rho_0, V)$ is an irreducible representation 
of $\OO(\cH_\R)$. 
In view of Theorem~\ref{thm:kredux}, it suffices to show that, 
if $\rho\res_{\OO_\infty(\cH_\R)}$ is irreducible and non-trivial, then 
$(\rho,V)$ is not inducible. 

According to Theorem~\ref{thm:oprep}, 
every irreducible bounded representation 
$(\rho,V)$ of $\OO_\infty(\cH_\R)$ is a highest weight 
representation $(\rho_\lambda,V)$ with finitely supported highest 
weight. Now \cite[Sect.~1, Thm.~IV.1]{NO98} implies that 
the kernel $Q_{\rho_\lambda}$ on $G/K$ is not positive definite 
for $\lambda\not=0$, so that $(\rho_\lambda, V)$ is not inducible. 
\end{prf}

\begin{rem} The results of the present section imply 
in particular that the unitary representations of 
the groups $G_1$ constructed in \cite{NO98} 
all extend to unitary representations 
of the full hermitian groups 
$\hat\U_{\rm res}(\cH_+, \cH_-)$, 
$\hat\OO^*_{\rm res}(\cH)$ and 
$\hat\Sp_{\rm res}(\cH)$. 
\end{rem}

\section{Hermitian groups with negative $JH^*$-triples} 
\mlabel{sec:8}

In this section we consider the case where 
$\fp$ is an infinite dimensional negative simple $JH^*$-triple, i.e., 
$G/K$ is $c$-dual to a symmetric Hilbert domain of type I-IV, i.e., 
a hermitian symmetric space of ``compact type'' in the sense 
of Kaup (cf.\ \cite{Ka81,Ka83}). 
We assume that $G$ is a simply connected Lie group with Lie algebra 
$\g = \hat\g(\fp)$ (Remark~\ref{rem:3.5b}) and recall that 
this implies that $K$ has the product structure required in 
Theorem~\ref{thm:kredux}.

\begin{thm} \mlabel{thm:8.1} Suppose that $G$ is full and simply connected 
and that $\fp$ is an infinite dimensional negative $JH^*$-triple. 
Then the following are equivalent for a bounded unitary irreducible 
representation $(\rho, V)$ of $K$: 
\begin{description}
\item[\rm(i)] $(\rho, V)$ is holomorphically inducible. 
\item[\rm(ii)] $\dd\rho([z^*,z]) \geq 0$ for $z \in \fp^+$ 
($\rho$ is anti-dominant). 
\end{description}
\end{thm}

\begin{prf} (i) $\Rarrow$ (ii) follows from Lemma~\ref{lem:neccond}. 

(ii) $\Rarrow$ (i): In view of Theorem~\ref{thm:redux}, we 
may w.l.o.g.\ assume that $\rho\res_{K_\infty}$ is irreducible, 
i.e., $\rho = \rho_0$. We have seen in Remark~\ref{rem:3.9} 
that, for all types, we have a surjective submersion 
$G_1 \rtimes K \to G$. In view of Remark~\ref{rem:5last}, 
it suffices to show that 
the function $f_\rho$ is positive definite 
in a $0$-neighborhood of the group $G_1$ if (ii) is satisfied. 
To verify this, we  consider all types separately, where for 
type I-III we take 
$d = \frac{i}{2}\diag(\1,-\1)$. 

(I$_\infty$) Here $G = \hat\U_{\rm res}(\cH_+ \oplus \cH_-)$ 
and $K = \T \times \U(\cH_+) \times \U(\cH_-)$ (Example~\ref{ex:3.5}(c)),  
so that $\rho$ has the form 
$\rho(z,k_1, k_2) = z^c \rho_{\lambda_+}(k_1) \otimes 
\rho_{\lambda_-}(k_2)$,  
where $c \in \Z$ and $\lambda_\pm \: J_\pm \to \Z$ 
are finitely supported functions, 
corresponding to the highest weights,  
and $(e_j)_{j \in J_\pm}$ in $\cH_\pm$ 
are orthonormal bases. 
From 
\[ G_1 = \SU_{1,2}(\cH_+ \oplus \cH_-), \quad 
K_1 = S(\U_{1}(\cH_+) \times \U_1(\cH_-))
\cong \T \rtimes (\SU(\cH_+) \times \SU(\cH_-))\] 
(cf.\ Example~\ref{ex:3.5} and Appendix~\ref{app:d.4}), 
we obtain on $K_1$ a representation with the bounded highest 
weight $\lambda = (c + \lambda_+, \lambda_-)$. 
Condition (ii) now leads with 
Lemma~\ref{lem:5.16} and \eqref{eq:weight1} in Example~\ref{ex:5.12b}(a) to 
\begin{equation}
  \label{eq:condifinfin}
c + \max(\lambda_+)  \leq \min(\lambda_-).
\end{equation}
Hence there exists a linear order $\preceq$ on $J = J_+ \dot\cup J_-$ 
with $J_+ \prec J_-$ for which $\lambda \: J \to \Z$ is increasing. 
Let $(\pi_\lambda, \cH_\lambda)$ denote the 
corresponding bounded highest weight representation 
of $\SU(\cH)$ whose existence follows from \cite[Prop.~III.7]{Ne98}. 

In view of \eqref{eq:condifinfin}, 
any highest weight vector $v_\lambda$ generates the $K_1$-representation 
$(\rho_\lambda,V_\lambda)$ of highest weight $\lambda$ 
annihilated by $\tilde\fp^- := \fp^- \cap \gl_1(\cH)$. 
For the orthogonal projection 
$p_V \: \cH_\lambda \to V_\lambda$ we therefore find 
\[ f_\rho(g) = p_V \pi_\lambda(g) p_V \in B(V_\lambda)\] 
because this relation holds for $g \in K_1 \subeq \SU(\cH)$ and 
both sides are annihilated by the differential operators 
$L_z$, $z\in \tilde\fp^-$. We conclude that the function $f_\rho$ is  
positive definite on the dense subset 
$W_G W_G^{-1} \cap \SU(\cH)$ of $G_1$, so that its continuity 
implies that it is positive definite. 

(I$_{\rm fin}$) Here we have $G = \U(\cH_+ \oplus \cH_-)$ 
and $K = \U(\cH_+) \times \U(\cH_-),$ 
so that $\rho$ has the form 
$\rho(k_1, k_2) = \rho_{\lambda_+}(k_1) \otimes \rho_{\lambda_-}(k_2)$ 
and $\lambda = (\lambda_+, \lambda_-)$ has finite support. 
In view of Lemma~\ref{lem:5.16} and \eqref{eq:weight1} in 
Example~\ref{ex:5.12b}(b), 
condition (ii) here means that 
\begin{equation}
  \label{eq:condifin}
\max(\lambda_+)  \leq \min(\lambda_-).
\end{equation}
Hence there exists a linear order $\preceq$ on $J = J_+ \dot\cup J_-$  
with $J_+ \prec J_-$ for which $\lambda$ is increasing. 
Let $(\pi_\lambda, \cH_\lambda)$ denote the 
corresponding bounded highest weight representation 
of $G$ (cf.\ Appendix~\ref{app:e}). 
As a consequence of \eqref{eq:condifin}, any highest weight vector 
$v_\lambda$ generates the $K$-representation 
$(\rho_\lambda,V_\lambda)$ of highest weight $\lambda$, 
and the subspace $V_\lambda$ is annihilated by $\fp^-$.
As above, it now follows that 
\[ f_\rho(g) = p_V \pi_\lambda(g) p_V\quad \mbox{ for } \quad g \in G, \] 
and hence that $\rho$ is inducible. 

(II) Here $G = \hat\OO_{\rm res}(\cH_\R)$ 
(cf.\ Example~\ref{ex:d3} and Appendix~\ref{app:d.3}) 
and $K \cong \T \times \U(\cH)$ (cf.\ Example~\ref{ex:3.6}),  
so that $\rho$ has the form 
$\rho(z,k) = z^c \rho_{\mu}(k)$ 
with $c \in \Z$ and highest weight $\mu\: J\to \Z$. 
From $G_1 = \tilde\OO_{1,2}(\cH_\R)$ and the fact that 
$K_1 \cong  \hat \U_1(\cH)$ is the unique $2$-fold covering 
of $\U_1(\cH) \cong \T \ltimes \SU(\cH)$, 
we obtain on $K_1$ the representation with the bounded highest 
weight $\lambda := c/2 + \mu$ (cf.\ Example~\ref{ex:3.6}). 
With Lemma~\ref{lem:5.16} and Example~\ref{ex:d3} 
condition (ii) now translates into 
\begin{equation}
  \label{eq:condii}
\lambda_j + \lambda_k =  c + \mu_j + \mu_k \leq 0 
\quad \mbox{ for } \quad j\not= k \in J. 
\end{equation}
As $\mu$ is integral with finite support, we get $c \leq 0$. 
We pick a linear order on $J$ 
for which $\lambda \: J \to \Z$ is increasing. 
From \cite[Sect.~VII]{Ne98} 
we obtain a bounded highest weight representation 
$(\pi_\lambda, \cH_\lambda)$ of 
$\tilde\OO_1(\cH_\R)$, resp., a holomorphic representation of 
its universal complexification  
$\tilde\OO_1(\cH_\R)_\C \cong \OO_1(\cH^2,\beta)$ 
(cf.~Appendix~\ref{app:d.1}). 
The remaining argument is similar as for 
type I$_\infty$ and uses the density of 
$\tilde\OO_1(\cH_\R)$ in $G_1 = \tilde\OO_{1,2}(\cH_\R)$. 

(III) Here $G = \hat\Sp_{\rm res}(\cH_\H) 
\subeq G_\C \cong \hat\Sp_{\rm res}(\cH^2,\omega)$ 
(cf.\ Example~\ref{ex:d2}), 
and the argument is quite similar to type III. 
Condition (ii) translates into 
\[ \lambda_j = c + \mu_j \leq 0 
\quad \mbox{ for } \quad j\in J, \] 
$c \leq 0$, and  $\lambda_j \in \shalf \Z$ for every $j \in J$. 
Since $\lambda$ is bounded, we can argue with 
a holomorphic highest weight representation 
of the group $\Sp_1(\cH^2,\omega)$ whose existence follows 
from \cite[Sect.~VI]{Ne98}. 

(IV) Here $G = \OO(\R^2 \oplus \cH_\R)$ is a full orthogonal 
group and 
\[ K 
= \SO_2(\R) \times \OO(\cH_\R) \cong \T \times \OO(\cH_\R), \] 
so that $\rho$ has the form 
$\rho(z,k) = z^c \rho_{\mu}(k)$ 
with $c \in \Z$ and a bounded highest weight 
representation $\rho_\mu$ of $\OO(\cH_\R)$. 
Let $J$ parameterize an orthonormal basis of 
$\C \oplus (\cH_\R,I)$, where $I$ is a complex structure on 
$\cH_\R$, such that $j_0 \in J$ corresponds to the basis element 
$e_{j_0} := (1,0)$. Then $\lambda_{j_0} := c$ and 
$\lambda\res_{J\setminus\{j_0\}} := \mu$ defines a function 
$\lambda \: J \to \Z$. 
With Lemma~\ref{lem:5.16} and Example~\ref{ex:d4}, we see that 
condition (ii) means that 
\begin{equation}
  \label{eq:typeivcond}
c\pm \mu_j \geq 0 \quad \mbox{ for } \quad j \not=j_0,
\end{equation}
i.e., $c \geq 0$ and $|\mu_j| \leq c$ (cf.\ \cite[Sect.~I]{NO98}). 
From \cite[Sect.~VII]{Ne98} we now obtain a 
bounded highest weight representation 
$(\pi_\lambda, \cH_\lambda)$ of $G$. 
As a consequence of \eqref{eq:typeivcond}, any highest weight vector 
$v_\lambda$ generates the $K$-representation 
$(\rho_\lambda,V_\lambda)$ of highest weight $\lambda$, 
and the subspace $V_\lambda$ is annihilated by $\fp^-$.
With similar arguments as for type I$_{\rm fin}$, it now follows that 
\[ f_\rho(g) = p_V \pi_\lambda(g) p_V \in B(V_\lambda)\] 
is positive definite, so that $\rho$ is inducible. 
\end{prf}

\begin{rem} \mlabel{rem:8.2} Suppose that, 
in the preceding proof, 
for one of the groups 
\[ \hat\U_{\rm res}(\cH), \quad\hat\OO_{\rm res}(\cH_\R), \quad 
\hat\Sp_{\rm res}(\cH_\H)\] 
we have $c=0$. Then the corresponding highest weight 
$\lambda$ is finitely supported and 
we derive from \cite[Sects.~V-VII]{Ne98} that the corresponding 
unitary highest weight representation 
of the corresponding non-extended group 
\[ \U_{\rm res}(\cH), \quad \OO_{\rm res}(\cH_\R), \quad
\Sp_{\rm res}(\cH_\H)\] 
extends to a bounded representation of the corresponding full group 
\[ \U(\cH),\quad \OO(\cH_\R), \quad \Sp(\cH_\H).\] 
In particular, the corresponding 
highest weight representation is bounded. 

Suppose, conversely, that $c \not=0$. We want to show that 
the corresponding highest weight representation 
$\pi_\lambda$ is unbounded. To this end, it suffices to show that 
for the corresponding set $\cP_\lambda$ of weights, the set of all 
values $\cP_\lambda(d)$ is infinite. 

For type I$_\infty$ we have 
$\lambda = (c + \lambda_+, \lambda_-) \: J \to \Z$ 
and the group $\cW = S_{(J)}$ of all finite permutations of $J$ 
acts on the weight set $\cP_\lambda$ with respect to $\ft_1 \subeq \fk_1$. 
Let $F \subeq J_+$ be a finite subset not contained in the support 
of $\lambda_+$ and 
$w \in \cW$ be an involution with 
$w(F_+) \subeq J_- \setminus \supp(\lambda_-)$ 
and $w(j) = j$ for $j \not\in F_+ \cup w(F_+)$. 
Then 
\[ \lambda - w\lambda= c\sum_{j \in F_+} (\eps_j - \eps_{w(j)}) \] 
has on $-i\cdot d$ the value  $c |F_+|$, which can be arbitrarily 
large. 

For type II we have 
$\lambda = c/2 + \mu$ and the weight set 
is invariant under the Weyl group 
$\cW$ which contains in particular sign changes 
$s(h_j) = (\eps_j h_j)$, where 
$\eps_j \in \{\pm 1\}$, and 
$F = \{ j \in J \: \eps_j = -1\}$ can be any finite subset of $J$ 
with an even number of elements (cf.\ Example~\ref{ex:d.1c}). 
For $F \cap \supp(\mu) = \eset$, 
we then obtain 
\[ (\lambda -w\lambda)(-i\cdot d) = c |F|,\] 
which can be arbitrarily large. 

For type III we have 
$\lambda = c + \mu$ and a similar  argument applies.
\end{rem}

\begin{ex} (Finite type I) 
If one of the spaces $\cH_+$ or $\cH_-$ in $\cH = \cH_+ + \cH_-$ 
is finite dimensional, 
then $G = \U_{\rm res}(\cH) = \U(\cH)$ is the full unitary 
group. Since every open invariant cone in $\fu(\cH)$ intersects 
the center (Theorem~\ref{thm:coneinuh}), all irreducible semibounded representations 
of $G$ are bounded (Proposition~\ref{prop:4.5}). 

According to Theorem~\ref{thm:kredux}, every 
bounded irreducible representation 
$(\rho, V)$ of $\U(\cH)$ is a tensor product 
of a representation 
$(\rho_0,V_0)$ whose restriction to 
$\U_\infty(\cH)$ is irreducible and of a representation 
$(\rho_1, V_1)$ with $\U_\infty(\cH) \subeq \ker \rho_1$. 
The representations $\rho_0$ are easily classified 
by their highest weights, as described in Appendix~\ref{app:e} 
but for the representations $\rho_1$ there is no concrete 
classification available (cf.\ Remark~\ref{rem:5.14}). 
\end{ex}

\begin{ex} (Infinite type I) 
Consider the group $\hat G = \hat\U_{\rm res}(\cH)$ 
with $\dim \cH_+ = \dim \cH_- = \infty$. 
For representations which are trivial on the center, i.e., 
$c =0$, we obtain for $\lambda = (\lambda_+, \lambda_-)$ 
the necessary condition 
$\lambda_+ \leq 0 \leq  \lambda_-$, which leads to a 
representation of $K \cong \U(\cH_+) \times \U(\cH_-)$ on a subspace
\[ \cH_{\lambda_+}\otimes \cH_{\lambda_-}
\subeq (\cH_+^*)^{\otimes n} \otimes \cH_-^{\otimes m} 
\subeq (\cH^*)^{\otimes n} \otimes \cH^{\otimes m}\] 
(cf.\ Appendix~\ref{app:e}).

On the latter space we have a unitary representation of 
$\U(\cH)$. We claim that it contains a subrepresentation 
which is holomorphically induced from the representation 
of $\U(\cH_+) \times \U(\cH_-)$ on 
$\cH_{\lambda_+}\otimes \cH_{\lambda_-}$. 
In fact, the natural action of $\fp^- \cong B_2(\cH_+, \cH_-)$ 
annihilates the subspace 
$\cH_{\lambda_+}\subeq (\cH_+^*)^{\otimes n}\subeq (\cH^*)^{\otimes n}$ 
and likewise, it annihilates the subspace 
$\cH_{\lambda_-}
\subeq \cH_-^{\otimes m} \subeq \cH^{\otimes m}$. 
Now the assertion follows from \cite[Cor.~3.9]{Ne10d}. 
\end{ex}

\appendix

\newcommand{\ruler}{\hbox{\vrule width0.5pt height 5mm depth 3mm}}

\section{Analytic operator-valued positive definite functions} 
\mlabel{app:a}

In this appendix we discuss operator-valued positive 
definite functions on Lie groups. The main result 
is Theorem~\ref{thm:extension}, asserting that local analytic 
positive definite functions extend to global ones. 
This generalizes the corresponding result for the scalar case 
in \cite{Ne10b}. 

\begin{defn} \mlabel{def:a.1} 
Let $X$ be a set and $\cK$ be a Hilbert space.

\par (a)  A function $Q \: X \times X \to B(\cK)$ is called a 
{\it $B(\cK)$-valued kernel}. 
It is said to be {\it hermitian} if 
$Q(z,w)^* = Q(w,z)$ holds for all $z, w \in X$. 

\par (b) A hermitian $B(\cK)$-valued kernel $K$ on $X$ is said to be 
{\it positive definite} if 
for every finite sequence $(x_1, v_1), \ldots, (x_n,v_n)$ in $X \times \cK$
we have 
\[ \sum_{j,k = 1}^n \la Q(x_j, x_k)v_k, v_j \ra \geq 0. \] 

\par (c) If $(S,*)$ is an involutive semigroup, then a 
function $\phi \: S \to B(\cK)$ is called {\it positive definite} 
if the kernel $Q_\phi(s,t) := \phi(st^*)$ is positive definite. 

\par (d) Positive definite kernels can be characterized as those 
for which there exists a Hilbert space $\cH$ and a 
function $\gamma \: X \to B(\cH,\cK)$ such that 
\[ Q(x,y) = \gamma(x)\gamma(y)^* \quad \mbox{ for } \quad x,y \in X \] 
(cf.\ \cite[Thm.~I.1.4]{Ne00}). 
Here one may assume that the vectors 
$\gamma(x)^*v$, $x \in X, v \in \cK$, span a dense subspace of 
$\cH$. Then the pair $(\gamma,\cH)$ is called a {\it realization 
of $K$}. 
The map $\Phi \: \cH \to \cK^X, \Phi(v)(x) := \gamma(x)v$, 
then realizes $\cH$ as a Hilbert subspace of $\cK^X$ 
with continuous point evaluations $\ev_x \: \cH \to \cK$. 
It is the unique Hilbert subspace in $\cK^X$ with this property 
for which $Q(x,y) = \ev_x \ev_y^*$ for $x,y \in X$. 
We write $\cH_Q \subeq \cK^X$ for this subspace and call it 
the {\it reproducing kernel Hilbert space with kernel~$Q$}.
\end{defn}

\begin{rem} \mlabel{rem:a.2} 
Let $Q \: X \times X \to B(\cK)$ 
and $f \: X \to \GL(\cK)$ be functions. It is obvious that the kernel 
$Q$ is positive definite if and only if the modified kernel 
\[ (x,y) \mapsto f(x) Q(x,y) f(y)^* \] 
is positive definite. 
\end{rem}

\begin{thm} \mlabel{thm:kerana} Let $M$ be an analytic Fr\'echet manifold and 
$\cH$, $\cK$ be Hilbert spaces. Then a function $\gamma \: M \to B(\cH,\cK)$ 
is analytic if and only if the kernel 
$Q(x,y) := \gamma(x)\gamma(y)^* \in B(\cK)$ is analytic on an open 
neighborhood of the diagonal in $M \times M$. 
\end{thm}

\begin{prf} Since the map 
$B(\cH,\cK) \times B(\cH,\cK) \to B(\cK), (A,B) \mapsto AB^*$ 
is continuous and real bilinear, it is real analytic. 
Therefore the analyticity of $\gamma$ implies that $Q$ is analytic. 

Conversely, we assume that $Q$ is analytic. 
On the analytic Fr\'echet manifold $X := M \times \cK$ we obtain 
the positive definite kernel 
\[ \tilde Q((m,v), (m',v')) 
:= \la Q(m,m')v', v\ra 
= \la \gamma(m)\gamma(m')^*v', v\ra
= \la \gamma(m')^*v', \gamma(m)^*v\ra  \] 
which is analytic by assumption. 
Therefore \cite[Thm.~5.1]{Ne10b} implies that the map 
\[ \Gamma \: M \times \cK \to \cH, \quad (m,v) \mapsto 
\gamma(m)^*v \] 
is analytic. 

Since the assertion we have to prove is local, we may 
w.l.o.g.\ assume that $M$ is an open subset of a 
real Fr\'echet space $V$. Pick $x_0 \in M$. Then the definition of 
analyticity implies the existence of an 
open neighborhood $U_V$ of $x_0$ in the complexification $V_\C$ with 
$U_V \cap V \subeq M$ and an open $0$-neighborhood 
$U_\cK \subeq \cK$ such that 
$\Gamma$ extends to a holomorphic function 
\[ \Gamma_\C \: U_V \times U_\cK \to \cH.\] 
Then $\Gamma_\C(m,v)$ is linear in the second argument, hence of the 
form $\Gamma_\C(m,v) = \gamma_\C(m)^* v$, where 
$\gamma_\C(m)^*\in B(\cK,\cH)$. 
From $\|\gamma(m)\|^2 = \|Q(m,m)\|$ we further derive that 
$\gamma$ is locally bounded, so that we may w.l.o.g.\ assume that 
$\gamma_\C$ is bounded. Now \cite[Prop.~I.1.9]{Ne00} implies 
that $\gamma_\C \: U_V \to B(\cK)$ is holomorphic, and hence that 
$\gamma$ is real analytic in a neighborhood of~$x_0$.
\end{prf}

\begin{defn} Let $\cK$ be a Hilbert space, $G$ be a group, and 
$U \subeq G$ be a subset. A function $\phi \: UU^{-1} \to B(\cK)$ 
is said to be {\it positive definite} if the kernel 
\[ Q_\phi \: U \times U \to B(\cK), \quad (x,y) \mapsto \phi(xy^{-1})\] 
is positive definite. 
\end{defn}

\begin{defn}
A Lie group $G$ is said to be 
{\it locally exponential} 
if it has an exponential function $\exp \: \g=\L(G) \to G$ 
for which there is an open $0$-neighborhood 
$U$ in $\L(G)$ mapped diffeomorphically onto an 
open subset of $G$. If, in addition, $G$ is analytic and 
the exponential function is an analytic local diffeomorphism in $0$, 
then $G$ is called a {\it BCH--Lie group}. 
Then the Lie algebra $\L(G)$ is a 
{\it BCH--Lie algebra}, i.e., there exists an open 
$0$-neighborhood $U \subeq \g$ such that, for $x,y \in U$, the BCH series 
$$x * y = x + y + \frac{1}{2}[x,y] + \cdots  $$
converges and defines an analytic function 
$U \times U \to \g, (x,y) \mapsto x * y$. 
The class of BCH--Lie groups contains in particular all Banach--Lie groups 
(\cite[Prop.~IV.1.2]{Ne06}). 
\end{defn}

\begin{defn} \label{def:6.5} Let $(\rho, \cD)$ be a $*$-representation 
of $U_\C(\g)$ on the pre-Hilbert space~$\cD$, i.e., 
\[ \la \rho(D)v,w\ra = \la v, \rho(D^*)w \ra \quad \mbox{ for } \quad 
D \in U_\C(\g), v,w \in \cD.\] We call 
a subset $E \subeq \cD$ {\it equianalytic} if 
there exists a $0$-neighborhood 
$U \subeq \g$ such that 
$$ \sum_{n = 0}^\infty \frac{\|\rho(x)^n v\|}{n!} < \infty 
\quad \mbox{ for } \quad v \in E, x \in U. $$
This implies in particular that each $v \in E$ is an analytic 
vector for every $\rho(x)$, $x \in \g$.
\end{defn}

\begin{thm} {\rm(Extension of local positive definite analytic functions)} 
\mlabel{thm:extension} 
Let $G$ be a simply connected Fr\'echet--BCH--Lie group, 
$V \subeq G$ an open connected $\1$-neighborhood, 
$\cK$ be a Hilbert space and 
$\phi \: VV^{-1} \to B(\cK)$ be an analytic positive definite function.  
Then there exists a unique analytic positive definite 
function $\tilde\phi \: G \to B(\cK)$ extending~$\phi$. 
\end{thm}

\begin{prf} The uniqueness of $\tilde\phi$ follows from the connectedness 
of $G$ and the uniqueness of analytic continuation.

{\bf Step 1:} To show that $\tilde\phi$ exists, 
we consider the reproducing kernel 
Hilbert space $\cH_Q \subeq \cK^V$ defined by the kernel $Q$ via 
$f(g) = Q_g f$ for $g \in V$ and $Q(h,g) = Q_h Q_g^*$ for $g,h \in V$ 
(Definition~\ref{def:a.1}(d)). 
Then the analyticity of the function 
$$ V \times V \to B(\cK), \quad (g,h) \mapsto Q(g,h) = \phi(gh^{-1}) $$
implies that the map $\eta \: V \to B(\cH_Q,\cK), g \mapsto Q_g$ 
is analytic (Theorem~\ref{thm:kerana}). Here we use that 
$G$ is Fr\'echet. Hence all functions in
$\cH_Q$ are analytic, so that we obtain for each 
$x \in \g$ an operator 
\[ L_x  \: \cH_Q \to C^\omega(V,\cK), \quad 
(L_x f)(g) 
:= \derat0 f(g\exp(tx)) = \derat0 Q_{g\exp(tx)}f, \quad g \in V.\] 
For $v \in \cK$, we then have 
\begin{align*}
(L_x Q_h^* v)(g)  
&= \derat0 Q_{g\exp(tx)}Q_h^*v 
= \derat0 \phi(g\exp(tx)h^{-1})v\\
&= \derat0 Q_g Q_{h\exp(-tx)}^*v
= \derat0 (Q_{h\exp(-tx)}^*v)(g) 
\end{align*}
for $g,h \in V, x \in \g$, which means that 
\[  L_x Q_h^* = \derat0 Q_{h\exp(-tx)}^* \in B(\cK,\cH_Q).\] 
Iterating this argument, we see by induction that, 
for $x_1,\ldots, x_n \in \g$, 
\begin{equation}
  \label{eq:diffrel}
L_{x_1} \cdots L_{x_n} Q_h^* 
= \frac{\partial^n}{\partial t_1 \cdots \partial t_n}
\ruler_{\ t_1 = \ldots = t_n=0} 
Q_{h\exp(-t_nx_n)\cdots \exp(-t_1x_1)}^* 
\end{equation}
defines a bounded operator $\cK \to \cH_Q$. 

{\bf Step 2:} For an open subset $W \subeq V$, we thus obtain the subspace 
\[ \cD(W) := \Spann\{ L_{x_1} \cdots L_{x_n} Q_h^*v \: 
v\in \cK, h \in W, x_1,\ldots, x_n \in \g,n \in \N_0\} \] 
of $\cH_Q$ and operators 
\[ \rho(x) := L_x\res_{\cal D} \: \cD(V) \to \cD(V), \quad x \in \g,\] 
defining a representation $\rho \: 
\g \to \End(\cD(V))$.\begin{footnote}
{In \cite[Thm.~7.3]{Ne10b} we prove a version of the present theorem 
for the special case $\cK = \C$. In loc.cit.\ we claim 
that $\cH_Q^0$ is invariant under the Lie algebra $\g$, but this need not 
be the case. The argument given here, where $\cH_Q^0$ is enlarged to 
$\cD(V)$ corrects this point. }  
\end{footnote}
Here we use that $\cD(V) \subeq C^\omega(V,\cK)$ and the fact that 
$\g$ acts by left invariant 
vector fields on this space. Next 
we observe that 
\begin{align*}
&\ \ \ \la L_{x_1} \cdots L_{x_n} Q_h^*v, Q_g^* w \ra  \\
&= \frac{\partial^n}{\partial t_1 \cdots \partial t_n}
\ruler_{\ t_1 = \ldots = t_n=0} 
\la Q_gQ_{h\exp(-t_nx_n)\cdots \exp(-t_1x_1)}^*v,w\ra \\
&= \frac{\partial^n}{\partial t_1 \cdots \partial t_n}
\ruler_{\ t_1 = \ldots = t_n=0} 
\la \phi(g\exp(t_1 x_1) \cdots \exp(t_nx_n)h^{-1})v,w\ra \\ 
&= \frac{\partial^n}{\partial t_1 \cdots \partial t_n}
\ruler_{\ t_1 = \ldots = t_n=0} 
\la Q_h^*v, Q_{g\exp(t_1 x_1) \cdots \exp(t_nx_n)}^*w \ra\\
&= \la Q_h^*v, (-L_{x_n})\cdots (-L_{x_1})Q_g^*w\ra.
\end{align*}
Therefore $\rho(x) \subeq - \rho(x)^*$, and thus $\rho$ 
extends to a $*$-representation of $U_\C(\g)$ on~$\cD(V)$. 

{\bf Step 3:} Since $\eta$ is analytic, we derive from \eqref{eq:diffrel} 
for each $g \in V$ and sufficiently small $x \in \g$, the relation 
\begin{equation}
  \label{eq:step3}
Q_{g \exp(-x)}^* = \sum_{n = 0}^\infty \frac{1}{n!} 
\rho(x)^n Q_g^*. 
\end{equation} 
Since derivatives of analytic functions are also analytic, 
using \eqref{eq:diffrel} again implies that 
$\cD(V)$ consists of analytic vectors for the representation 
$\rho$ of~$\g$.

{\bf Step 4:} Let $W_G \subeq V$ be an open $\1$-neighborhood 
and $W_\g \subeq \g$ an open balanced $0$-neighborhood 
with $W_G \exp(W_\g) \subeq V$. Next we show that 
$\cD(W_G)$ is equianalytic and spans a dense subspace of $\cH_Q$. 

Since the map $W_G \times W_\g\to B(\cK,\cH_Q), (g,x) \mapsto 
Q_{g \exp(x)}^*$ is analytic, \cite[Lemma~7.2]{Ne10b} shows that, 
after shrinking $W_G$ and $W_\g$, we may assume that 
\begin{equation}
  \label{eq:12}
Q_{g \exp (-x)}^* = \sum_{n = 0}^\infty 
\frac{1}{n!} \rho(x)^n Q_g^* 
\quad \mbox{ for } \quad g \in W_G, x \in W_\g.
\end{equation}
This implies that $\cD(W_G)$ is equianalytic. 
To see that $\cD(W_G)$ is dense, 
we use the analyticity of $\eta$ to see that 
$\eta(V)^*\cK \subeq \oline{\cD(W_G)}$, which implies that 
$\cD(W_G)$ is dense in $\cH_Q$. 

Applying \cite[Thm.~6.8]{Ne10b} to the $*$-representation 
$(\rho,\cD(W_G))$, 
we now obtain a continuous unitary representation 
$(\pi, \cH_Q)$ of $G$ with 
$\pi(\exp(x)) = e^{\oline{\rho(x)}}$ for every $x \in \g$. 
Then 
\[ \tilde\phi(g) := Q_\1 \pi(g) Q_\1^* \] 
is a positive definite $B(\cK)$-valued function on $G$, 
and for $x \in W_\g$ we obtain from \eqref{eq:12} 
\begin{align*}
\tilde\phi(\exp x) v
&= Q_\1 \pi(\exp x) Q_\1^*v 
= Q_\1 \sum_{n = 0}^\infty \frac{\rho(x)^n}{n!} Q_\1^*v 
= Q_\1 Q_{-\exp x}^*v 
= \phi(\exp x)v.
\end{align*}
Since the kernel 
\[ (g,h) \mapsto \tilde\phi(gh^{-1}) 
= Q_\1 \pi(gh^{-1})Q_\1^* 
= (Q_\1 \pi(g)) (Q_\1\pi(h))^* \] 
is analytic in a neighborhood of the diagonal, 
Theorem~\ref{thm:kerana} implies that the 
map 
\[ G \to B(\cH_Q, \cK), \quad  g \mapsto Q_\1\pi(g)\] 
is analytic, and this in turn implies that 
$\tilde\phi$ is analytic. 
As $\tilde\phi$ coincides with $\phi$ in a $\1$-neighborhood, 
the analyticity of $\phi$ and $\tilde\phi$, together with 
the connectedness of $V$ lead to $\tilde\phi\res_V = \phi$.
\end{prf}

\section{Applications to holomorphically induced representations} 
\mlabel{app:c} 

Let $\g$ be a Banach--Lie algebra and 
$d \in \g$ be an elliptic element, i.e., 
$e^{\R \ad d}$ is bounded. We say that $d$ satisfies 
the {\it splitting condition} if $0$ is an isolated spectral value 
of $\ad d$ (cf.\ \cite{Ne10d}). With $\fh := \ker \ad d$, 
we then obtain an $\ad d$-invariant direct sum decomposition 
$\g_\C = \fp^+ \oplus \h_\C \oplus \fp^-$, 
where the spectrum of $\mp i\ad d$ on $\fp^\pm$ is 
contained in $]0,\infty[$. 
Let $G$ be a connected Lie group with Lie algebra $\g$. 
We consider a bounded representation 
$(\rho,V)$ of $H =\la \exp \fh\ra \subeq G$ and want to 
obtain criteria for the holomorphic inducibility of 
$(\rho, V)$. 

We recall the closed subalgebra 
$\fq = \fp^+ \rtimes \fh_\C \subeq \g_\C$. 
Let $U_\pm \subeq \fp^\pm$ and $U_0 \subeq \fh_\C$ 
be open convex $0$-neighborhoods for which the BCH-multiplication map
\[ U_+ \times U_0 \times U_- \to \g_\C, \quad 
(x_+,x_0,x_-) \mapsto x_+ * x_0 * x_- \] 
is biholomorphic onto an open subset $U$ of $\g_\C$. 
We then define a holomorphic map 
\[ F \: U \to B(V), \quad 
F(x_+*x_0*x_-) := e^{\dd\rho(x_0)}.\]

\begin{thm} \mlabel{thm:c.2} The following are equivalent: 
  \begin{description}
  \item[\rm(i)] $(\rho, V)$ is holomorphically inducible. 
  \item[\rm(ii)] $f_\rho(\exp x) 
:= F(x)$ defines a positive definite analytic 
function on 
a $\1$-neighborhood of $G$. 
  \end{description}
\end{thm}

\begin{prf} (i) $\Rarrow$ (ii): Let $(\pi, \cH)$ be the unitary 
representation of $G$ obtained by holomorphic induction from 
$(\rho, V)$. We identify $V$ with the corresponding closed subspace 
of $\cH$ and write $p_V\: \cH \to V$ for the corresponding orthogonal 
projection. For $v \in V \subeq (\cH^\omega)^{\fp^-}$ 
(\cite[Rem.~2.18]{Ne10d}) we let 
$f_\rho^v \: U_v \to \cH$ 
be a holomorphic map on an open convex $0$-neighborhood 
$U_v \subeq U$ satisfying $f_\rho^v(x) = \pi(\exp x)v$ for $x \in U_v \cap \g$. 
Then $\dd\pi(\fp^-)v = \{0\}$ implies that 
$L_z f_\rho^v = 0$ for $z \in \fp^-$. For $w \in V$ and $z \in \fp^+$, 
we also obtain 
\[ \la (R_z f_\rho^v)(x), w \ra 
= \la \dd\pi(z) f_\rho^v(x), w \ra
= \la f_\rho^v(x), \dd\pi(z^*) w \ra=0.\] 
This proves that $R_z (p_V \circ f_v) = 0$. 
We conclude that, for $x_\pm$ and $x_0$ sufficiently close to $0$, we  have 
\[ p_V f_\rho^v(x_+ * x_0 * x_-) = f\rho^v(x_0) = e^{\dd\rho(x_0)}v = F(x_+ * x_0 * x_-)v.\]
Therefore $p_V \circ f_\rho^v$ extends holomorphically to $U$ and 
\[ \la \pi(\exp x)v,w \ra  = \la F(x)v,w \ra \quad \mbox{ for } \quad 
x \in U_v \cap \g, v,w \in V.\] 
We conclude that $F(x) = p_V \pi(\exp x) p_V$ 
holds for $x$ sufficiently close to $0$, and hence that 
$f_\rho(\exp x) = p_V \pi(\exp x) p_V$ defines a positive definite 
function on a $\1$-neighborhood of the real Lie group~$G$. 

(ii) $\Rarrow$ (i): From Theorem~\ref{thm:extension} it 
follows that some restriction of $f_\rho$ to a possibly smaller 
$\1$-neighborhood in $G$ extends to a global analytic positive 
definite function $\phi$. Then the vector-valued GNS construction 
yields a unitary representation of $G$ on the corresponding reproducing 
kernel Hilbert space $\cH_\phi \subeq V^G$ 
for which all the elements of $\cH_\phi^0 = \Spann(\phi(G)V)$ 
are analytic vectors. In particular, 
$V \subeq \cH_\phi^\omega$ consists of smooth vectors, 
and the definition of $f_\rho$ implies that 
$\dd\pi(\fp^-)V = \{0\}$. Therefore \cite[Thm.~2.17]{Ne10d} 
implies that the representation 
$(\pi, \cH_\phi)$ is holomorphically induced from 
$(\rho,V)$. 
\end{prf}

\section{Classical groups of operators} \mlabel{app:d} 

In this appendix we review the zoo of classical groups 
of operators on Hilbert spaces showing up in this paper.

\subsection{Unitary and general linear groups}\mlabel{app:d.1}

For a Hilbert space $\cH$ over $\K \in \{\R,\C,\H\}$, we 
write $\GL(\cH) = \GL_\K(\cH)$ for the group of $\K$-linear topological 
isomorphisms of $\cH$, which is the unit group 
of the real Banach algebra $B(\cH)$ of bounded $\K$-linear 
operators on $\cH$. It contains the subgroup 
\begin{equation}
  \label{eq:ukh}
\U_\K(\cH) := \U(\cH) := \{ g \in \GL_\K(\cH) \: g^* = g^{-1}\} 
\end{equation}
of unitary operators with Lie algebra 
$\fu(\cH)= \{ x \in \gl(\cH)\: x^* = -x\}$. 
If $\cH$ is real, then we also write 
$\OO(\cH) := \U_\R(\cH)$, and if $\cH$ is quaternionic, we 
also write $\Sp(\cH) := \U_\H(\cH)$ for the corresponding 
unitary groups. 

In many situations it is convenient to describe 
real Hilbert spaces as pairs 
$(\cH,\sigma)$, where $\cH$ is a complex Hilbert space 
and $\sigma \: \cH\to \cH$ is a {\it conjugation}, i.e., 
an antilinear isometry with $\sigma^2 = \id_\cH$. 
Then $\cH^\sigma := \{ v \in \cH \: \sigma(v) = v\}$ 
is a real Hilbert space and, conversely, 
every real Hilbert space $\cH$ arises this way 
by the canonical conjugation $\sigma$ on 
$\cH_\C$ with $(\cH_\C)^\sigma =\cH$. Then 
\begin{equation}
  \label{eq:ohglh}
\OO(\cH) \cong \{ g \in \U(\cH_\C) \: g\sigma = \sigma g\}
\quad \mbox{ and }  \quad 
\GL_\R(\cH) \cong \{ g \in \GL(\cH_\C) \: g\sigma = \sigma g\}.
\end{equation}
For $x \in B(\cH_\C)$ we put $x^\top := \sigma x^* \sigma$ and obtain 
\[ \OO(\cH) =  \{ g \in \U(\cH_\C) \: g^\top = g^{-1}\}, \quad 
\fo(\cH) =  \{ x \in \fu(\cH_\C) \: x^\top = -x\}.\] 
For the complex symmetric bilinear form 
$\beta(x,y) := \la x,\sigma y\ra$ on $\cH_\C$,  the orthogonal group 
is 
\[ \OO(\cH_\C,\beta) = \{ g \in \GL(\cH_\C) \: g^\top = g^{-1}\}\] 
with Lie algebra $\fo(\cH_\C,\beta) \cong \fo(\cH)_\C$. 

In the following we write 
$\Herm(\cH) := \{ x \in B(\cH) \: x^* = x\}$, 
\[ \Sym(\cH) := \{ x \in B(\cH) \: x^\top = x\} 
\quad\mbox{ and } \quad 
\Skew(\cH) := \{ x \in B(\cH) \: x^\top = -x\}.\] 

A quaternionic Hilbert space $\cH$ can be considered as a 
complex Hilbert space $\cH^\C$ (the underlying complex Hilbert 
space), endowed with an {\it anticonjugation} $\sigma$, i.e., 
$\sigma$ is an antilinear isometry with $\sigma^2 = -\1$. 
Then 
\begin{equation}
  \label{eq:sphglh}
\Sp(\cH) =\{ g \in \U(\cH^\C) \: g\sigma = \sigma g\}
\quad \mbox{ and }  \quad 
\GL_\H(\cH) = \{ g \in \GL(\cH^\C) \: g\sigma = \sigma g\}.
\end{equation}
For $x \in B(\cH^\C)$ we put $x^\sharp := \sigma x^* \sigma$ and obtain 
\[ \Sp(\cH) =  \{ g \in \U(\cH^\C) \: g^\sharp = g^{-1}\}, \quad 
\sp(\cH) =  \{ x \in \fu(\cH^\C) \: x^\sharp = -x\}.\] 
For the complex skew-symmetric bilinear form 
$\omega(x,y) := \la x,\sigma y\ra$ on $\cH^\C$,  the symplectic group 
is 
\[ \Sp(\cH^\C,\omega) = \{ g \in \GL(\cH^\C) \: g^\sharp = g^{-1}\}\] 
with Lie algebra $\sp(\cH^\C,\omega) = \sp(\cH)_\C$. 
In \cite{Ne02} we also use the notation 
$\GL(\cH_\C,\sigma) = \OO(\cH_\C, \beta)$ in the real case 
and $\GL(\cH^\C,\sigma) = \Sp(\cH^\C, \omega)$ in the quaternionic 
case. 

\subsection{Root systems of classical groups} \mlabel{app:d.1a}

\begin{defn} \mlabel{def:basic} 
(a) Let $\g$ be a real Lie algebra and $\g_\C$ be its complexification. 
If $\sigma \: \g_\C \to \g_\C$ denotes the complex conjugation with 
respect to $\g$, we write $x^* := -\sigma(x)$ for $x \in \g_\C$, so that 
$\g = \{ x \in \g_\C \: x^* = -x\}$. 
Let $\ft \subeq \g$ be a maximal abelian subalgebra 
and $\fh := \ft_\C\subeq \g_\C$ be its complexification. 
For a linear functional $\alpha \in \fh^*$,  
\[ \g_\C^\alpha = \{ x \in \g_\C \: (\forall h \in \fh)\ 
[h,x]= \alpha(h)x\}\] 
is called the corresponding {\it root space}, and 
\[ \Delta :=  \{ \alpha \in \fh^* \setminus \{0\} \: \g_\C^\alpha
\not= \{0\}\}\] 
is the {\it root system} of the pair $(\g_\C,\fh)$. 
We then have $\g_\C^0 = \fh$ and $[\g_\C^\alpha, \g_\C^\beta] 
\subeq \g_\C^{\alpha+\beta}$, hence in particular 
$[\g_\C^\alpha, \g_\C^{-\alpha}] \subeq \h$. 

(b) If $\g$ is a Banach--Lie algebra, then we 
say that $\ft$ is {\it elliptic} if the group 
$e^{\ad \ft}$ is bounded in $B(\g)$. 
We then have 
\begin{itemize}
\item[\rm(I1)] $\alpha(\ft) \subeq i \R$ for $\alpha \in \Delta$, 
and therefore 
\item[\rm(I2)] $\sigma(\g_\C^\alpha) = \g_\C^{-\alpha}$ 
for $\alpha \in \Delta.$ 
\end{itemize}
\end{defn}

\begin{lem}
  \mlabel{lem:e.1} Suppose that $\ft \subeq\g$ is elliptic. 
For $0 \not= x_{\alpha} \in \g_\C^{\alpha}$, the subalgebra 
\[ \g_\C(x_\alpha) := \Spann_\C\{x_\alpha, x_{\alpha}^*,
[x_\alpha, x_{\alpha}^*]\} \] 
is $\sigma$-invariant and of one of the following types: 
\begin{description}
\item[\rm(A)] The abelian type: $[x_\alpha, x_{\alpha}^*] = 0$, i.e., 
$\g_\C(x_\alpha)$ is two dimensional abelian. 
\item[\rm(N)] The nilpotent type: $[x_\alpha, x_{\alpha}^*] \not= 0$ 
and $\alpha([x_\alpha, x_{\alpha}^*]) = 0$, i.e., 
$\g_\C(x_\alpha)$ is a three dimensional Heisenberg algebra. 
\item[\rm(S)] The simple type: $\alpha([x_\alpha, x_{\alpha}^*]) \not= 0$,
i.e., $\g_\C(x_\alpha) \cong \fsl_2(\C)$. In this case we distinguish 
two cases: 
\begin{description}
\item[\rm(CS)] $\alpha([x_\alpha, x_{\alpha}^*]) > 0$, i.e., 
$\g_\C(x_\alpha) \cap \g \cong \su_2(\C)$, and 
\item[\rm(NS)] $\alpha([x_\alpha, x_{\alpha}^*]) < 0$, i.e., 
$\g_\C(x_\alpha) \cap \g \cong \su_{1,1}(\C) \cong \fsl_2(\R)$. 
\end{description}
\end{description}
\end{lem} 

\begin{prf} (cf.\ \cite[App.~C]{Ne10c})
First we note that, in view of $x_\alpha^* \in \g_\C^{-\alpha}$,  
\cite[Lemma~I.2]{Ne98} applies, and we see that $\g_\C(x_\alpha)$ is of one of
the three types (A), (N) and (S). We note that 
$\alpha([x_\alpha, x_\alpha^*]) \in \R$ 
because of (I2) and $[x_\alpha, x_\alpha^*] \in i\ft$. 
Now it is easy to check that 
$\g_\C(x_\alpha) \cap \g$ is of type (CS), resp., (NS), according to the sign
of this number. 
\end{prf} 

\begin{defn} \label{def:coroot} 
Assume that $\g_\C^\alpha = \C x_\alpha$ is one dimensional 
and that $\g_\C(x_\alpha)$ is of type (S). 
Then there exists a unique element 
$\check \alpha \in \h \cap [\g_\C^\alpha, \g_\C^{-\alpha}]$ with 
$\alpha(\check \alpha) = 2$. It is called the {\it coroot of $\alpha$}. 
The root $\alpha \in \Delta$ is said to
be {\it compact} if for $0\not= x_\alpha \in \g_\C^\alpha$ we have 
$\alpha([x_\alpha,x_\alpha^*]) > 0$ and {\it non-compact} otherwise. 
We write $\Delta_c$ for the set of compact and $\Delta_{nc}$ for the set
of non-compact roots. Lemma~\ref{lem:e.1} implies that 
\begin{equation}
  \label{eq:corootsign}
\check \alpha \in \R^+ [x_\alpha,x_\alpha^*] \quad \mbox{ for } \quad 
\alpha \in \Delta_c
\quad \mbox{ and } \quad 
\check \alpha \in \R^+ [x_\alpha^*,x_\alpha] \quad \mbox{ for } \quad 
\alpha \in \Delta_{nc}.
\end{equation}

The Weyl group $\cW \subeq \GL(\fh)$ is the subgroup generated by 
all reflections 
\[ r_\alpha(x) := x - \alpha(x) \check \alpha. \] 
It acts on the dual space by the adjoint maps 
\[ r_\alpha^*(\beta) := \beta - \beta(\check \alpha) \alpha.\] 
\end{defn}

We now describe the relevant root data for the three types of 
unitary Lie algebras over $\K \in \{\R,\C,\H\}$. 

\begin{ex} \mlabel{ex:d.1a} (Root data of unitary Lie algebras)  
Let $\cH$ be a complex Hilbert space with 
orthonormal basis $(e_j)_{j \in J}$ 
and $\ft \subeq \g := \fu(\cH)$ be the 
subalgebra of all diagonal operators with respect to the $e_j$. 
Then $\ft$ is elliptic and maximal abelian, 
$\fh = \ft_\C \cong \ell^\infty(J,\C)$, and the set of 
roots of $\g_\C= \gl(\cH)$ with respect to $\fh$  is given by 
\[ \Delta = \{ \eps_j - \eps_k \: j\not= k \in J \}.\]
Here the operator 
$E_{jk} e_m := \delta_{km} e_j$ is an $\fh$-eigenvector 
in $\gl(\cH)$ generating the corresponding eigenspace 
and $\eps_j(\diag(h_k)_{k \in J}) = h_j$. 
From $E_{jk}^* = E_{kj}$ it follows that 
\[  (\eps_j - \eps_k)\,\check{}= E_{jj} - E_{kk} 
= [E_{jk}, E_{kj}] = [E_{jk}, E_{jk}^*],\] 
so that $\Delta = \Delta_c$, i.e., all roots are compact.

The Weyl group $\cW$ is isomorphic to the group $S_{(J)}$ of finite 
permutations of $J$, acting in the canonical  way on $\fh$. 
It is generated by the reflections $r_{jk} := r_{\eps_j - \eps_k}$ 
corresponding to the transpositions of $j \not= k \in J$. 
\end{ex}

\begin{ex} \mlabel{ex:d.1b} (Root data of symplectic Lie algebras)  
For a complex Hilbert space $\cH$ with a conjugation $\sigma$, 
we consider the quaternionic Hilbert space $\cH_\H := \cH^2$, where the 
quaternionic structure is defined by the anticonjugation 
$\tilde\sigma(v,w) := (\sigma w, -\sigma v)$. 
Then $\g := \sp(\cH_\H) = \{ x \in \fu(\cH^2) \: \tilde\sigma x = x \tilde\sigma \}$ 
and 
\[ \sp(\cH_\H)_\C = \left\{ \pmat{ A & B \cr C & -A^\top \cr} 
\in B({\cal H}^2) \: B^\top = B,  C^\top  = C\right\}. \] 
Let $(e_j)_{j \in J}$ be an orthonormal basis of $\cH$ and 
$\ft \subeq \g \subeq \fu(\cH^2)$ be the 
subalgebra of all diagonal operators with respect to the basis elements  
$(e_j,0)$ and $(0,e_k)$ of $\cH^2$. Then $\ft$ 
is elliptic and maximal abelian in $\g$. 
Moreover, $\fh = \ft_\C \cong \ell^\infty(J,\C)$, 
consists of diagonal operators of the form 
$h = \diag((h_j), (-h_j))$, and the set of 
roots of $\g_\C$ with respect to $\fh$  is given by 
\[ \Delta = \{ \pm 2 \eps_j, \pm (\eps_j \pm \eps_k) \: j \not= k, j,k
\in J \},\] 
where $\eps_j(h)  = h_j$. If we write 
$E_j \in \fh$ for the element defined by 
$\eps_k(E_j) = \delta_{jk}$, then the coroots are given by 
\begin{equation}
  \label{eq:CJcoroot}
 (\eps_j \pm  \eps_k)\,\check{}= E_j \pm  E_k 
\quad \mbox{ for } \quad j \not=k 
\quad \mbox{ and } \quad 
(2\eps_j)\,\check{}= E_j.
\end{equation}
Here the roots $\eps_j - \eps_k$ correspond to block diagonal operators, 
 the roots $\eps_j + \eps_k$ to strictly upper triangular operators, 
and the roots $-\eps_j - \eps_k$ to strictly lower triangular operators. 
Again, all roots are compact, and the 
Weyl group $\cW$ is isomorphic to the group $N \rtimes S_{(J)}$,
where $N \cong \{\pm 1\}^{(J)}$ is the group of finite 
sign changes on $\ell^\infty(J,\R)$. In fact, 
the reflection $r_{\eps_j - \eps_k}$ acts as a transposition and 
the reflection $r_{2\eps_j}$ changes the sign of the $j$th component. 
\end{ex}

\begin{ex} \mlabel{ex:d.1c} (Root data of  orthogonal Lie algebras)  
Let $\cH_\R$ be an  infinite dimensional real Hilbert space 
and $(e_j)_{j \in \tilde J}$ be an orthonormal basis of $\cH_\R$. 
Since $\tilde J$ is infinite, 
it contains a subset $J$ for which there exists an involution 
$\eta \: \tilde J \to \tilde J$ with 
$\eta(J) = \tilde J \setminus J$ and $\tilde J = J \dot\cup \eta(J)$. 
Then 
\[ I e_j := 
\begin{cases}
e_{\eta(j)} & \text{ for } j \in J \\ 
-e_{\eta(j)} & \text{ for } j \in \eta(J) 
\end{cases}\]
defines an orthogonal complex structure on $\cH_\R$. 
This complex structure defines on $\cH_\R$ the structure of a complex 
Hilbert space $\cH := (\cH_\R,I)$. 
We write $\sigma$ for the conjugation on $\cH$ defined by 
$\sigma(e_j) = e_j$ for $j \in J$. 

Then $\iota \: \cH \to \cH \oplus \cH, 
v \mapsto \frac{1}{\sqrt 2}(v,  \sigma(v))$ 
is real linear and isometric. Since its image is a totally 
real subspace, $\iota$ extends to a unitary 
isomorphism $\iota_\C \: \cH_\C \to \cH^2$. 
Let $\beta \: \cH^2 \to \C$ denote the complex bilinear 
extension of the scalar product of $\cH_\R$, so that 
$\fo(\cH_\R)_\C \cong \fo(\cH^2,\beta).$ 
It is given by 
\[ \beta((x,y), (x',y')) 
= \beta(x,y') + \beta(x',y) 
= \la x, \sigma(y')\ra + \la x', \sigma(y)\ra \] 
because the right hand side is complex bilinear and 
has the correct restriction to $\iota(\cH)$. 
This implies that 
\begin{align*}
\fo(\cH_\R)_\C  &\cong\fo(\cH^2, \beta) 
= \Big\{ X \in \gl(\cH^2) \: X^\top \pmat{ 0 & \1 \cr \1 & 0\cr} 
 + \pmat{ 0 & \1 \cr \1 & 0\cr}X = 0\Big\} \\ 
&= \left\{ \pmat{ A & B \cr C & -A^\top \cr} 
\in B({\cal H}^2) \: B^\top = -B,  C^\top  = -C\right\}. 
\end{align*}
For the conjugation 
$\tilde\sigma(v,w) = (\sigma(w), \sigma(v))$ with respect to $\iota(\cH)$, 
we have 
\[ \g := \fo(\cH_\R) \cong  \{ x \in \fu(\cH^2) \: \tilde\sigma x 
= x \tilde\sigma \}. \] 
Now $(e_j)_{j \in J}$ is an orthonormal basis of the complex 
Hilbert space $\cH$ and the subalgebra 
$\ft \subeq \g \subeq \fu(\cH^2)$ 
of all diagonal operators with respect to the basis elements  
$(e_j,0)$ and $(0,e_k)$ of $\cH^2$ is elliptic and maximal 
abelian in~$\g$. Again, $\fh = \ft_\C \cong \ell^\infty(J,\C)$  
consists of diagonal operators of the form 
$h = \diag((h_j), (-h_j))$, and the set of 
roots of $\g_\C$ with respect to $\fh$  is given by 
\[ \Delta = \{ \pm (\eps_j \pm \eps_k) \: j \not= k, j,k\in J \}.\] 
As in Example~\ref{ex:d.1b}, all roots are compact, 
but since $2\eps_j$ is not a root and 
the reflection $r_{\eps_j + \eps_k}$ changes the sign of the $j$- and the 
$k$-component, the 
Weyl group $\cW$ is isomorphic to the group $N_{\rm even} \rtimes S_{(J)}$,
where $N_{\rm even}$ is the group of finite even sign changes. 
\end{ex}

\subsection{$c$-duality and complexification of $\gl_\K(\cH)$} 
\mlabel{app:d.1b}

In the preceding discussion we have seen $3$ types of 
unitary groups $\U_\K(\cH)$: 
$\OO(\cH)$, $\U(\cH)$ and $\Sp(\cH)$ for $\K = \R$,$\C$ and $\H$, 
and their complexifications $\OO(\cH_\C,\beta)$, 
$\GL(\cH)$ and $\Sp(\cH^\C,\omega)$. 
We have also seen the subgroups 
$\GL_\K(\cH)= \U_\K(\cH)\exp(\Herm_\K(\cH))$ of $\U_\K(\cH)_\C$ 
which are symmetric Lie groups with respect to the involution 
$\theta(g) = (g^*)^{-1}$ and the corresponding decomposition 
$\gl_\K(\cH) = \fu_\K(\cH) \oplus \Herm_\K(\cH)$ of the Lie algebra. 
The corresponding $c$-dual symmetric Lie algebras correspond 
to unitary groups: 
\begin{description}
\item[$(\R)$] $\gl_\R(\cH)^c = \fo(\cH) \oplus i \Sym(\cH) \cong \fu(\cH_\C)$. 
\item[$(\C)$] $\gl_\C(\cH)^c = \fu(\cH) \oplus i \Herm(\cH) \cong\fu(\cH)^2$. 
\item[$(\H)$] $\gl_\H(\cH)^c = \sp(\cH) \oplus i \Herm_\H(\cH) 
\cong \fu(\cH^\C)$. 
\end{description}
The complexifications of these Lie algebras are 
\[ \gl_\R(\cH)_\C \cong \gl(\cH_\C), \quad
\gl_\C(\cH)_\C \cong \gl(\cH)^2 \quad\mbox{ and } \quad 
\gl_\H(\cH)_\C \cong \gl(\cH^\C).\]

\subsection{Variants of hermitian groups} \mlabel{app:d.1c}

We shall also need the following variant of the hermitian 
groups. 

\begin{ex}\mlabel{ex:d1}
If $\cH = \cH_+ \oplus \cH_-$ is an orthogonal decomposition 
of the complex Hilbert space $\cH$, then 
$h((z_+, z_-), (w_+, w_-)) := \la z_+, w_+\ra - \la z_-, w_-\ra$ 
defines a hermitian form on $\cH$. 
We write 
\begin{equation}
  \label{eq:uhpm}
\U(\cH_+, \cH_-) := \U(\cH,h)\subeq\GL(\cH)
\end{equation}
for the corresponding group of complex linear $h$-isometries. 
Its Lie algebra $\fu(\cH_+, \cH_-)$ is a real form 
of $\gl(\cH)$. 

Let $J := J_+ \dot\cup J_-$ be such that 
$(e_j)_{j \in J_\pm}$ is an ONB of $\cH_\pm$. 
We have seen in Example~\ref{ex:d.1a}, that the 
root system $\Delta$ of $\fu(\cH_+, \cH_-)_\C\cong \gl(\cH) \cong\fu(\cH)_\C$ with respect to 
the subalgebra $\fh\cong \ell^\infty(J,\C)$ 
of diagonal matrices with respect to the $e_j$ is given by 
\[ \Delta = \{ \eps_j - \eps_k \: j\not= k \in J \}\]
(cf.\ Example~\ref{ex:d.1a}). 
For $d := \frac{i}{2}\diag(\1,-\1)$, 
we obtain the compact roots 
\[ \Delta_c = \Delta_k  := \{ \alpha \in \Delta \: \alpha(d) = 0\} 
= \{ \eps_j - \eps_k \: j\not= k \in J_+; j\not= k \in J_-\} \] 
corresponding to the complexification 
$\gl(\cH_+) \oplus \gl(\cH_-)$ of the centralizer 
\[ \fu(\cH_+, \cH_-) \cap \fu(\cH) 
\cong \fu(\cH_+) \oplus \fu(\cH_-)\] 
of $d$ in $\fu(\cH_+, \cH_-)$, 
and the non-compact roots $\Delta_{nc} =\Delta_p^+ \cup \Delta_p^-$, 
where 
\begin{equation}
  \label{eq:deltapplus}
\Delta_p^\pm  := \{ \alpha \in \Delta \: \alpha(-i\cdot d) = \pm 1\} 
= \pm \{ \eps_j - \eps_k \: j \in J_+, k \in J_-\} 
\end{equation}
correspond to the $\pm i$-eigenspaces 
of $\ad d$ in $\gl(\cH)$ (cf.\ Example~\ref{ex:5.12b}(b)). 
\end{ex}

\begin{ex}\mlabel{ex:d2} 
For a complex Hilbert space $\cH$ with a conjugation $\sigma$, 
we define the corresponding {\it symplectic group} by 
\begin{equation}
  \label{eq:sph}
\Sp(\cH)
:= \left\{ g \in \U(\cH,\cH) \: g^\top
\pmat{ 0 & \1 \cr -\1 & 0\cr} g = \pmat{ 0 & \1 \cr -\1 & 0\cr}
\right\}
\end{equation}
with the Lie algebra 
\[ \sp(\cH) = \left\{ \pmat{ A & B \cr B^* & -A^\top \cr} 
\in B({\cal H}^2) \: A^* = - A, B^\top = B\right\} \] 
and note that 
\[ \sp(\cH)_\C = \left\{ \pmat{ A & B \cr C & -A^\top \cr} 
\in B({\cal H}^2) \: B^\top = B,  C^\top  = C\right\}. \] 
The anticonjugation $\tilde\sigma(v,w) := (\sigma w, -\sigma v)$ defines a 
quaternionic structure on $\cH^2$. We write 
$\cH_\H$ for the so obtained quaternionic Hilbert space. 
The skew-symmetric complex bilinear form on $\cH^2$ defined by 
$\tilde\sigma$ is given by 
\[\omega((v,w),(v',w')) = \la (v,w), (-\sigma w', \sigma v')\ra 
= \la w,\sigma v'\ra - \la v,\sigma w'\ra\]  
and $\sp(\cH)_\C 
\cong \sp(\cH^2,\omega) \cong \sp(\cH_\H)_\C$ 
(cf.\ Example~\ref{ex:d.1b}). 

Let $(e_j)_{j \in J}$ be an ONB of $\cH$. 
In Example~\ref{ex:d.1b} we have seen that 
the root system $\Delta$ of $\sp(\cH^2,\omega)$ with respect to 
the subalgebra $\fh\cong \ell^\infty(J,\C)$ 
of diagonal matrices is 
\[ \Delta = \{ \pm 2 \eps_j, \pm (\eps_j \pm \eps_k) \: j \not= k, j,k
\in J \},\] 
For $d := \frac{i}{2} \diag(\1,-\1)\in \fh$, we obtain 
the compact roots 
\[ \Delta_c = \Delta_k  := \{ \alpha \in \Delta \: \alpha(d) = 0\} 
= \{ \eps_j - \eps_k \: j\not= k \in J\} \] 
corresponding to the complexification of 
$\sp(\cH) \cap \fu(\cH^2) \cong \fu(\cH),$ 
and the non-compact roots 
\[ \Delta_p^\pm  := \{ \alpha \in \Delta \: \alpha(-i\cdot d) = \pm 1\} 
= \pm\{ \eps_j + \eps_k \: j, k \in J\}.\]  
\end{ex}

\begin{ex}
  \mlabel{ex:d3}
In a similar fashion, we define for a complex Hilbert space 
$\cH$ with a conjugation $\sigma$ 
\begin{equation}
  \label{eq:ostarh}
\OO^*(\cH)
:= \left\{ g \in \U(\cH,\cH) \: g^\top
\pmat{ 0 & \1 \cr \1 & 0\cr} g = \pmat{ 0 & \1 \cr \1 & 0\cr}
\right\}
\end{equation}
with the Lie algebra 
\[ \fo^*(\cH) = \left\{ \pmat{ A & B \cr B^* & -A^\top \cr} 
\in B({\cal H}^2) \: A^* = - A, B^\top = -B\right\} \] 
and 
\[ \fo^*(\cH)_\C = \left\{ \pmat{ A & B \cr C & -A^\top \cr} 
\in B({\cal H}^2) \: B^\top = -B, C^\top = -C\right\}. \] 
The conjugation $(v,w) \mapsto (\sigma w, \sigma v)$ defines a 
symmetric complex bilinear form 
\[\beta((v,w),(v',w')) = \la (v,w), (\sigma w', \sigma v')\ra 
= \la v,\sigma w'\ra + \la w,\sigma v'\ra \]  
on $\cH^2$ with $\fo^*(\cH)_\C \cong \fo(\cH^2,\beta)$. 

Let $(e_j)_{j \in J}$ be an ONB of $\cH$. 
According to 
Example~\ref{ex:d.1c}, the 
root system $\Delta$ of $\fo(\cH^2,\beta)$ with respect to 
the subalgebra $\fh\cong \ell^\infty(J,\C)$ 
of diagonal matrices is given by 
\[ \Delta = \{ \pm \eps_j \pm \eps_k \: j \not= k \in J \},\] 
For $d := \frac{i}{2}\diag(\1,-\1) \in \fh$, we obtain 
the compact roots 
\[\Delta_c= \Delta_k  := \{ \alpha \in \Delta \: \alpha(d) = 0\} 
= \{ \eps_j - \eps_k \: j\not= k \in J\} \] 
corresponding to the complexification of 
$\fo^*(\cH) \cap \fu(\cH^2) \cong \fu(\cH)$ and the non-compact roots  
\[ \Delta_p^\pm  := \{ \alpha \in \Delta \: \alpha(-i\cdot d) = \pm 1\} 
= \pm\{ \eps_j + \eps_k \: j\not= k \in J\}. \]  
\end{ex}

\begin{ex}
  \mlabel{ex:d4} 
For a real Hilbert space $\cH_\R$, we consider the pseudounitary 
group $\OO(\R^2, \cH_\R)$ of the indefinite quadratic form 
\[ q(x,v) := \|x\|^2 - \|v\|^2 
= \la Q(x,v), (x,v)\ra, \quad Q = \pmat{\1 & 0 \\ 0 & -\1},\]  
on $\R^2 \oplus \cH_\R$. 
Let $\cH := (\R^2 \oplus \cH_\R, I)$ be the complex 
Hilbert space, where $I$ is a complex structure on 
$\R^2 \oplus \cH_\R$ given by $I(x,y) = (-y,x)$ on $\R^2$ and 
an by an orthogonal complex structure on $\cH_\R$. 
Let $(e_j)_{j \in J}$ be an orthonormal basis of this complex 
Hilbert space and $j_0 \in J$ with $e_{j_0} = (1,0) \in \R^2$ 
and define a conjugation 
$\sigma$ on $\cH$ by $\sigma(e_j) = e_j$ for $j \in J$. 

We realize $\OO(\R^2, \cH_\R)$ as a subgroup of $\OO(\cH^2,\beta)$, where 
$\beta$ is the complex bilinear extension 
of $\beta$ to $\cH^2 \cong \cH_\C$ as follows. 
The conjugation $(v,w) \mapsto (\sigma w, \sigma v)$ defines a 
symmetric complex bilinear form 
\[\beta((v,w),(v',w'))= \la (Qv,Qw), (\sigma w', \sigma v')\ra 
= \la Qv,\sigma w'\ra + \la Qw,\sigma v'\ra \]  
on $\cH^2$ with $\fo(\R^2,\cH_\R)_\C \cong \fo(\cH^2,\beta)$. 
Then the corresponding complex orthogonal group is 
\begin{equation}
  \label{eq:o2}
\OO(\cH^2,\beta) 
= \left\{ g \in \GL(\cH^2) \: g^\top
\pmat{ 0 & Q \cr Q & 0\cr} g = \pmat{ 0 & Q \cr Q & 0\cr}
\right\}
\end{equation}
with the Lie algebra 
\begin{align*}
\fo(\cH^2,\beta) 
&= \left\{ X \in B(\cH^2) \: X^\top
\pmat{ 0 & Q \cr Q & 0\cr} + \pmat{ 0 & Q \cr Q & 0\cr} X = 0\right\} \\
&= \left\{ \pmat{ A & B \cr C & D \cr} 
\in B({\cal H}^2) \: B^\top =-QBQ, C^\top =-QCQ, D = -QA^\top Q\right\}.
\end{align*}
The subalgebra $\fh\cong \ell^\infty(J,\C)$ 
of diagonal matrices 
$\diag((h_j), (-h_j)) \in \gl(\cH^2)$ 
is maximal abelian in $\fo(\cH^2,\beta)$ 
and the root system $\Delta$ of $\fo(\cH^2,\beta)$ with respect to 
$\fh$ is given by 
\[ \Delta = \{ \pm \eps_j \pm \eps_k \: j \not= k \in J \}.\] 
Then 
\[ d := \diag\Big(\pmat{0 & -1 \\ 1&0},0\Big) \in \fo(\R^2, \cH_\R)\]  
corresponds to $h  = \diag(h_j)\in \fh$ with 
$h_j = \delta_{j,j_0} i$, so that 
\[ \Delta_k  := \{ \alpha \in \Delta \: \alpha(d) = 0\} 
= \{ \pm \eps_j \pm \eps_k \: j \not= k \in J\setminus \{j_0\} \}\] 
and 
\[ \Delta_p^\pm  := \{ \alpha \in \Delta \: \alpha(-i\cdot d) = \pm 1\} 
= \{ \pm\eps_{j_0} + \eps_j, \pm\eps_{j_0} - \eps_j \: j_0 \not= j\in J\}. \]  
\end{ex}

\subsection{Groups related to Schatten classes} \mlabel{app:d.2}

For the algebra $B(\cH)$ of bounded operators on the $\K$-Hilbert space 
$\cH$, the ideal of compact operators 
is denoted $K(\cH) = B_\infty(\cH)$, and for $1 \leq p < \infty$, we write 
\[ B_p(\cH) := \{ A \in K(\cH) \: \tr((A^*A)^{p/2}) < \infty\}\]  
for the {\it Schatten ideals}. In particular, 
$B_2(\cH)$ is the space of {\it Hilbert--Schmidt operators} 
and $B_1(\cH)$ the space of {\it trace class operators}. 
For two Hilbert spaces $\cH_\pm$, we put 
\[ B_2(\cH_-, \cH_+) := \{ A \in B(\cH_-, \cH_+) \: \tr(A^*A) < \infty\}.\]
For $1 \leq p \leq \infty$, we obtain Lie groups 
\[ \GL_p(\cH) := \GL(\cH) \cap (\1 + B_p(\cH)) 
\quad \mbox{ and } \quad 
\U_p(\cH) := \U(\cH) \cap \GL_p(\cH) \] 
with the Lie algebras 
\[ \gl_p(\cH) := B_p(\cH) \quad \mbox{ and } \quad 
\fu_p(\cH) := \fu(\cH) \cap \gl_p(\cH). \] 
For $\K = \C$, we have a {\it determinant homomorphism} 
\[ \det \: \GL_1(\cH) \to \C^\times \quad \mbox{ with }\quad 
\det(\U_1(\cH)) = \T.\] 
With $\SU(\cH):= \ker(\det\res_{\U_1(\cH)})$, we then obtain 
\begin{equation}
  \label{eq:ustilde}
\U_1(\cH) \cong \T \ltimes \SU(\cH) 
\quad \mbox{ and } \quad 
\tilde\U_1(\cH) \cong \R \ltimes \SU(\cH) 
\end{equation}
for the simply connected covering group (cf.\ \cite[Prop.~IV.21]{Ne04}). 

\subsection{Restricted groups} \mlabel{app:d.3}

Let $\cH$ be a complex 
Hilbert space which is a direct sum $\cH = \cH_+ \oplus \cH_-$. 
Then 
\[  B_{\rm res}(\cH) := \Big\{ \pmat{a & b \\ c &d} \in B(\cH) \: 
b\in B_2(\cH_-, \cH_+), c\in B_2(\cH_+, \cH_-)\Big\} \] 
is a complex Banach-$*$ algebra. Its unit group is 
\[ \GL_{\rm res}(\cH) = \GL(\cH) \cap B_{\rm res}(\cH).\] 
Intersecting with $\GL_{\rm res}(\cH)$, we obtain various 
{\it restricted classical groups}: 
\begin{description}
\item  $\U_{\rm res}(\cH) := \U(\cH) \cap \GL_{\rm res}(\cH)$ 
(the {\it restricted unitary group}).
\item  $\U_{\rm res}(\cH_+,\cH_-) := \U(\cH_+,\cH_-) \cap \GL_{\rm res}(\cH)$
(the {\it restricted pseudo-unitary group}).
\item  $\Sp_{\rm res}(\cH) := \Sp(\cH) \cap \GL_{\rm res}(\cH \oplus \cH)$
(the {\it restricted symplectic group}).
\item  $\OO^*_{\rm res}(\cH) := \OO^*(\cH) \cap \GL_{\rm res}(\cH \oplus \cH)$.
\end{description}

The unitary group in $c$-duality with 
$\U_{\rm res}(\cH_+,\cH_-)$ is $\U_{\rm res}(\cH_+ \oplus \cH_-)$ 
(cf.\ Example~\ref{ex:1.7}). 
The groups 
$\Sp_{\rm res}(\cH)$ and $\OO^*_{\rm res}(\cH)$ 
also have corresponding $c$-dual unitary groups which 
can be realized as follows. 
Let $\cH_\H$ denote $\cH^2$, endowed with its canonical 
quaternionic structure given by $\tilde\sigma(v,w) = (-\sigma w,\sigma v)$, 
so that 
$\gl_\H(\cH_\H)_\C \cong \gl(\cH^2)$ (cf.\ Appendix~\ref{app:d.1b}) 
and $\sp(\cH_\H)_\C \cong \sp(\cH^2,\omega)$. Then 
\[ \Sp_{\rm res}(\cH_\H) := \Sp(\cH_\H) \cap \GL_{\rm res}(\cH \oplus \cH)\] 
is a group whose Lie algebra 
$\sp_{\rm res}(\cH_\H)$ is $c$-dual to $\sp_{\rm res}(\cH)$. 

For a real Hilbert space $\cH_\R$ with complex structure $I$ 
and the corresponding complex Hilbert space 
$\cH = (\cH_\R,I)$, the realization of $\OO(\cH_\R)$ as 
$\U(\cH^2) \cap \OO(\cH^2, \beta)$ leads to the 
{\it restricted unitary group} 
\[ \OO_{\rm res}(\cH_\R) := \OO(\cH_\R) \cap \GL_{\rm res}(\cH \oplus \cH)\] 
of those orthogonal operators $g = g_0 + g_1$ on $\cH_\R$ for which the 
antilinear part $g_1$ with respect to $I$ is Hilbert--Schmidt 
and  the subgroup of $I$-linear elements in $\OO_{\rm res}(\cH_\R)$ 
is the unitary group $\U(\cH)$ of the complex Hilbert space $\cH$. 

\subsection{Doubly restricted groups} \mlabel{app:d.4}

If $\cH = \cH_+ \oplus \cH_-$ is an orthogonal decomposition of the 
complex Hilbert space $\cH$, then we  write 
operators on $\cH$ as $(2 \times 2)$-matrices. 
The subspace 
\begin{equation}
  \label{eq:b1,2}
B_{1,2}(\cH) := \Big\{ x  = \pmat{ x_{11} & x_{12} \\ x_{21} & x_{22}} 
\in B(\cH) \: 
\|x_{11}\|_1, \|x_{12}\|_2, \|x_{21}\|_2, 
\|x_{22}\|_1 < \infty \Big\} 
\end{equation}
is a Banach algebra with the associated group 
\begin{equation}
  \label{eq:gl1,2}
\GL_{1,2}(\cH) := \{ g \in \GL(\cH) \: \1 - g \in B_{1,2}(\cH) \}. 
\end{equation}

For a subgroup $G \subeq \GL(\cH)$, the corresponding 
{\it doubly restricted group} is now defined as 
$G_{1,2} := G \cap \GL_{1,2}(\cH)$. We shall need these groups 
for $G = \U_{\rm res}(\cH_+, \cH_-)$, 
$\OO^*_{\rm res}(\cH)$, $\Sp_{\rm res}(\cH)$, and 
also for the corresponding unitary $c$-dual groups 
$\U_{\rm res}(\cH)$, $\OO_{\rm res}(\cH_\R)$, $\Sp_{\rm res}(\cH_\H)$. 

\subsection{Polar decomposition} 

\begin{rem} \mlabel{rem:polar}
In \cite[Thm.~II.6, Prop.~III.8]{Ne02} it is shown that for all groups 
$G$ discussed in this appendix, the Lie algebra $\g$ decomposes as 
\[ \g = \fk \oplus \fp \quad \mbox{ with } \quad 
\fk = \{ x \in \g \: x^* = -x\}\quad \mbox{ and } \quad 
\fp = \{ x \in \g \: x^* = x\},\] 
and that the corresponding polar map 
$K \times \fp \to G, (k,x) \mapsto k \exp x$ 
is a diffeomorphism (see also the appendix in \cite{Ne04}). 
For the doubly restricted group 
$G_{1,2}$, see in particular \cite[Lemma~III.6]{NO98}. 
\end{rem}

\section{Bounded representations 
of $\U_p(\cH)$, $1 < p \leq \infty$}\label{Sect2} 
\mlabel{app:e} 

In this section we completely describe the bounded 
representations of the unitary groups 
$G := \U_p(\cH)_0$, $1 < p \leq \infty$, where $\cH$ is an 
infinite dimensional Hilbert space over $\K \in \{\R,\C,\H\}$. 
For $\K =\C$, we have shown in 
 \cite[Thm.~III.14]{Ne98} that all bounded 
unitary representations of $G$, resp., 
all holomorphic $*$-representations of the complexified group 
$G_\C = \GL_p(\cH)$ are direct sums of irreducible ones 
and that the irreducible ones are 
classified by their ``highest weights''. 
In this section we explain this classification in some detail 
and extend it to real and quaternionic groups. 

To this end, we work with the root systems 
$\Delta$ from Examples~\ref{ex:d.1a}, \ref{ex:d.1b} and \ref{ex:d.1c}, 
where $\fh = \ft_\C \cong \ell^p(J,\C) \subeq \g_\C$ now stands for the 
Lie algebra of diagonal operators in $\g_\C$ and 
$\ft = \g \cap \fh$. The root data does not depend on the parameter~$p$.

\begin{defn}\label{def_weights} 
A continuous linear functional $\beta \in \fh'$
is called a 
{\it weight} if $\beta(\check \alpha) \in \Z$ holds for each 
$\alpha \in \Delta$. We write $\cP$ for the additive group of weights 
and $\cQ \subeq \cP$ for the subgroup generated by $\Delta$. 
Note that the Weyl group $\cW$ acts naturally on $\cP$ 
(cf.\ \cite{Ne98}, \cite{Ne04}). 
\end{defn}

\begin{rem} (a) For $q \in [1,\infty[$ defined by 
$\frac{1}{q} + \frac{1}{p} = 1$, we have 
$\fh' \cong \ell^q(J,\C)$, so that we consider elements 
of this space as functions $\beta \: J \to \C$. 
Since the root system $\Delta$ contains 
all roots of the form $\eps_j - \eps_k$, $j\not=k \in J$, 
it follows from \cite[Prop.~III.10]{Ne98} that 
each weight $\beta \in \cP$ is finitely supported with 
$\beta(J) \subeq \Z$. Conversely, 
the description of the coroots in 
\eqref{eq:CJcoroot} implies for $\K = \H$ and $\R$ 
(Examples~\ref{ex:d.1b} and \ref{ex:d.1c})  
that any finitely supported $\Z$-valued function on 
$J$ is a weight, i.e., 
$\cP \cong \Z^{(J)}$ is a free group over $J$. 

(b) The weight group $\cP$ 
can be identified with the character group of the 
Banach--Lie group $T = \exp \ft$ 
by assigning to $\lambda$ the character given by 
$\chi_\lambda(t) := \prod_{j \in J} t_j^{\lambda_j}$. 
\end{rem} 

\begin{defn} \mlabel{def:tensrep} 
From \cite[Prop.~III.10, Sects.~VI, VIII]{Ne98} we know that, 
for each $\lambda\in\cP$, there exists a unique bounded 
unitary representation $(\pi_\lambda, \cH_\lambda)$ 
of $\U_p(\cH)$ whose weight set is given by 
\begin{equation}
  \label{eq:weightset}
\cP_\lambda = \conv(\cW \lambda) \cap (\lambda + \cQ). 
\end{equation}

From \cite[Thm.~III.15, Sects.~VI,VII]{Ne98}, 
we even know that these representation 
extend to bounded representations of the 
full unitary group $\U(\cH)$. The uniqueness of the extension
 to $\U_\infty(\cH)$ follows from the density of 
$\U_p(\cH)$ in $\U_\infty(\cH)$, and the uniqueness of 
the extension to $\U(\cH)$ follows from 
the perfectness of the Lie algebra $\fu(\cH)$ 
(\cite[Lemma~I.3]{Ne02}) because 
any two extensions $\U(\cH) \to \U(\cH_\lambda)$ 
differ by a continuous homomorphism 
\[ \U(\cH) \to \U(\cH_\lambda) \cap \pi_\lambda(\U_p(\cH))' = \T \1\] 
by Schur's Lemma. We also write 
$\pi_\lambda$ for the unique extension to $\U(\cH)$. 
\end{defn} 

\begin{rem} (a) Since the representation $(\pi_\lambda, \cH_\lambda)$ 
of $\U_\infty(\cH)$ is uniquely determined by \eqref{eq:weightset}, 
and the right hand  side only depends on the Weyl group orbit, which 
in turn coincides with the set of extreme points of its 
convex hull (\cite[Lemma~I.19]{Ne98}), it follows that 
$\pi_\mu \cong \pi_\lambda$ if and only if $\mu \in \cW\lambda$. 
Hence the equivalence classes of these representations are 
parameterized by the orbit space $\cP/\cW$. 

(b) For $\K = \C$, we write every $\lambda\in\cP$ as 
$\lambda = \lambda_+ - \lambda_-$ for 
$\lambda_\pm := \max(\pm\lambda, 0)$. 
With \cite[Thm.~2.2]{BN11} we then obtain the factorization 
\[ \pi_\lambda\cong\pi_{\lambda_+}\otimes\pi_{\lambda_-}^*\] 
(see also \cite{Ki73} for the case where $\cH$ is separable). 
\end{rem}

The first part of the following theorem is a variation of 
\cite[Thm.~III.14]{Ne98} which is the corresponding 
result for the groups $\U_p(\cH)$, where $\cH$ is a complex 
Hilbert space. From \cite{Ne98} we know that this result 
does not extend to~$p =1$.

\begin{thm} \mlabel{thm:oprep} 
Let $\cH_\R$ be an infinite dimensional real  
Hilbert space. For $1 < p \leq \infty$, every bounded 
unitary representation of the simply connected 
covering group 
$\tilde\OO_p(\cH_\R)_0$ of $\OO_p(\cH_\R)_0$ 
is a direct sum of bounded 
irreducible representations. The bounded irreducible representations 
are highest weight representations $(\rho_{\lambda}, V)$ with 
finitely supported highest weights $\lambda \: J \to \Z$. 
In particular, they have a unique extension to the 
full orthogonal group $\OO(\cH_\R)$. 
\end{thm}

\begin{prf} Let $G := \tilde\OO_p(\cH_\R)_0$. 
Passing to the derived representation 
of the complexified Lie algebra 
$\g_\C = \fo_p(\cH_\R)_\C \cong \fo_p(\cH^2,\beta)$, we obtain 
a $*$-representation $(\rho,V)$ of $\g_\C$. 
Applying \cite[Lemma~III.13]{Ne98} to $\fh \subeq \g_\C$, it follows that 
each bounded $*$-representation $(\rho,V)$ of $\g_\C$ 
is a direct sum of $\fh$-weight spaces and that the weight set 
$\cP_V$ of $V$ is contained in $\cP \cong \Z^{(J)}$. 
For $\mu \in \Z^{(J)}$,  we have 
\[ \|\mu\|_1 
= \sum_j |\mu_j| \leq  \sum_j |\mu_j|^q = \|\mu\|_q^q,\] 
so that the boundedness of $\rho$ implies the boundedness of 
the set $\cP_V \subeq \fh'$ of weights of $(\rho,V)$ as a subset of 
$\ell^1(J,\C) \subeq \ell^\infty(J,\C)'$.  

Let $\gl_p(\cH) \into \fo_p(\cH^2,\beta), x\mapsto \diag(x,-x^{\top})$ 
denote the canonical 
embedding and $d := \frac{i}{2}\diag(\1,-\1)$. Then 
\cite[Thm.~III.14]{Ne98} implies  that 
the restriction of $\rho$ to $\gl_p(\cH)$ is a 
direct sum of irreducible representations. 
From the decomposition into weight spaces 
for $\fh \cong \ell^p(J,\C)$ and the boundedness of 
the weight set $\cP_V$ in $\ell^1$, we obtain an extension 
of $\rho$ from $\fh$ to a bounded representation 
of $\fh_\infty := \ell^\infty(J,\C)$, hence to a 
bounded representation $(\hat\rho,V)$ of the semidirect 
sum $\g_\C \rtimes \fh_\infty$. 
Since the operator $-i\hat\rho(d)$ is bounded with integral eigenvalues, 
it has only finitely many eigenspaces, and all its eigenspaces 
are $\gl_p(\cH)$-invariant. 

Let $(\rho_0, V_0)$ be an irreducible 
$\gl_p(\cH)$-subrepresentation of the eigenspace 
corresponding to the maximal eigenvalue and 
$\fp^\pm \subeq \g_\C$ denote the $\pm i$-eigenspaces of $\ad d$. 
Then $\hat\rho(\fp^+)V_0 = \{0\}$. 
Next we observe that  the representation 
of $\g$ on the closed subspace 
$\hat V_0$ generated by $V_0$ is also $\hat\rho(d)$-invariant, hence 
a representation of the real Lie algebra $\fg + \R d$. 
We conclude with \cite[Thm.~2.17]{Ne10d} 
that it is holomorphically induced from 
the irreducible representation 
$(\rho_0, V_0)$ of $\gl_p(\cH) + \C d$, hence irreducible 
by \cite[Cor~2.14]{Ne10d}. 
This shows that every 
bounded representation $(\rho,V)$ of $\g$ on a non-zero Hilbert space 
contains an irreducible one, and therefore Zorn's Lemma implies 
that $\rho$ is a direct sum of irreducible representations. 

Now we suppose that $(\rho,V)$ is irreducible 
and that $V_0$ is chosen as above.  
For each weight $\mu \: J \to \Z$ in $\cP_V$ we have 
$\mu(-i\cdot d) = \sum_{j \in J} \mu_j.$ 
Since $\cP_V$ is invariant under the Weyl group $\cW$ 
which contains all finite sign changes (Example~\ref{ex:d.1b}), 
the maximality of $\mu(-i\cdot d)$ among 
$(\cW\mu)(-i\cdot d)$ implies that each weight $\mu \in \cP_{V_0}$ satisfies 
$\mu_j \in \N_0$ for each $j \in J$. 
Since the representation $(\rho_0, V_0)$ of 
$\gl(\cH)$ is a representation with 
some highest weight $\lambda$ (\cite[Thm.~III.14]{Ne98}), it follows 
from \cite[Prop.~VII.2]{Ne98} that 
$(\rho,V)$ contains the  representation 
with the highest weight $\lambda$, hence is a highest weight 
representations because it is irreducible. 
\end{prf}

Almost the same arguments as for the orthogonal 
groups apply to the unitary groups $\Sp_p(\cH_\H)$ of a 
quaternionic Hilbert space: 

\begin{thm} \mlabel{thm:sprep} Let $\cH_\H = \cH^2$ be the 
quaternionic Hilbert space canonically associated to the complex 
Hilbert space $\cH$ and a conjugation $\sigma$ on $\cH$. 
For $1 < p \leq \infty$, every bounded unitary representation of 
$\Sp_p(\cH_\H)$ is a direct sum of bounded 
irreducible representations. The bounded irreducible representations 
are highest weight representations $(\rho_{\lambda}, V)$ with 
finitely supported highest weights $\lambda \: J \to \Z$. 
In particular, they have a unique extension to the 
full symplectic group $\Sp(\cH_\H)$. 
\end{thm}

For a separable Hilbert space, the 
representations discussed in this section 
coincide with those representations of 
$\U_\infty(\cH)$ extending to strongly continuous 
representations on the full unitary groups 
$\U(\cH)$. Their restrictions to the direct limit groups 
$\U_\infty(\K) := \indlim \U_n(\K)$ are precisely 
the tame representations (\cite[Sect.~3]{Ol90}, \cite{Ol78}). 

\section{Separable representations of infinite dimensional unitary groups} 
\mlabel{app:g}

\begin{thm}{\rm(Kirillov, Olshanski, Pickrell)} \mlabel{thm:8.1b}
For a separable continuous unitary representation 
$\pi$ of the unitary group $\U(\cH)$ of a separable 
infinite dimensional Hilbert space $\cH$ over $\K \in \{\R,\C,\H\}$, 
the following assertions hold: 
\begin{description}
\item[\rm(a)] $\pi$ is also continuous with respect to the strong 
operator topology. 
\item[\rm(b)] $\pi$ is a direct sum of irreducible representations.  
\item[\rm(c)] All irreducible representations occur 
in a tensor product $\cH_\C^{\otimes n}$ 
for some $n \in \N_0$. In particular they are bounded. 
\item[\rm(d)] $\pi$ extends uniquely to a continuous 
representation of the group $\U(\cH)^\sharp$ with the same commutant, 
where 
\[ \U(\cH)^\sharp \cong  
\begin{cases}
\U(\cH_\C) & \text{ for } \K = \R,  \\ 
\U(\cH)^2 & \text{ for } \K = \C,  \\ 
\U(\cH^\C) & \text{ for } \K = \H, 
\end{cases} \] 
is the group with Lie algebra 
$\L(\U(\cH)^\sharp) = \fu(\cH) + i \Herm(\cH)$. 
Here $\cH^\C$ denotes the complex Hilbert space underlying the quaternionic 
Hilbert space $\cH$, and, for $\K = \C$, the inclusion 
$\eta \: \U(\cH) \to \U(\cH)^2$ is the diagonal embedding. 
\end{description}
\end{thm}

\begin{prf} (a) is due to Pickrell (\cite{Pi88}). 

(b) and (c) are claimed by Kirillov in \cite{Ki73} but detailed 
proofs have been provided by Olshanski in \cite{Ol78} 
(see also \cite{Ol90}). Note that, for $\K = \C$, 
$\cH_\C \cong \cH \oplus \cH^*$, so that one could also say that 
the irreducible representations occur in some
$\cH^{\otimes n} \otimes (\cH^*)^{\otimes m}$. 

(d) Here the main idea is that (a) and (c) implies that 
the representation extends to a representation of the full 
contraction semigroup $C(\cH)$ in which $\U(\cH)$ is strongly 
dense. Then one applies a holomorphic extension 
argument which yields a representation of the 
$c$-dual group of $\GL(\cH)$, and the compatibility with the 
embedding $\eta \: \U(\cH) \to\U(\cH)^\sharp$ 
shows that we actually obtain a representation 
of $\U(\cH)^\sharp$. 
\end{prf}

\begin{cor} \mlabel{cor:e.2} Let $K$ be a quotient of a product 
$K_1 \times \cdots \times K_n$, where each $K_j$ 
is compact or a quotient of 
some group $\U(\cH)$, where $\cH$ is a separable $\K$-Hilbert space.  
Then every separable continuous unitary representation 
$\pi$ of $K$ is a direct sum of irreducible representations 
which are bounded. 
\end{cor}

The preceding corollary means that the separable representation 
theory of $K$ very much resembles the representation theory of a compact 
group.

\section{Perspectives and open problems} 
\mlabel{app:h}

\subsection{Positive energy without semiboundedness} 

\begin{prob} To clarify the precise relation between semiboundedness 
and the positive energy condition, one has to answer the following question: 
Is every irreducible  positive energy representation 
of a full hermitian Lie group $G$ semibounded?
From Theorem~\ref{thm:posen}(iv) we know that the converse is true.  

Here is what we know: 
Let $(\pi, \cH)$ be an irreducible positive energy representation 
of the hermitian Lie group $(G,\theta, d)$. 
From \cite[Prop.~3.6]{Ne10d} we derive that 
$\cH$ is generated by the closed subspace 
$V := \oline{(\cH^\infty)^{\fp^-}}$. Now the problem is to 
show that the $K$-representation $\rho$ on $V$ is bounded because 
then Theorem~\ref{thm:3.15} applies. From 
Lemma~\ref{lem:6.2b} we know that 
$\cH$ decomposes into eigenspaces of $\dd\pi(d)$. 

If we assume, in addition, that $\cH$ and hence also $V$ are 
separable, then Corollary~\ref{cor:e.2} applies to the $K$-representation 
on $V$, so that we obtain a bounded irreducible 
$K$-representation $W \subeq V$, but a priori we do not know 
if $W$ can be chosen in such a way that 
$W \cap \cH^\infty$ is dense in $W$. Hence we cannot 
apply the tools from \cite{Ne10d}. 
What is missing at this point are tools to decompose the 
$K$-representation on the Fr\'echet space $\cH^\infty$ 
(cf.\ \cite{Ne10a}) that would lead to the 
existence of an irreducible subrepresentation. 
In a certain sense we are asking for a Peter--Weyl theory 
for groups such as $\U(\cH)$. 
\end{prob}

\begin{prob} Classify semibounded representations of 
Fr\'echet--Lie groups such as $\Diff(\bS^1)$, the Virasoro group  
and affine Kac--Moody groups. 
This requires in particular to extend the 
tools developed in \cite{Ne10d} to strongly 
continuous automorphism groups $(\alpha_t)_{t \in \R}$ 
on Banach--Lie groups, where the infinitesimal generator 
$\alpha'(0)$ is unbounded. One also has to develop suitable 
direct integral techniques which permit to identify smooth vectors. 
\end{prob}

\begin{prob} Use the L\"ucher-Mack Theorem \cite{MN11} 
to determine which irreducible representations we obtain by 
restriction from hermitian groups to 
automorphism groups of real forms $M_\R$ of hermitian symmetric 
spaces~$M$. 
It seems that many of the irreducible semibounded 
representations remain irreducible when restricted to such subgroups. 

On the level of the group $K$, typical examples 
where restrictions to rather small subgroups remain irreducible 
arise in the situation of Theorem~\ref{thm:8.1b}(d). 
In this context all highest weight representation 
of the group $\U(\cH)^\sharp$ with non-negative highest weight 
remain irreducible when restricted to the subgroup 
$\U(\cH)$. 
\end{prob}

\begin{prob} It is also of some interest to identify the 
irreducible semibounded representations for 
hermitian groups $(G,\theta, d)$ which are not irreducible. 
Here the main issue is to understand the situation where 
$\fp$ contains infinitely many simple $JH^*$-ideals. 
\end{prob}

\subsection{Other restricted hermitian groups} 

\begin{ex}\mlabel{ex:h.1}
For smooth manifolds $M$ with $\dim M > 1$, 
the group $C^\infty(M,K)$, $K$ a compact Lie group, 
has natural homomorphisms into 
the group 
\[ G := \U_{\rm res,p}(\cH) := \Big\{ g \in \U(\cH) \: 
g_{12} \in B_p(\cH_-,\cH_+), 
g_{21} \in B_p(\cH_+,\cH_-)\Big\}\]  
for some $p > 2$ (cf.\ \cite{Mi87, Mi89}, \cite{Pi89}). 
The Lie algebra of this group contains the element 
$d := \frac{i}{2}\diag(\1,-\1)$ 
for which $\theta := e^{\pi \ad d}$ defines an involution and we can also 
ask for its semibounded unitary representations. 
Note that $(G,\theta,d)$ is not a hermitian Lie group 
in the sense of Definition~\ref{def:1.1} because 
$\fp \cong B_p(\cH_-, \cH_+)$ carries for $p > 2$ no 
Hilbert space structure invariant under 
the adjoint action of the group 
$K = G^\theta \cong\U(\cH_+) \times \U(\cH_-)$. 

First we claim that $H^2_c(\g,\R) = \{0\}$, i.e., that all 
central extensions of $\g$ are trivial. 
In fact, with the same arguments as in the proof of 
Lemma~\ref{lem:cocyc}, we see 
that every continuous cocycle $\omega \: \g \times \g \to \R$ 
is equivalent to one vanishing on $\fk \times \g$. Then 
$\omega$ defines on $\fp \times \fp$ an $\Ad(K)$-invariant 
skew-symmetric bilinear form, and Lemma~\ref{lem:cocyc} 
implies that its restriction to the dense subspace 
$B_2(\cH_-,\cH_+)$ is a multiplies 
of the canonical cocycle $\Im \tr(xy^*)$. Since this cocycle 
is not continuous with respect to $\|\cdot\|_p$, we obtain 
$\omega =0$, and therefore $H^2_c(\g,\R) = \{0\}$. 

Next we observe that every semibounded unitary 
representation $\pi$ of $G$ restricts to a semibounded 
representation of the subgroup $\U_{\rm res}(\cH)$, whose 
central charge $c$ vanishes (Theorem~\ref{thm:8.1}). 
If $\dd\pi$ is non-zero on $\fp$, then 
Corollary~\ref{cor:e.2} further implies that 
$\pi$ extends to a bounded highest weight representation 
$(\pi_\lambda, \cH_\lambda)$ of the full unitary 
group $\U(\cH)$. Therefore the group $\U_{\rm res,p}(\cH)$ has 
no unbounded semibounded irreducible unitary representations. 
If $\cH$ is separable, then Pickrell show in 
\cite[Prop.~7.1]{Pi90} that all separable continuous unitary 
representations of $\U_{\rm res,p}(\cH)$ extend uniquely 
to strongly continuous representations of the full unitary group 
$\U(\cH)$, hence are direct sums of bounded representation. 
\end{ex}

\begin{prob} \mlabel{ex:h.2} For the Schatten norms, we have 
the estimate 
\[  \|AB\|_p \leq \|A\|_{p_1} \|B\|_{p_2} \quad 
\hbox{ for } \quad 
\frac{1}{p} \leq \frac{1}{p_1} + \frac{1}{p_2} \] 
(\cite[Th.~IV.11.2]{GGK00}). 
This leads in particular to 
\[  \|AB\|_{p/2} \leq \|A\|_{p} \|B\|_{p}.\] 
For each $p \geq 2$ we thus obtain for 
$\cH = \cH_+ \oplus \cH_-$ a Banach algebra 
\begin{equation}
  \label{eq:bp}
B_{p/2,p}(\cH) 
:= \Big\{ x  = \pmat{ x_{11} & x_{12} \\ x_{21} & x_{22}} 
\in B(\cH) \: 
\|x_{11}\|_{p/2}, \|x_{12}\|_p, \|x_{21}\|_p, 
\|x_{22}\|_{p/2} < \infty \Big\} 
\end{equation}
and a corresponding unitary group 
\begin{equation}
  \label{eq:up} 
\U_{p/2,p}(\cH) := \U(\cH) \cap (\1 + B_{p/2,p}(\cH)).
\end{equation}

Accordingly, we obtain for each $p \geq 2$ a variant 
of a hermitian Lie algebra by 
\[ \g := \R d + \fu_{p/2,p}(\cH) 
\cong \fu_{p/2,p}(\cH) \rtimes_{\ad d} \R \quad \mbox{ for } \quad 
d = \frac{i}{2}\diag(\1,-\1)\]  
and a corresponding group $G := \U_{p/2,p}(\cH) \rtimes \R$. 
Then 
\[ \fk = \fz_\g(d) \cong \fu_{p/2}(\cH_+) \oplus \fu_{p/2}(\cH_-)\oplus \R.\] 
With a similar argument as in Lemma~\ref{lem:cocyc}, we 
see that each cohomology class in $H^2_c(\g,\R)$ contains a 
$d$-invariant cocycle~$\omega$. 
Then $\omega(\fk,\fp)= \{0\}$, and a similar argument as 
in Example~\ref{ex:h.1} implies that it vanishes 
on $\fp \times \fp$. However, for $p > 2$, the Lie algebra 
$\fu_{p/2}(\cH_+)$  has non-trivial cocycles which can be written as 
$\omega_D(x,y) := \tr([x,y]D)$ for $D \in B(\cH_+)$ if $p \leq 4$, and 
for $D \in B_q(\cH_+)$ if $p > 4$ and $\frac{1}{q} + \frac{4}{p} = 1$ 
(\cite[Prop.~III.19]{Ne03}). 

Suppose that $p > 2$. One can show that  the Lie algebra $\fk$ also has the 
property that $\fk/\z(\fk)$ contains no open invariant cones. 
In particular, all irreducible semibounded representations of 
$G$ are holomorphically induced from bounded representations of 
$K = (G^\theta)_0$ (Theorem~\ref{thm:6.2b}). 
For $p > 4$, one can apply the results 
from Appendix~\ref{app:e} to the groups 
$\U_{p/4}(\cH_\pm)$ to show that, for all bounded 
unitary representations of central extensions of 
the groups $\U_{p/2}(\cH_\pm)$, the center acts trivially, 
which leaves only representations extending to highest weight 
representations of the full unitary groups $\U(\cH_\pm)$. 
From that one concludes that the central extensions of 
$\U_{p/2,p}(\cH)$ do not lead to new semibounded representations, 
and that all semibounded representations extend to the 
larger group $\U_{\rm res,p}(\cH)$, hence even to 
$\U(\cH)$ (cf.~Example~\ref{ex:h.1}). 

For $2 < p \leq 4$ the preceding method does not work. It is an 
interesting problem whether in this case the central extensions 
of $\U_{p/2}(\cH_\pm)$ have bounded irreducible representations 
which are non-trivial on the center. We do not expect that this 
is the case. 

We also note that, since $\U_{p/2}(\cH)$ is contained 
in $\U_{p/2,p}(\cH)$, the structure of the bounded representations of 
$G$ follows from \cite[Thm.~III.14]{Ne98}. They all extend to 
highest weight representations of $\U(\cH)$ (Definition~\ref{def:tensrep}). 

For $p = 2$, the  group $\U_{p/2}(\cH_\pm) = \U_1(\cH_\pm)$ 
has no non-trivial central 
extensions, but it has a bounded representation theory 
which is not of type I (cf.\ \cite{Ne98}). Furthermore 
$\fk/\fz(\fk)$ contains non-trivial open invariant cones. 
Therefore one can expect that in this case $G$ has a rich but also 
more complicated variety of semibounded unitary representations. 
\end{prob}

\end{document}